\documentclass[12pt]{article}
\usepackage{amssymb,amsmath,amsfonts,relsize}
\usepackage{graphicx}
\usepackage[labelsep=period]{caption}

\newtheorem{theorem}{Theorem}[section]
\newtheorem{definition}[theorem]{Definition}
\newtheorem{corollary}[theorem]{Corollary}
\newtheorem{proposition}[theorem]{Proposition}
\newtheorem{remark}[theorem]{Remark}
\newtheorem{lemma}[theorem]{Lemma}

\newtheorem{problem}[theorem]{Problem}
\newtheorem{statement}[theorem]{Statement}
\newtheorem{assumption}[theorem]{Assumption}
\newcommand {\Kc}      {{\mathcal K}}
\newcommand {\Gc}      {{\mathcal G}}
\newcommand {\Lc}      {{\mathcal L}}

\newcommand {\Wc}      {{\mathcal W}}
\newcommand {\Hc}      {{\mathcal H}}
\newcommand {\Ac}      {{\mathcal A}}
\newcommand {\Tc}      {{\mathcal T}}

\newcommand {\Cc}      {{\mathcal C}}
\newcommand {\Dc}      {{\mathcal D}}
\newcommand {\Pc}      {{\mathcal P}}
\newcommand {\Ec}      {{\mathcal E}}
\newcommand {\Qc}      {{\mathcal Q}}

\newcommand {\Ic}      {{\mathcal I}}
\newcommand {\Sc}      {{\mathcal S}}
\newcommand {\Fc}      {{\mathcal F}}

\newcommand {\Nc}      {{\mathcal N}}
\newcommand {\LPO}     {L_p(\Omega)}
\newcommand {\WMP}     {W^m_p(\RN)}
\newcommand {\LMP}     {L^m_p(\RN)}

\newcommand {\LMT}     {L^m_p(\RT)}
\newcommand {\WMPO}    {W^m_p(\Omega)}
\newcommand {\LMPO}    {L^m_p(\Omega)}
\newcommand {\LMPW}     {L^m_{\tp}(\RN)}
\newcommand {\LMPWO}    {L^m_{\tp}(\Omega)}
\newcommand {\W}       {W^m_p}
\newcommand {\LM}      {L^m_p}

\newcommand {\tA}      {\widetilde{A}}
\newcommand {\tB}      {\widetilde{B}}
\newcommand {\tF}      {\widetilde{F}}
\newcommand {\tH}      {\widetilde{H}}
\newcommand {\tlh}     {\widetilde{h}}
\newcommand {\tG}      {\widetilde{G}}
\newcommand {\tK}      {\widetilde{K}}
\newcommand {\tC}      {\widetilde{C}}

\newcommand {\tW}      {\widetilde{W}}
\newcommand {\tw}      {\widetilde{w}}
\newcommand {\tp}      {\widetilde{p}}
\newcommand {\tfi}     {\widetilde{\varphi}}
\newcommand {\thw}     {\bar{h}}
\newcommand {\tz}      {\widetilde{z}}
\newcommand {\ta}      {\widetilde{a}}
\newcommand {\tb}      {\widetilde{b}}
\newcommand {\tdl}     {\widetilde{\delta}}

\newcommand {\tgm}     {\widetilde{\gamma}}
\newcommand {\tal}     {\widetilde{\alpha}}
\newcommand {\TV}      {\widetilde{\Tc}}
\newcommand {\NW}      {\widetilde{\Nc}}
\newcommand {\GW}      {\widetilde{\Gc}}

\newcommand {\R}       {{\bf R}}
\newcommand {\N}       {{\bf N}}
\newcommand {\RN}      {\R^n}
\newcommand {\RT}      {\R^2}
\newcommand {\DO}      {\partial\Omega}
\newcommand {\DS}      {\partial \SB}
\newcommand {\DG}      {\partial G}
\newcommand {\DHA}     {\partial H}
\newcommand {\ve}      {\varepsilon}
\newcommand {\emp}     {\emptyset}

\newcommand {\intl}    {\int\limits}
\newcommand {\WPT}     {{\cal WP}_\Omega}
\newcommand {\NPT}     {{\cal NP}_\Omega}
\newcommand {\br}      {\bar{x}}
\newcommand {\sm}      {\setminus}
\newcommand {\SCL}     {\SB^{\cl}}

\newcommand {\KCL}     {K^{\cl}}
\newcommand {\KZC}     {K(z)^{\cl}}
\newcommand {\DKZ}     {\partial K(z)}

\newcommand {\SB}      {\widetilde{S}}
\newcommand {\hS}      {\widehat{S}}
\newcommand {\hF}      {\widehat{F}}
\newcommand {\cs}      {\tilde{c}}
\newcommand {\rs}      {R}
\newcommand {\GB}      {G_B}
\newcommand {\TB}      {\Tc_B}
\newcommand {\TBZ}     {\Tc_{B,z}}
\newcommand {\TBP}     {\TB^\oplus}
\newcommand {\TBM}     {\TB^\ominus}
\newcommand {\ETB}     {e_{\TB}}
\newcommand {\bz}      {\bar{z}}
\newcommand {\bg}      {\bar{\gamma}}
\newcommand {\CE}      {\theta}
\newcommand {\hdl}     {\hat{\delta}}

\newcommand {\rj}      {\rho_{w,j}}
\newcommand {\rp}      {\rho_{w,1}}
\newcommand {\rv}      {\rho_{w,2}}
\newcommand {\pj}      {\varphi_j}
\newcommand {\pp}      {\varphi_1}
\newcommand {\pv}      {\varphi_2}
\newcommand {\FG}      {rapidly }
\newcommand {\x}       {\bar{x}}
\newcommand {\y}       {\bar{y}}
\newcommand {\bu}      {\bar{u}}
\newcommand {\TF}      {\Fc}
\newcommand {\FF}      {\widetilde{\Fc}}
\newcommand {\oi}[2]   {(#1,#2)_{\DS}}
\newcommand {\ci}[2]   {[#1,#2]_{\DS}}
\newcommand {\sol}[2]  {(#1,#2]_{\DS}}
\newcommand {\sor}[2]  {[#1,#2)_{\DS}}
\newcommand {\smed}    {\mathlarger{\sum}}
\newcommand {\sbig}    {\mathlarger{\mathlarger{\sum}}}

\newcommand {\usm}    {\mathsmaller{\bigcup}}
\newcommand {\ism}    {\mathsmaller{\bigcap}}

\newcommand {\diam}    {\operatorname{diam}}
\newcommand {\dist}    {\operatorname{dist}}
\newcommand {\cl}      {\operatorname{\,cl}}
\newcommand {\PR}      {\operatorname{Pr}}
\newcommand {\Lip}     {\operatorname{Lip}}
\newcommand {\len}     {\operatorname{len}}
\newcommand {\length}  {\operatorname{length}}
\newcommand {\VST}     {\vspace*{1mm}}
\newcommand {\bx}      {\hspace{10mm}$\Box$}

\newcommand {\rbx}     {\hspace{10mm}$\vartriangleleft$}
\newcommand {\rbf}     {\hspace{10mm}\vartriangleleft}
\newcommand {\nn}      {\nonumber}
\newcommand {\rf}[1]    {(\ref{#1})}      
\newcommand {\reff}[1] {\ref{#1}}         
\newcommand{\lbl}[1]      {\label{#1}}       
\newcommand{\be}          {\begin{eqnarray}}
\newcommand{\bel}[1]      {\begin{eqnarray} \label{#1}}
\newcommand{\ee}           {\end{eqnarray}}
\newcommand {\SECT}[2] {\section*{\centerline{\normalsize
{\bf #1}}} \setcounter{section}{#2}
\setcounter{theorem}{0}\setcounter{equation}{0}}
\begin{document}
\parindent 1em
\parskip 0mm
\medskip
\centerline{\Large{\bf On planar Sobolev $L^m_p$-extension domains}}
\vspace*{10mm} \centerline{By~
 Pavel Shvartsman}\vspace*{2 mm}
\centerline {\it Department of Mathematics, Technion - Israel Institute of Technology}\vspace*{1 mm}
\centerline{\it 32000 Haifa, Israel}\vspace*{1 mm}
\centerline{\it e-mail: pshv@tx.technion.ac.il}
\vspace*{3 mm}
\centerline{and}
\vspace*{3 mm}
\centerline{Nahum Zobin}\vspace*{2 mm}
\centerline {\it Department of Mathematics, College of William and Mary,}\vspace*{1 mm}
\centerline {\it PO Box 8795, Williamsburg, VA 23187-8795, USA}\vspace*{1 mm}
\centerline{\it e-mail: nxzobi@wm.edu}
\vspace*{6 mm}
\renewcommand{\thefootnote}{ }
\footnotetext[1]{{\it\hspace{-6mm}Math Subject
Classification} 46E35\\
{\it Key Words and Phrases} Sobolev space,
extension operator, domain, subhyperbolic metric.\smallskip
\par Parts of this research were carried out in August 2010 during the Workshop ``Differentiable Structures on Finite Sets'' at the American Institute of Mathematics in Palo Alto, CA,  in September 2012 during the Thematic Program on Whitney Problems at the Fields Institute in Toronto, Canada, and in April 2013 during the Whitney Problems Workshop at the Banff International Research Station, Canada. The authors were generously supported by the AIM, by the Fields Institute and by the BIRS. The second author was also partially supported by the NSF grant DMS 713931.}
\begin{abstract}
For each $m\ge 1$ and $p>2$ we characterize bounded simply connected Sobolev $L^m_p$-extension domains $\Omega\subset \RT$. Our criterion is expressed in terms of certain intrinsic subhyperbolic metrics in $\Omega$.
Its proof is based on a series of results related to the existence of special chains of squares joining given points $x$ and $y$ in $\Omega$.
\par An important geometrical ingredient for obtaining these results is a new ``Square Sepa\-ration Theorem''. It states that under certain natural assumptions on the relative positions of a point $x$ and a square $S\subset\Omega$ there exists a similar square $Q\subset\Omega$ which touches $S$ and has the property that $x$ and $S$ belong to distinct connected components of $\Omega\setminus Q$.
\end{abstract}
\renewcommand{\contentsname}{ }
\tableofcontents
\addtocontents{toc}{{\centerline{\sc{Contents}}}
\vspace*{4mm}\par}
\SECT{1. Introduction}{1}
\addtocontents{toc}{~~~1. Introduction.\hfill \thepage\par\VST}
\indent
\par {\bf 1.1. Main definitions and main results.}
\addtocontents{toc}{~~~~1.1. Main definitions and main results. \hfill \thepage\par}
Let $\Omega$ be an open  subset of $\RN$. We
recall that, given $m\in\N$ and $p\in[1,\infty]$, the homogeneous Sobolev space $\LMPO$ consists of all functions $f\in L_{1,loc}(\Omega)$ whose distributional partial derivatives on $\Omega$ {\it of order $m$} belong to $L_p(\Omega)$. See, e.g., Maz'ya \cite{M}. $\LMPO$ is seminormed by
$$
\|f\|_{\LMPO}:= \smed\{\|D^\alpha
f\|_{L_p(\Omega)}:|\alpha|=m\}.
$$
\par As usual, we let $\WMPO$ denote the corresponding Sobolev space of all functions $f\in \LPO$ whose distributional partial derivatives on $\Omega$ of {\it all orders up to $m$} belong to $L_p(\Omega)$. This space is normed by
$$
\|f\|_{\WMPO}:= \smed\{\|D^\alpha
f\|_{L_p(\Omega)}:|\alpha|\le m\}.
$$
\begin{definition}{\em We say that a domain $\Omega\subset\RN$ has {\it the Sobolev $\LM$-extension property} if there exists a constant $\theta\ge 1$ such that the following condition is satisfied: for every $f\in\LMPO$ there exists a function $F\in\LMP$ such that
$F|_\Omega=f$ and
\bel{IN-FCR}
\|F\|_{\LMP}\le\, \theta\,\|f\|_{\LMPO}\,.
\ee
\par We refer to any domain $\Omega$ which has this property as a {\it Sobolev $L_{p}^{m}$-extension domain}.}
\end{definition}
\medskip
\par Note that in this definition we may omit the requirement of the existence of a constant $\theta$ satisfying inequality \rf{IN-FCR}. (This follows easily from the Banach Inverse Mapping Theorem, see Subsection 7.1). Nevertheless for our purpose it will be convenient to introduce the parameter $\theta$ and the following ``index'' associated with this parameter
\bel{E-MP}
e(\LMPO):=\inf\,\theta
\ee
which provides us with a way of quantifying the Sobolev extension property of $\Omega$.
\par We define Sobolev $\W$-extension domains in an analogous way. (For various equivalent definitions of Sobolev extension domains we refer the reader to Subsection 7.1.)
\par  In this paper we study the following
\begin{problem} {\em Given $p\in[1,\infty]$ and $m\in\N$ find a geometrical characterization of
the class of Sobolev $\LM$-extension domains in $\RN$.\vspace*{2mm}}
\end{problem}
\par We give a complete solution to this problem for the family of bounded simply connected  domains in $\RT$ whenever $p>2$ and $m\in\N$. Our main result is the following
\begin{theorem} \lbl{MAIN-EXT} Let $2<p<\infty$ and let $m\in\N$. Let $\Omega\subset\RT$ be a bounded simply connected domain. Then $\Omega$ is a Sobolev $\LM$-extension domain if and only if for some constant $C>0$ the following condition is satisfied: for every $x,y\in\Omega$ there exists a rectifiable curve $\gamma\subset\Omega$ joining $x$ to $y$ such that
\bel{IN-C}
\intl_{\gamma } \dist(u,\DO)
^{\frac{1}{1-p}}\,ds(u) \le\, C\,\|x-y\|^{\frac{p-2}{p-1}}\,.
\ee
Here $ds$ denotes arc length measure along $\gamma$.
\end{theorem}
\medskip
\par Inequality \rf{IN-C} motivates us to express the statement of Theorem \reff{MAIN-EXT} in terms of certain intrinsic metrics. Following Buckley and Stanoyevitch \cite{BSt3}, given $\alpha \in [0,1]$ and a rectifiable curve $\gamma \subset \Omega $, we define the {\it subhyperbolic length} of $\gamma $ by
\bel{SH-LN}
\len_{\alpha ,\Omega }(\gamma ) :=\intl_{\gamma }
\dist(u,\DO)^{\alpha -1}\,ds(u).
\ee
Then we let $d_{\alpha ,\Omega}$ denote the corresponding {\it subhyperbolic metric} on $\Omega$ given, for each $x,y\in\Omega$, by 
\bel{DEF-D}
d_{\alpha ,\Omega }(x,y):=\inf_{\gamma}\,\len_{\alpha ,\Omega}(\gamma)
\ee
where the infimum is taken over all rectifiable curves
$\gamma \subset \Omega $ joining $x$ to $y$.
\par The metric $d_{\alpha ,\Omega }$ was introduced and studied by Gehring and Martio in \cite{GM}. Note that $\len_{0,\Omega }$ and $d_{0,\Omega }$ are the well-known {\it quasihyperbolic length} and {\it quasihyperbolic distance}, and $\len_{1,\Omega }$ and $d_{1,\Omega}$ are the length of a curve and the
{\it geodesic metric} on $\Omega$ respectively. For various equivalent definitions and other properties of subhyperbolic metrics we refer the reader to \cite{BKos,BSt2,BSt3,BSt,L,S10,S11}. See also Subsection 7.2.
\par   Now inequality \rf{IN-C} can be reformulated in the form
$$
d_{\alpha ,\Omega }(x,y)\le C\,\|x-y\|^\alpha~~~\text{with}~~~\alpha=\tfrac{p-2}{p-1}
$$
which leads us to work with a certain class of domains, essentially those which were introduced in \cite{GM}. See also \cite{BKos,BSt2,BSt3,BSt,L}. In our context here, it seems convenient to use the following terminology which is different from that of \cite{GM} and other papers.
\begin{definition}\lbl{ASHD} {\em For each
$\alpha\in(0,1]$, the domain $\Omega \subset \RN$
is said to be $\alpha$-{\it sub\-hyper\-bolic} if there exists a constant $C_{\alpha,\Omega}>0$ such that for every $x,y\in \Omega$ the following inequality
\bel{CAOM}
d_{\alpha,\Omega}(x,y)\le \,C_{\alpha,\Omega}
\| x-y\|^{\alpha}
\ee
holds.}
\end{definition}
\par For instance, a domain $\Omega$ is a $1$-subhyperbolic if and only if $\Omega$ is a {\it quasiconvex} domain, i.e., if {\it the geodesic metric in $\Omega$ is equivalent to the Euclidean distance.}\smallskip
\par Given an $\alpha$-sub\-hyper\-bolic domain $\Omega \subset \RN$ we define a measure of its subhyperbolicity by letting
\bel{SH-MS}
s_\alpha(\Omega):=\sup_{x,y\in\Omega,\,x\ne y}\,
\frac{d_{\alpha,\Omega}(x,y)}{\| x-y\|^{\alpha}}\,.
\ee
\smallskip
\par Now Theorem \reff{MAIN-EXT} can be reformulated as
follows: {\it For each $p>2$ and each $m\in\N$, a simply connected bounded domain $\Omega\subset\R^2$ is a Sobolev
$L_{p}^{m}$-extension domain if and only if $\Omega$ is a
$\frac{p-2}{p-1}$\,-\,subhyperbolic domain.}\smallskip
\par Actually we prove a slightly stronger version of this result which reveals a universal quantitative connection between Sobolev extension properties of a simply connected bounded domains and their interior subhyperbolic geometry.
\begin{theorem} \lbl{Q-V} Let $2<p<\infty$ and let $m\in\N$. Let $\Omega\subset\RT$ be a bounded simply connected domain. Then $\Omega$ is a Sobolev $L^m_p$-extension domain if and only if $s_{\alpha}(\Omega)$ is finite. In that case $\Omega$
also satisfies
\bel{C-MC}
\tfrac{1}{C}\,e(\LMPO)\le s_\alpha(\Omega)\le C\, e(\LMPO)^{\frac{3p}{p-1}}~~~\text{where}~~~ \alpha=\tfrac{p-2}{p-1}
\ee
and $C>0$ is a constant depending only on $p$ and $m$.
\end{theorem}
\smallskip
\par An approach which we develop in this paper when combined with certain results which were obtained earlier in \cite{S10} enables us to prove the following interesting self-improvement property of Sobolev extension domains. Its proof can be found in Subsection 7.2.
\begin{theorem} \lbl{OP-END} Let $2<p<\infty$ and let $m\in\N$. Let $\Omega\subset\RT$ be a bounded simply connected domain. Suppose that $\Omega$ is a Sobolev $L^m_p$-extension domain.
\par Then $\Omega$ is a Sobolev $L^k_q$-extension domain for all $q>\tp$ and $k\in \N$ where $\tp\in (2,p)$ is a constant depending only on $m$, $p$ and $\Omega$.
\end{theorem}
\par We refer to this result as an ``open ended property'' of planar Sobolev extension domains.\bigskip
\par {\bf 1.2. Historical remarks.}
\addtocontents{toc}{~~~~1.2. Historical remarks.\hfill \thepage\par}
Before we discuss the main ideas of the proof of Theorem \reff{MAIN-EXT} let us recall something of the history of Sobolev extension domains. It is well known that if $\Omega$ is a Lipschitz domain, i.e., if its boundary $\DO$ is locally the graph of a Lipschitz function, then $\Omega$ is a $W_{p}^{m}$-extension domain for every $p\in[1,\infty]$ and every $m\in\N$ (Calder\'{o}n \cite{C2}, $1<p<\infty$, Stein \cite{St}, $p=1,\infty$). Jones \cite{Jn} introduced a wider class of $(\varepsilon ,\delta )$-domains and proved that every $(\varepsilon ,\delta )$-domain is a Sobolev $W_{p}^{m}$-extension domain in ${\bf R}^{n}$ for every $m\ge 1$ and every $p\ge 1$. Burago and Maz'ya \cite{BM}, \cite{M}, Ch. 6, described extension domains for the space $BV(\RN)$ of functions whose distributional derivatives of the first order are finite Radon measures.
\smallskip
\par Let us list several results related to Theorem \reff{MAIN-EXT}. An analogue of Theorem \reff{MAIN-EXT} for the space $W^1_p(\RT)$ has been earlier noted in the literature, see \cite{S10}. In particular, the necessity part of this result was proved by Buckley and Koskela \cite{BKos}, and the sufficiency part by Shvartsman \cite{S10}.\smallskip
\par For $p=\infty$ inequality \rf{IN-C} is equivalent to the quasiconvexity of the domain $\Omega$. In particular, it can be easily seen that the class of bounded $L^1_\infty$-extension domains coincides with the class of quasiconvex bounded domains. The situation is much more complicated for $m>1$. This case has been studied by Whitney \cite{W3} and Zobin \cite{Zob2} who proved the following:\medskip
\par {\it (i). (Whitney) Let $m\ge 1$ and let $\Omega$ be a bounded quasiconvex domain in $\RN$. Then $\Omega$ is an $L^m_\infty$-extension domain;\smallskip
\par (ii). (Zobin) Every finitely connected bounded planar $L^m_\infty$-extension domain is quasiconvex.}
\medskip
\par Zobin \cite{Zob1} also proved  that for every $m>1$ there exists an infinitely connected bounded planar domain $\Omega_m$ which is an $L^m_\infty$-extension domain but it is not an $L^k_\infty$-extension domain for any $k, 1\le k<m$. In particular, $\Omega_m$ is not an $L^1_\infty$-extension domain, so that it is not quasiconvex.\smallskip
\par The first result related to description of Sobolev extension domains in $\R^2$ for $1<p<\infty$ was obtained by Gol'dstein, Latfullin and Vodop'janov \cite{GLV,GV1,GV2} who proved that a {\it simply connected bounded planar} domain $\Omega$ is a Sobolev $L^{1}_{2}$-extension domain if and only if its boundary is a {\it quasicircle}, i.e., if it is the image of a circle under a quasiconformal mapping of the plane onto itself. See also \cite{GVR}. Jones \cite{Jn} showed that every {\it finitely connected} domain $\Omega\subset\RT$ is a $W^{1}_{2}$-extension domain if and only it if its boundary consists of finite number of points and quasicircles; the latter is equivalent to the fact that $\Omega$ is an $(\ve,\delta)$-domain for some positive  $\ve$ and $\delta$. Christ \cite{Chr} proved that the same result is true for  $W^{2}_{1}$-extension domains.
\par Maz'ya \cite{M} gave an example of a simply connected domain $\Omega\subset\RT$ such that $\Omega$ is a $W^1_p$-extension domain for every $p\in[1,2)$, while $\R^2\setminus \Omega^{\cl}$ is a $W^1_p$-extension domain for all $p>2$. However the boundary of $\Omega$ is not a quasicircle. See also \cite{KYZ}.
\par Koskela, Miranda and Shanmugalingam \cite{KMS} showed that a bounded simply connected  planar domain $\Omega$ is a $BV$-extension domain if and only if {\it the complement of $\,\Omega$
is quasiconvex}. (This result partly relies on the above-mentioned work of Burago and Maz'ya \cite{BM}.)
\par We refer the reader to \cite{Chr,HKT,HKT1,K1,K,M,MP,Zh,Zh1} and  references therein for other results related to Sobolev extension domains and techniques for obtaining them.\bigskip
\par {\bf 1.3. Our approach: ``The Wide Path'' and ``The Narrow Path''.}
\addtocontents{toc}{~~~~1.3. Our approach: ``The Wide Path'' and ``The Narrow Path''.\hfill \thepage\par\VST}
Let us briefly indicate the main ideas of the proof of Theorem \reff{MAIN-EXT}.
\par Shvartsman \cite{S10} proved that $e(\WMPO)\le C(m,p)\, s_\alpha(\Omega)$ provided that $p>n>1$, $\alpha=\frac{p-n}{p-1}$, and $\Omega$ is an arbitrary
{\it locally}\, $\alpha$\,-\,subhyperbolic domain in $\RN$. (The locality means that $\Omega$ satisfies inequality \rf{CAOM} for all $x,y\in\Omega$ such that $\|x-y\|\le \delta$ where $\delta$ is a positive constant depending only on $\alpha$ and $\Omega$.)
\par Trivial changes in the proof of this result (mostly related to omitting calculation of $L_p$-norms of derivatives of order less than $m$) lead
us to a similar statement for the space $\LMP$ which we now formulate.
\begin{theorem} \lbl{SOB-EXT-RN} Let $n<p<\infty$ and let $\Omega\subset\RN$ be an $\alpha$\,-\,subhyperbolic domain where $\alpha=\tfrac{p-n}{p-1}$. Then $\Omega$ is a Sobolev $L_{p}^{m}$-extension domain for every $m\ge 1$. 
\par Furthermore, $e(\LMPO)\le C s_\alpha(\Omega)$ where $C$ is a  constant depending only on $n,m$ and $p$.
\end{theorem}
\par Applying this theorem to an arbitrary bounded simply connected domain $\Omega\subset\RT$ we obtain {\it the sufficiency part} of Theorem \reff{MAIN-EXT} and the first inequality in \rf{C-MC}.
\medskip
\par We turn to the proof of {\it the necessity part} of Theorem \reff{MAIN-EXT} and the second inequality in \rf{C-MC}. These statements are equivalent to the following 
\begin{theorem} \lbl{MAIN-NEC} Let $2<p<\infty$, $m\in\N$, and let $\alpha=\frac{p-2}{p-1}$. Let $\Omega\subset\RT$ be a bounded simply connected domain. Suppose that there exists a constant $\CE\ge 1$ such that
 every function $f\in\LMPO$ extends to a function $F\in\LMT$ for which $\|F\|_{\LMT}\le \CE\,\|f\|_{\LMPO}$.
\par  Then for every $\x,\y\in\Omega$ the following inequality
\bel{IN-C-NEC}
d_{\alpha,\Omega}(\x,\y)\le \,C\,\|\x-\y\|^{\alpha}
\ee
holds. Here $C=\tC\,\theta^{\frac{3p}{p-1}}$ where $\tC$ is a positive constant depending only on $m$ and $p$.
\end{theorem}
\par Let us describe the main steps of the proof of inequality \rf{IN-C-NEC}. Let $\Omega$ be a domain satisfying the hypothesis of Theorem \reff{MAIN-NEC}. Suppose that $\x$ and $\y$ are a pair of points in $\Omega$ for which there exists a function $F_m \in \LMPO$ (depending on $\x$ and $\y$) which has the following properties:
\bel{T-A}
D^\beta F_{m}(\x)=0~~~\text{for all}~~~\beta,\, |\beta|=m-1,
\ee
\bel{T-0}
\|F_{m}\|_{\LMPO}\le C_1
\ee
and
\bel{T-1}
d_{\alpha,\Omega}(\x,\y)^{1-\frac1p} \le\, C_2\,\smed_{|\beta|=m-1}
|D^\beta F_{m}(\y)|
\ee
where $C_1$ and $C_2$ are certain positive constants depending only on $m$, $p$ and $\CE$. We shall prove that the existence of such a function $F_m$ implies that
\bel{ID-1}
d_{\alpha,\Omega}(\x,\y)\le \,C\,\|\x-\y\|^{\alpha}~~~\text{with}~~~ \alpha=\tfrac{p-2}{p-1}
\ee
and $C=C(m,p,\CE)$.
\par In fact, since $\Omega$ is an $\LM$-extension domain, the function $F_{m}$ extends to a function $\TF\in\LMT$ with
\bel{TF-4}
\|\TF\|_{\LMT}\le \CE\,\|F_{m}\|_{\LMPO}\le C_1\,\CE\,.
\ee
By the Sobolev-Poincar\'{e} inequality, the partial derivatives of $\TF$ of order $m-1$ satisfy the H\"older condition of order $1-\tfrac{2}{p}$, i.e.,
\bel{SP-IN}
|D^{\beta }\TF(u)-D^{\beta }\TF(v)| \le
C_3\,\|\TF\|_{\LMT}\|u-v\|^{1-\frac{2}{p}}
\ee
for all $\beta$ with $|\beta|=m-1$ and all $u,v\in\RT$. Here $C_3=C_3(m,p)$. See, e.g., \cite{M} or \cite{MP}.
\par By \rf{T-A},
$$
\smed_{|\beta|=m-1}|D^\beta F_{m}(\y)|=\smed_{|\beta|=m-1}
|D^{\beta}\TF(\x)-D^{\beta }\TF(\y)|
$$
so that applying \rf{SP-IN} to $\x$ and $\y$ we obtain
\bel{SM-ED}
\smed_{|\beta|=m-1}|D^\beta F_{m}(\y)|\le C_4\,C_3\,
\|\TF\|_{\LMT}\,\|\x-\y\|^{1-\frac{2}{p}}\le C_4\,C_3\,C_1\,\CE\,
\|\x-\y\|^{1-\frac{2}{p}}.
\ee
Here $C_4=C_4(m)$. Hence, by \rf{T-1},
\bel{DA-OM}
d_{\alpha,\Omega}(\x,\y)^{1-\frac1p} \le\, C_2\,\smed_{|\beta|=m-1}
|D^\beta F_{m}(\y)|\le C_1\,C_2\,C_3\,C_4\,\CE\,
\|\x-\y\|^{1-\frac{2}{p}}
\ee
proving \rf{ID-1}.
\par These observations enable us to reduce the proof of Theorem \reff{MAIN-NEC} to constructing a function $F_{m}=F_{m}(\cdot:\x,\y)\in\LMPO$ satisfying conditions \rf{T-A}, \rf{T-0} and \rf{T-1}. This must be done for each pair of points $\x$ and $\y$ in $\Omega$ (subject of course to the requirement that $\Omega$ satisfies the hypotheses of the theorem). We refer to $F_{m}$ as a {\it   ``\FG growing'' function} associated with the points $\x$ and $\y$.
\par As we have mentioned above, two particular cases of Theorem \reff{MAIN-NEC} were proved earlier by Zobin \cite{Zob1} (for the space $L^m_\infty(\RT)$, $m\in \N$,) and by Buckley and Koskela \cite{BKos} (for the space $L^1_p(\RT)$, $2<p<\infty$). In \cite{Zob1} a construction of the ``\FG growing'' function $F_{m}$ suggested by Zobin relies on the existence of a certain chain of subdomains of $\Omega$, so-called {\it ``rooms''} and {\it ``enfilades''}, which joins $\x$ to $\y$ in $\Omega$. In \cite{BKos} Buckley and Koskela construct the function $F_{m}$ using another approach which involves cutting the domain $\Omega$ into certain disjoint pieces of suitable geometry (so-called {\it ``slices''}). See \cite{Zob1} and \cite{BKos} for the details. These two approaches are very different. We were not able to find a direct and simple generalization of either of them to the case of the Sobolev space $\LM(\Omega)$ for {\it arbitrary} $p>2$ and $m\in\N$.
\par In this paper we suggest a new method for constructing the ``\FG growing'' functions defined on bounded simply connected planar domains. In a similar spirit to \cite{Zob1} and \cite{BKos}, given $\x,\y\in\Omega$ we also construct the function $F_{m}=F_m(\cdot:\x,\y)$ using a special chain of touching subdomains of $\Omega$ joining $\x$ to $\y$. A convenient feature of our construction is that each subdomain of this chain has a very simple geometrical structure - it is an open {\it square} lying in $\Omega$.\smallskip
\par Let us describe our approach in more detail. It is based on the existence of two geometrical objects associated with the points $\x,\y\in\Omega$. We refer to these objects as {\it ``The Wide Path''} and {\it ``The Narrow Path''}. Both ``The Wide Path'' and ``The Narrow Path'' are open subsets of $\Omega$ and they both have a rather simple geometrical structure. More specifically, each of these sets is a {\it chain of open touching subsquares of $\Omega$ joining $\x$ to $\y$}.
\par We describe the geometrical structure of ``The Wide Path'' more precisely in the next theorem. In its formulation and everywhere below the word ``square'' will mean an {\it open} square in $\RT$ whose sides are {\it parallel to the coordinate axes}. By $E^{\cl}$ we denote the closure of a set $E\subset\RT$, and by $E\hspace{0.1mm}^{\circ}$ its interior.
\begin{theorem} \lbl{W-PATH}(``The Wide Path Theorem'') Let $\Omega$ be a simply connected bounded domain in $\RT$, and let $\x,\y\in\Omega$. There exists a finite family
$$
\Sc_\Omega(\x,\y)=\{S_1, S_2,...,S_k\}
$$
of pairwise disjoint squares in $\Omega$ such that \smallskip 
\par (i). $\x\in S_1$ and $\y\in S_k^{\cl}$;\medskip
\par (ii). $S_i^{\cl}\cap S_{i+1}^{\cl}\cap\Omega \ne\emptyset$ for all $i=1,...,k-1$, but
$S_i^{\cl}\cap S_j^{\cl}\cap\Omega=\emp$ for all $1\le i,j\le k$ such that $|i-j|>1$;
\medskip
\par (iii). For every $i=2,...,k-1$ the open set $\Omega \setminus S_i^{\cl}$ is not connected, and the sets
$$\bigcup_{j<i}S_{j}~~~\text{and}~~~\bigcup_{j>i}S_{j}$$ belong to distinct connected components of $\Omega \setminus S_i^{\cl}$.
\end{theorem}
\par This result is the main ingredient of our geometrical construction. We consider the proof of Theorem \reff{W-PATH}, which we present in Sections 2 and 3, to be the most difficult technical part of this paper.
\par It may happen that for certain $i\in\{1,...,k-1\}$ the intersection $S_i^{\cl}\cap S_{i+1}^{\cl}$ is exactly a singleton $\{w_i\}$. In this case we define an additional square $\hS_i$ centered at $\{w_i\}$ of diameter $2\delta$ where $\delta$ is a sufficiently small positive number. See Definition \reff{DF-SH}. We put $\hS_i:=\emp$ whenever $i=k$ or when $S_i^{\cl}\cap S_{i+1}^{\cl}$ is not a singleton and $1\le i<k$.
\par Let
\bel{WP-DEF}
\WPT^{(\x,\y)}:=\left(\,\bigcup_{i=1}^{k} \left(S_i^{\cl}\,\usm\, \hS_i\,\right)\right)^{\circ}\,.
\ee
We refer to the open set $\WPT^{(\x,\y)}$ as a ``Wide Path'' joining $\x$ to $\y$ in $\Omega$.\medskip
\par See Figure 1 for an example of a domain $\Omega$, points $\x,\y\in\Omega$ and a ``Wide Path'' joining $\x$ to $\y$ in $\Omega$ which consists of twelve  consecutively touching squares $S_i, i=1,...,12$.
\begin{figure}[h!]
\center{\includegraphics[scale=0.9]{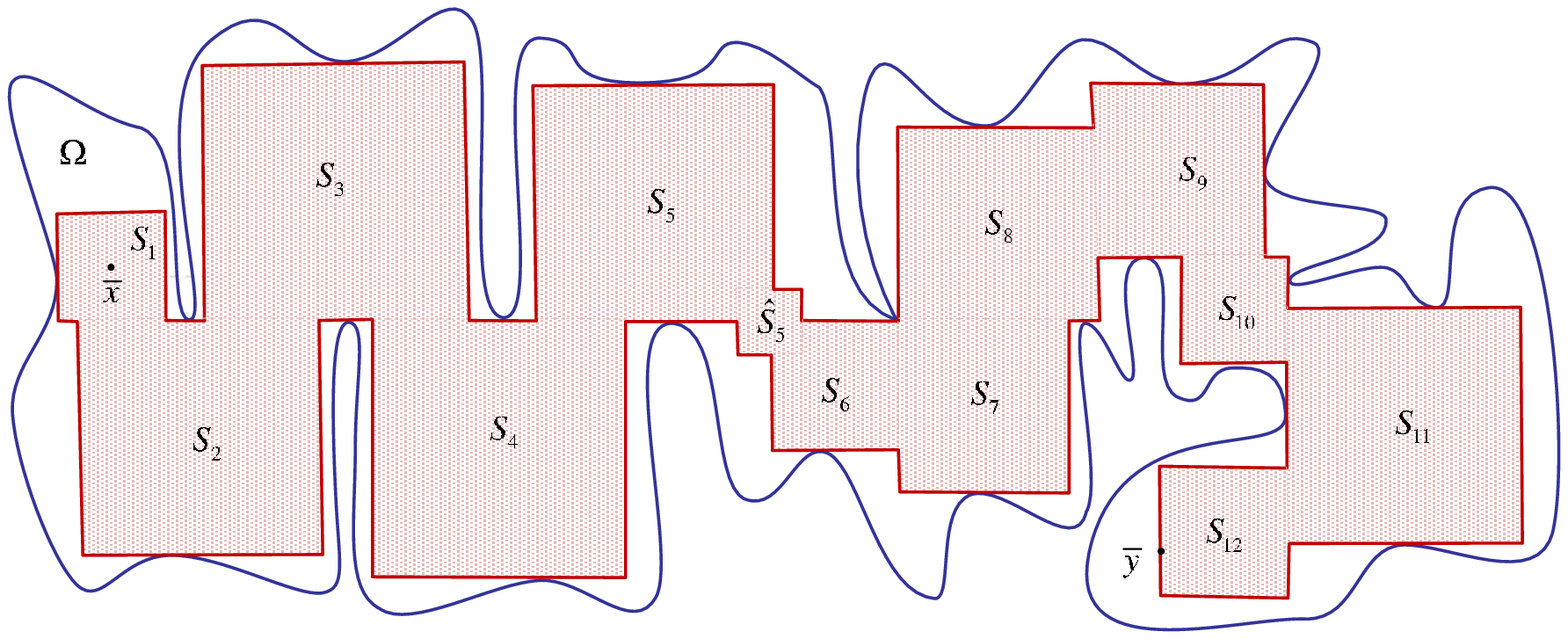}}
\caption{An example of a ``Wide Path'' joining $\x$ to $\y$ in $\Omega$.}
\end{figure}
\par The set $\WPT^{(\x,\y)}$ is an open subset of $\Omega$ possessing a number of pleasant properties which we present and prove in Sections 3. In Section 4 we study Sobolev extension properties of ``The Wide Path''. The following extension theorem is the main result of that section.
\begin{theorem} \lbl{WP-EXT} Let $p>2$ and $m\in\N$. Let $\x,\y\in\Omega$ where $\Omega$ is a simply connected bounded domain in $\RT$. If $\,\Omega$ is a Sobolev $L^m_p$-extension domain, then any ``Wide Path'' $\Wc=\WPT^{(\x,\y)}$ joining $\x$ to $\y$ in $\Omega$ has the Sobolev $L^m_p$-extension property.
\par Furthermore,
\bel{AD}
e(\LM(\Wc))\le C \,e(\LMPO)
\ee
where $C$ is a constant depending only on $m$ and $p$.
\end{theorem}
\smallskip
(See \rf{E-MP} for the definition of the indices appearing in \rf{AD}.)
\bigskip
\par Our next step is to construct ``The Narrow Path''. More specifically, in Section 5, given any ``Wide Path''  $\Wc=\WPT^{(\x,\y)}$ generated from a family $\{S_1, S_2,...,S_k\}$ of squares, we prove  the existence of a family $\Qc_\Omega(\x,\y)=\{Q_1, Q_2,...,Q_k\}$ of pairwise disjoint squares having several ``nice'' properties. Let us list some of them:
\smallskip
\par (i). $Q_1=S_1$, $Q_k=S_k$, and $Q_i\subset S_i$, $1\le i\le k$;\smallskip
\par (ii). $Q_i^{\cl}\cap Q_{i+1}^{\cl}\ne\emp$, $1\le i\le k-1$; \smallskip
\par (iii). $\diam Q_{i+1}\le 2\dist(Q_{i},Q_{i+2})$ provided $Q_{i}^{\cl}\cap S_{i+2}^{\cl}=\emp$ and $1\le i\le k-2$.\medskip
\par For additional properties of the family  $\Qc_\Omega(\x,\y)$ we refer the reader to Proposition \reff{3-SQ}.
\medskip
\par Let
\bel{NP}
\NPT^{(\x,\y)}:=\left(\,\bigcup_{i=1}^{k} \left(Q_i^{\cl}\,\usm\, \hS_i\,\right)\right)^{\circ}\,.
\ee
We refer to the open set $\NPT^{(\x,\y)}$ as a ``Narrow Path'' joining $\x$ to $\y$ in $\Omega$.\medskip
\par\noindent Figure 2 shows a ``Narrow Path'' corresponding to ``The Wide Path'' shown in Figure 1.
\begin{figure}[h]
\center{\includegraphics[scale=0.87]{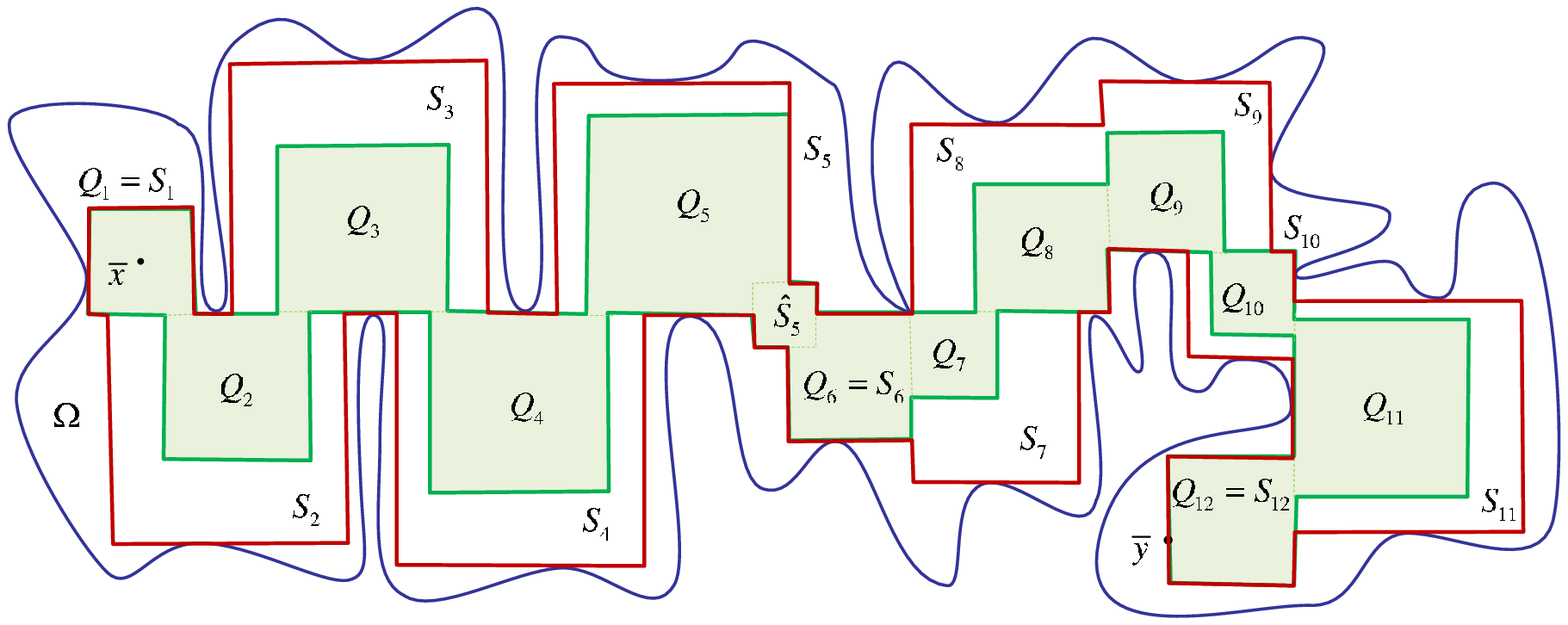}}
\caption{A ``Narrow Path'' joining $\x$ to $\y$ in $\Omega$.}
\end{figure}
\par ``The Narrow Path'' $\Nc=\NPT^{(\x,\y)}$ has a simpler geometrical structure than ``The Wide Path'' $\Wc=\WPT^{(\x,\y)}$. Furthermore its  extension properties are similar to those of $\WPT^{(\x,\y)}$. In particular Theorem \reff{NP-EXT}, which is proved in Section 5, states that {\it every function $f\in\LM(\Nc)$
extends to a function $F\in\LMPO$ such that}
\bel{NP-K}
\|F\|_{\LMPO}\le C(m,p)\,\CE^2\|f\|_{\LM(\Nc)}
\ee
provided $\Omega$ satisfies the hypothesis of Theorem \reff{MAIN-NEC}.
\bigskip
\par In Section 6 we construct the ``\FG growing'' function $F_{m}$. We do this in two steps. In the first step  {\it we define a function $h_{m}$ on ``The Narrow Path''} $\Nc=\NPT^{(\x,\y)}$ (see Definition \reff{FGR-F}).
We prove that
\bel{HS-A}
D^\beta h_{m}(\x)=0,~~~ |\beta|=m-1,
\ee
\bel{HS-B}
\|h_{m}\|_{\LM(\Nc)}^p\le C\smed_{|\beta|=m-1}|D^\beta h_{m}(\y)|
~~~\text{and}~~~
d_{\alpha,\Omega}(\x,\y) \le\, C\smed_{|\beta|=m-1}
|D^\beta h_{m}(\y)|
\ee
where $C$ is a constant depending only on $m$ and $p$.
(See Proposition \reff{P-HSA}.)
\par In the second step of this procedure, using Theorem \reff{NP-EXT}, we extend $h_{m}$ to a function $H_{m}\in\LMPO$ such that
$$
\|H_{m}\|_{\LMPO}\le C(m,p,\CE)\,\|h_{m}\|_{\LM(\Nc)}.
$$
(See inequality \rf{NP-K}.) In Proposition \reff{P-HS} we prove that properties similar to \rf{HS-A} and \rf{HS-B} also hold for the function $H_m$. (See \rf{T-AH}, \rf{T-0H} and \rf{T-1H}.)
\par Finally, we define the function $F_m$ by
$$
F_m(u:\x,\y):=\left(\smed_{|\beta|=m-1}
|D^\beta H_m(\y)|\right)^{-\frac{1}{p}}
\,\cdot\,H_m(u:\x,\y),~~~u\in\Omega\,.
$$
\par It can be readily seen that the above-mentioned properties of $H_{m}$ imply \rf{T-A}, \rf{T-0} and \rf{T-1} proving that $F_{m}$ is a ``\FG growing'' function associated with $\x$ and $\y$.
\par This completes the proof of inequality \rf{IN-C-NEC} and therefore also the necessity part of Theorem \reff{MAIN-EXT}.
\medskip
\par {\bf Acknowledgements.} We are very thankful to
M.\ Cwikel, C.\ Fefferman and V.\ Gol'd\-shtein for useful suggestions and remarks. We are also very grateful to all participants of the ``Whitney Problems Workshops" in Toronto, August 2012, Banff, April 2013, and Williamsburg, August 2014, for stimulating discussions and valuable advice.
\SECT{2. ``The Square Separation Theorem'' in simply connected domains}{2}
\setcounter{equation}{0}
\addtocontents{toc}{2. ``The Square Separation Theorem'' in simply connected domains.\hfill \thepage\par\VST}
\indent\par {\bf 2.1. Notation and auxiliary lemmas.}
\addtocontents{toc}{~~~~2.1. Notation and auxiliary lemmas. \hfill \thepage\par}
Let us fix some additional notation. Throughout the paper $C,C_1,C_2,...$ will be generic positive constants which depend only on $m$ and $p$. These constants can change even in a single string of estimates. The dependence of a constant on certain parameters is expressed, for example, by the notation $C=C(p)$. We write $A\sim B$ if there is a constant $C\ge 1$ such that $A/C\le B\le CA$.
\par As is customary, the word ``domain'' means an {\it open connected subset of $\RT$}. By $\Sc(\RT)$ we denote the family of all open squares in $\RT$ whose sides are parallel to the coordinate axis.
Given a square $S\in \Sc(\RT)$ by $c_S$ we denote its center and by $r_S$ half of its side length. Given $\lambda>0$ we let $\lambda\, S$ denote the dilation of $S$ with respect to its center by a factor of $\lambda$. We let $S(c,r)$ denote the square in $\RT$ centered at $c$ with side length $2r$. We refer to $r=r_S$ as the ``radius'' of the square $S(c,r)$. Thus $S=S(c_S,r_S)$ and $\lambda S=S(c_S,\lambda r_S)$ for every constant $\lambda>0$.
\par We say that squares $S_1$ and $S_2$ are {\it touching} squares 
$$
\text{if}~~~S_1\cap S_2=\emp~~\text{but}~~S_1^{\cl}\cap S_2^{\cl}\ne\emp.
$$
\par We denote the coordinate axes by $Oz_1$ and $Oz_2$. We also refer to the axis $Oz_j$ as the $z_j$-axis, $j=1,2$.
Given $z=(z_1,z_2)\in\RT$ by
\bel{N-U}
\|z\|:=\max \{|z_1|,|z_2|\}
\ee
and by $\|z\|_2:=(|z_1|^2+|z_2|^2)^{\frac12}$ we denote the uniform and the Euclidean norms in $\RT$ respectively.
\par Let $A,B\subset \RT$. We put
$\diam A:=\sup\{\|a-a'\|:~a,a'\in A\}$ and
$$
\dist(A,B):=\inf\{\|a-b\|:~a\in A, b\in B\}.
$$
\par Given $\ve>0$ and a set $A\subset\RT$ by $[A]_\ve$ we denote the $\ve$-neighborhood of $A$:
\bel{EP-N}
[A]_{\ve}:=\{z\in\RT:\dist(z,A)<\ve\}.
\ee
The Lebesgue measure of a measurable set $A\subset \RT$ will be denoted by $\left|A\right|$. By $\# A$ we denote the number of elements of a finite set $A$.
\par Let $t_0=0<t_1<t_2<...<t_m=1$, and let $\Psi:[0,1]\to\RT$ be a continuous mapping which is {\it linear} on every subinterval $[t_i,t_{i+1}]$. We refer to the curve $\gamma=\Psi([0,1])$ as a {\it polygonal curve}. Thus $\gamma$ is the union of a finite number of line segments $[\Psi(t_i),\Psi(t_{i+1})]$, $i=0,..,m-1$. We refer to these line segments as {\it edges}. An endpoint of an edge is called a {\it vertex}.
\par In what follows the word ``path'' will mean a {\it polygonal curve}. We say that a path is {\it simple} if it does not self intersect. We also refer to a simple  {\it closed} path as a {\it simple polygon}.
\par Finally, for each pair of points $z_1$ and $z_2$ in $\RT$ we let $[z_1,z_2]$, $(z_1,z_2)$, $[z_1,z_2)$, $(z_1,z_2]$ denote respectively the closed, open and semi-open line segments joining them.
\par Let us present several auxiliary geometrical results which we use in the sequel. First of them relates to certain properties of squares in $\RT$. Recall that we measure distances in $\RT$ with respect to the uniform norm in $\RT$, see \rf{N-U}.
\begin{lemma}\lbl{SQ-PR} Let $S_1=S(c_1,r_1)$ and  $S_2=S(c_2,r_2)$ be squares in $\RT$. Then:
\medskip
\par (i). $S_1\subset S_2$ if and only if $r_1\le r_2$ and
$\|c_1-c_2\|\le r_2-r_1$;
\medskip
\par (ii). $S_1\cap S_2\ne\emp$ if and only if $\|c_1-c_2\|<r_1+r_2$\,;
\medskip
\par (iii). $S_1$ and $S_2$ are touching squares if and only if $\|c_1-c_2\|=r_1+r_2$. In this case
$S_1^{\cl}\cap S_2^{\cl}=\partial S_1\cap \partial S_2,$
and the set $S_1^{\cl}\cap S_2^{\cl}$ is either a line segment or a point.
Furthermore,
\bel{IN-S}
[c_1,c_2]\cap S_1^{\cl}\cap S_2^{\cl}=\{A\}
\ee
where $A:=\alpha c_1+(1-\alpha)c_2$ with $\alpha:=r_2/(r_1+r_2)$.
\end{lemma}
\par An elementary proof of the lemma we leave to the reader as an easy exercise.
\par The following statement is well known in geometry.
\begin{lemma}\lbl{S-PH} Let $\Omega$ be a domain in $\RT$.
\medskip
\par (i). Every two point in $\Omega$ can be joined by a simple path;
\par (ii). Let $x,y\in\Omega$ and let $\Gamma$ be a path
connecting $x$ to $y$ in $\Omega$. Then there exists a simple path $\gamma\subset\Gamma$ which joins $x$ to $y$.
\end{lemma}
\par We will be also needed certain well known results related to the Jordan curve theorem for polygons and certain properties of simply connected planar domains. We recall these results in the next statements. See, e.g \cite{CR} and \cite{E}.
\begin{statement}\lbl{ST-J}
\par (i). Consider a simple polygon $P$ in the plane. Its complement $\RT\setminus P$ has exactly two connected components.
One of these components is bounded (the interior) and the other is unbounded (the exterior), and the polygon $P$ is the boundary of each component;\smallskip
\par (ii). Let $\Omega$ be a simply connected planar domain. Then the interior of any simple polygon $P\subset\Omega$ lies in $\Omega$.
\end{statement}
\smallskip
\begin{definition}\lbl{S-C} {\em Let $y',y''\in\RT$ and let $P\subset\RT$ be a simple polygon. We say that the line segment $[y',y'']$  {\it strictly crosses} $P$ if $[y',y'']\cap P=\{A\}$ for some $A\in P$, and one of the following conditions is satisfied:  \smallskip
\par (i).  $A$ is not a vertex of $P$;
\smallskip
\par (ii). If $A$ is a common vertex of edges $[z',A]$ and $[A,z'']$ in the polygon $P$, then the straight line $\ell$ passing through $y'$ and $y''$ strictly separates $z'$ and $z''$. (I.e., $z'$ and $z''$ lie in distinct open half-planes generated by $\ell$.)}
\end{definition}
\begin{statement}\lbl{S-CR} Let $y',y''\in\RT$ and let $P\subset\RT$ be a simple polygon. If $[y',y'']$  strictly crosses $P$, then $y'$ and $y''$ lie in distinct connected components of $\RT\setminus P$.
\par In particular, let $\gamma$ be a simple path with ends at points $x$ and $y$. If $\gamma$ crosses $P$ exactly once at a point which is not a vertex of $P$ and not a vertex of $\gamma$, then $x$ and $y$ lie in distinct components of $\RT\setminus P$.
\medskip
\end{statement}
\par We turn to the proof of Theorem \reff{W-PATH}. Its main ingredient is the following statement.
\begin{theorem}\lbl{N-SQ} (``The Square Separation Theorem'') Let $\Omega$ be a simply connected domain in $\RT$. Let $\SB\subset \Omega$ be a square such that
$$
\DS \cap \DO \ne \emp.
$$
\par Let $B\in \Omega\setminus \SB^{\cl}$. Then there exists a square $Q\subset \Omega \setminus \SB^{\cl}$ satisfying the following conditions:\medskip
\par (i). $Q^{\cl}\cap \SB^{\cl}\cap\Omega \ne \emp$;
\smallskip
\par (ii). Either  $B\in Q^{\cl}$ or
\bel{Q-R}
\SB~~\text{and}~~B~~\text{lie in different connected components of}~~\Omega\setminus Q^{\cl}.
\ee
\end{theorem}
\medskip
See Figure 3.
\begin{figure}[h]
\center{\includegraphics[scale=1]{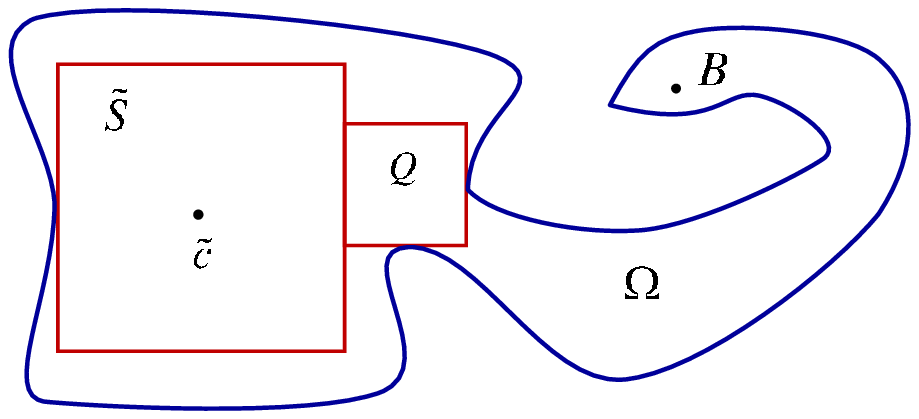}}
\caption{~}
\end{figure}
\par In the sequel we let $\cs$ and $\rs$ denote the center and the ``radius'' of $\SB$ respectively; thus
$$
\SB=S(\cs,\rs).
$$
\par The proof of Theorem \reff{N-SQ} relies on a series of auxiliary results. Towards their formulation let us  introduce several definitions and notations.
\begin{definition}\lbl{ORD} {\em Fix a point $w\in\DS\cap \partial \Omega$. By $\prec$ we denote the total ordering on the set $\DS \setminus \{w\}$ induced by the clockwise direction on $\DS$.
\par Given $a,b \in \DS \setminus \{w\}$, $a\ne b$, we define the open interval $\oi{a}{b}$, closed interval $\ci{a}{b}$ and semi-open intervals  $\sol{a}{b}$ and $\sor{a}{b}$ by letting
$$
\oi{a}{b}=\{x\in \DS \setminus \{w\}: a\prec x\prec b, \,x\ne a,b\},
$$
$\ci{a}{b}=\oi{a}{b}\cup\{a,b\}$ and $\sol{a}{b} = \oi{a}{b}\cup \{b\}$, $\sor{a}{b}=\oi{a}{b}\cup \{a\}.$}
\end{definition}
\par In particular, every {\it connected component} $\Tc$ of $\DS\setminus \DO=\DS\cap\Omega$ is an open interval in $\DS\setminus \{w\}$, completely determined by its beginning $b_\Tc$ and its end $e_\Tc$. Thus $b_\Tc,e_\Tc\in \DS\cap\DO$, $b_\Tc \prec e_\Tc$ and $\Tc = \oi{b_\Tc}{e_\Tc}$\,.
\par It is also clear that for every two {\it distinct} connected components $\Tc_0$ and $\Tc_1$ of $\DS\setminus \DO$ either $e_{\Tc_0} \prec b_{\Tc_1},$  or $e_{\Tc_1} \prec b_{\Tc_0}$. We also notice the following important properties of the components $\Tc_0$ and $\Tc_1$:
$$
\oi{e_{\Tc_0}}{b_{\Tc_1}}\cap \DO\ne\emp~~~\text{provided}~~~e_{\Tc_0} \prec b_{\Tc_1}.
$$
\begin{lemma}\lbl{CM-X} (i). Let $G$ be a connected component of $\Omega\setminus \SB^{\cl}$. There exists a unique connected component  $\Tc=\Tc(G)$ of $\DS\setminus\DO$ having the following property:
\bel{P-T}
\text{Every}~x\in G~\text{and every}~y\in\Tc~\text{can be joined by a path}~~\gamma~\text{such that}~\gamma\setminus\{y\}\subset G
\ee
\par (ii). For every connected component $\Tc$ of $\DS\setminus\DO$ there exists a unique  connected component $G$ of $\Omega\setminus \SB^{\cl}$ which satisfies condition \rf{P-T}.
\end{lemma}
\par {\it Proof.} First we prove the following  \bigskip
\par {\it Statement A}: Let $G$ be a connected component of $\Omega\setminus \SB^{\cl}$ and let  $\Tc$ be a connected component  of $\DS\setminus\DO$. Let $x_0\in G$ and let $p_0\in\Tc$. Suppose that
\bel{G0}
\text{there exists a path}~~\gamma_0~~\text{which joins}~~ x_0~~\text{to}~~p_0~~\text{such that}~~\gamma_0\setminus\{p_0\}\subset G.
\ee
\par Then condition \rf{P-T} holds.
\bigskip
\par Let us prove this statement. Since every $x\in G$ can be connected to $x_0$ by a path in $G$, to prove \rf{P-T} it suffices to show that for each $y\in\Tc$ there exists a path $\gamma$ which joins $x_0$ to $y$ such that $\gamma\setminus\{y\}\subset G$.
\par Without loss of generality we can assume that $p_0$ and $y$ belong to the same side of the square $\SB$. In other words, we can assume that $I:=[p_0,y]\subset \Tc$. Since $I$ is a compact subset of $\Omega$, we have $\ve:=\dist(I,\DO)/2>0$.
\par Recall that $[I]_\ve$ denotes the $\ve$-neighborhood of $I$, see \rf{EP-N}. Then, by definition, $[I]_{\ve}\subset \Omega$. Furthermore, the set $D_\ve:=[I]_{\ve}\setminus \SB^{\cl}$ is an open rectangle.
\par Since $\gamma_0$ is a continuous curve which joins $x_0$ to $p_0$, there exists a point $\tp\in \gamma_0\cap [I]_{\ve}$. Since $\gamma_0\setminus\{p_0\}\subset G\subset \Omega\setminus \SB^{\cl}$, we conclude that $\tp\in D_\ve$. Let $\gamma_1:=[\tp,y]$ and let $\gamma_2$ be the union of $\gamma_1$ and the subarc of $\gamma_0$ from $x_0$ to $\tp$.  Since the rectangle $D_{\ve}$ is convex, $\gamma_1\setminus\{y\}\subset D_\ve\subset\Omega$ so that $\gamma_2\setminus\{y\}\subset\Omega$. Since $x_0\in\gamma_2$ we conclude that  $\gamma_2\setminus\{y\}\subset G$ proving Statement A.
\medskip
\par Let us prove part (i) of the lemma. Fix a point $x_0\in G$. By $\gamma_0$ we denote a path in $\Omega$ which connects $x_0$ to the point $\cs$, the center of the square $\SB$. See Lemma \reff{S-PH}.
\par Since $x_0\notin \SB^{\cl}$ and $\cs\in\SB^{\cl}$, there exists a point
$$p_0\in \DS\setminus \DO=\DS\cap\Omega$$ such that  $\gamma_0\setminus\{p_0\}\subset G$. Let $\Tc=\Tc(G)$ be a connected component of $\DS\setminus \DO$ which contains $p_0$. Since condition \rf{G0} is satisfied, by Statement A, condition \rf{P-T} holds.\smallskip
\par Prove {\it the uniqueness} of the component $\Tc=\Tc(G)$. Suppose that the set $\DS\setminus\DO$ contains {\it two distinct connected components}, $\Tc'$ and $\Tc''$, $\Tc'\ne\Tc''$, such that for every $x\in G$ and  every $y'\in\Tc',y''\in\Tc''$ there exist paths $\gamma'$ and $\gamma''$ joining $x$ to $y'$ and $y''$ respectively such that
\bel{P-TU}
\gamma'\setminus\{y'\}\subset G~~~\text{and}~~~\gamma''\setminus\{y''\}\subset G.
\ee
\par Fix a point $\br\in G$ and points $p'\in\Tc'$ and
$p''\in\Tc''$. Without loss of generality we can assume that $p'\prec p''$. Let
$$
V_0:=\oi{p'}{p''}~~~\text{and}~~~V_1:=\DS\setminus \ci{p'}{p''}.
$$
Then $V_0\cup V_1=\DS\setminus \{p',p''\}$.
\par Since $p'\in\Tc'$, $p''\in\Tc''$ and $\Tc'\ne\Tc''$, we have $V_0\nsubseteq\Omega$. In fact, if $V_0\subset\Omega$, then $p'$ and $p''$ belong to the same connected component of $\DS\setminus \DO$ so that $\Tc'=\Tc''$, a contradiction. In the same way we prove that $V_1\nsubseteq \Omega$.
\par Thus there exist points $y_0\in V_0\setminus\Omega$ and $y_1\in V_1\setminus\Omega$. Prove that the existence of these points leads us to a contradiction. By \rf{P-TU}, there exist paths $\Gamma'$ and $\Gamma''$ which connects $\br$ to $p'$ and $p''$ respectively, and such that the sets $\Gamma'\setminus\{p'\}$ and $\Gamma''\setminus\{p''\}$ lie in $G$. Hence  $\Gamma:=\Gamma'\cup \Gamma''$ is a path which joins $p'$ to $p''$ such that $\Gamma\setminus \{p',p''\}\subset G$.
\par By part (ii) of Lemma \reff{S-PH}, there exists a {\it simple} path $\gamma_1\subset \Gamma$ which connects $p'$ to $p''$. Hence, $\gamma_1\setminus\{p',p''\}\subset G$ so that $\gamma_1\setminus\{p',p''\}\subset \RT\setminus \SB^{\cl}$.
\par Let $\gamma_2:=[p',\cs]$ and let $\gamma_3:=[\cs,p'']$. Then the loop
$\tgm:=\gamma_1\cup\gamma_2\cup\gamma_3$
is a simple closed path in $\Omega$, i.e., $\tgm$ is a simple  polygon. See Figure 4.
\bigskip
\begin{figure}[h!]
\center{\includegraphics[scale=0.8]{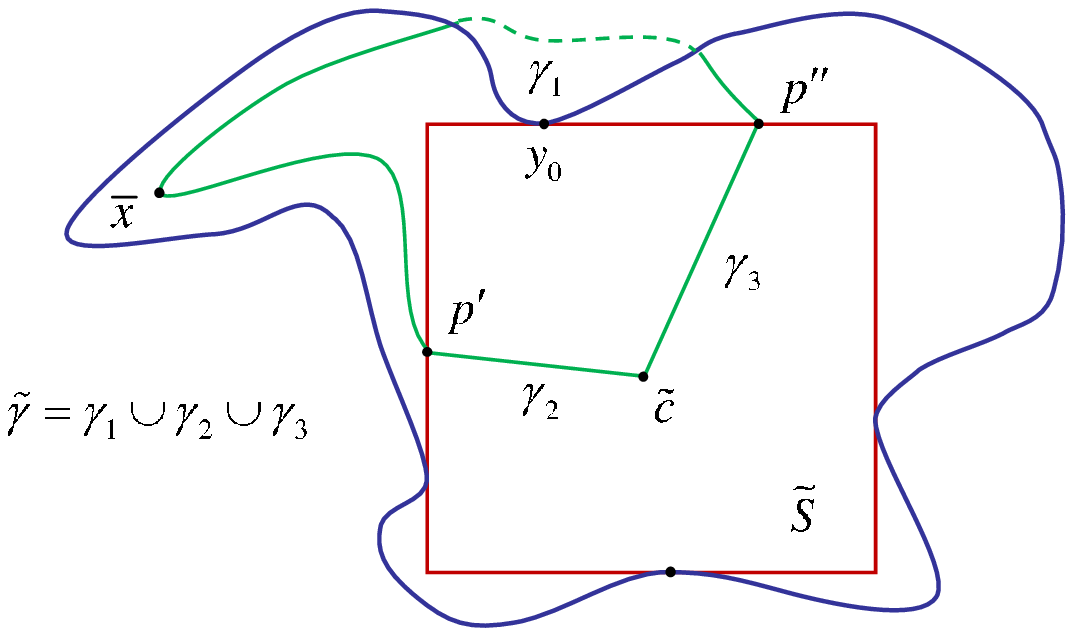}}
\caption{$p'\in\Tc',p''\in\Tc''$ and $\gamma_1$ joins $p'$ to $p''$ in $\RT\setminus\SB^{\cl}$.}
\end{figure}
\medskip
\par By the Jordan curve theorem, see part (i) of Statement \reff{ST-J}, the complement of $\tgm$, the set $\RT\setminus\tgm$, consists of exactly two connected components - the interior component (which is a bounded set), and the exterior component (which is an unbounded set).
We denote these components by $D_{int}$ and $D_{ext}$ respectively. The polygon $\tgm$ is the boundary of these domains, i.e.,
$$
\tgm=\partial D_{int}=\partial D_{ext}.
$$
Furthermore, since $\Omega$ is a simply connected domain and $\tgm\subset\Omega$ is a simple polygon, by part (ii) of  Statement \reff{ST-J},
\bel{D-INO}
D_{int}\subset\Omega.
\ee
\par Clearly, there exists a polygonal path $\gamma'$ (with at most two edges) which joins $y_0$ to $y_1$ in $\SB$ and crosses $(p',\cs]\cup[\cs,p'')$ exactly once at a point which is not $\cs$ or a vertex of $\gamma'$. Since
$$
\gamma_1\setminus\{p',p''\}\subset \RT\setminus \SB^{\cl},
$$
the path $\gamma'$ has no common points with $\gamma_1$, so that $\gamma'$ crosses the simple polygon
$$
\tgm:=\gamma_1\cup\gamma_2\cup\gamma_3
$$
exactly once at a point which is not a vertex of $\tgm$ or $\gamma'$. Hence, by Statement \reff{S-CR}, the points  $y_0$ and $y_1$ lie in different components of $\RT\setminus \tgm$.
\medskip
\par Thus the component $D_{int}$ contains either $y_0$ or $y_1$. But $y_0,y_1\in\RT\setminus\Omega$ so that $D_{int}\nsubseteq\Omega$. On the other hand, by \rf{D-INO}, $D_{int}\subset \Omega$. We have obtained a contradiction which proves part (i) of the lemma. \smallskip
\par Prove (ii). Let $\Tc$ be a connected component of $\DS\setminus \DO$ and let $p_0\in\Tc$. Since the point $p_0\in\DS\setminus \DO=\DS\cap\Omega$, for an $\ve>0$ small enough the square $S(p_0,\ve)\subset\Omega$. Clearly, 
$$
J_{\ve}:=S(p_0,\ve)\setminus \SB^{\cl}
$$
is a non-empty connected set. Also there exists a point $x_0\in J_{\ve}$ such that the line segment $[x_0,p_0)\subset J_{\ve}$.
\smallskip
\par Let $G$ be a connected component of $\Omega\setminus \SB^{\cl}$ which contains $x_0$, and let $\gamma_0:=[x_0,p_0]$. Since $J_{\ve}$ is a connected subset of $\Omega\setminus \SB^{\cl}$ containing $x_0$, we have $J_{\ve}\subset G$. Hence $\gamma_0\setminus\{p_0\}\subset G$ so that condition \rf{G0} is satisfied. On the other hand, by Statement A,  condition \rf{G0} implies \rf{P-T} proving the existence of a connected component $G$ satisfying part (ii) of the lemma.
\par This proof also enables us to show the uniqueness of the component $G$. In fact, let $G'$ be a connected component of $\Omega\setminus \SB^{\cl}$ such that any $x\in G'$ and any $y\in \Tc$ can be joined by a path $\gamma$ with $\gamma\setminus\{y\}\subset G'$. Let $\gamma$ be such a path which connects $x\in G'$ to $y=p_0$. Since $\gamma$ is a continuous curve, there exists a point $z\in\gamma\cap S(p_0,\ve)$. But $\gamma\subset\Omega\setminus \SB^{\cl}$ so that
$$
z\in S(p_0,\ve)\setminus \SB^{\cl}=J_{\ve}.
$$
\par The proof of the lemma is complete.\bx
\medskip
\par Lemma \reff{CM-X} shows that $\Tc=\Tc(G)$ is a one-to-one mapping between the families of connected components of $\Omega\sm \SB^{\cl}$ and the families of connected components of $\DS\sm\DO$.
\par We let $\Gc$ denote the mapping which is inverse to $\Tc(G)$.
Thus for every connected component $\Tc$ of $\DS\sm\DO$ the set $G=\Gc(\Tc)$ is the (unique) connected component of $\Omega\sm \SB^{\cl}$ such that \rf{P-T} is satisfied.
\par We also notice a simple connection between  $\Tc$ and $G=\Gc(\Tc)$:
$$
\Tc(G)=\DG\sm\DO=\DG\cap\Omega.\smallskip
$$
\par We turn to the next step of the proof of Theorem \reff{N-SQ}.\bigskip\bigskip
\par {\bf 2.2. A parameterized family of separating squares and its main properties.}
\addtocontents{toc}{~~~~2.2. A parameterized family of separating squares and its properties. \hfill \thepage\par}
\medskip
\begin{definition}\lbl{DF-GTB} {\em Let $B\in\Omega\sm \SB^{\cl}$. By $\GB$ we denote the connected component of $\Omega\sm \SB^{\cl}$ containing $B$, and by  $\TB=\Tc(\GB)$ we denote the corresponding connected component of $\DS\sm\DO$ associated with $\GB$. We represent $\TB$ in the form $\Tc=\oi{b_{\TB}}{e_{\TB}}$ where $b_{\TB},e_{\TB}\in\DS$, $b_{\TB}\prec e_{\TB}$. See Definition \reff{ORD}.}
\end{definition}
\par By Lemma \reff{CM-X}, the component $\TB$ is well defined.\medskip
\par Our aim at this step of the proof is to introduce a certain parametrization of squares touching the square $\SB=S(\cs,\rs)$ and lying in $\GB$. Let $z\in\DS$ and let $r>0$. By $K_r(z)$ we denote a square with ``radius'' $r$ and center
$$
z_r:=z+\tfrac{r}{\rs}(z-\cs).
$$
\par  Since $\|z-\cs\|=\rs$, we have
$$
\|z_r-\cs\|=\|z+\tfrac{r}{\rs}(z-\cs)-\cs\|=
(1+\tfrac{r}{\rs})\|z-\cs\|=\rs+r
$$
so that, by part (iii) of Lemma \reff{SQ-PR},
\bel{T-SK}
K_r(z)~~\text{and}~~\SB~~\text{are touching squares}.
\ee
\par Furthermore, if $0<r_1\le r_2$, then
$$
\|z_{r_1}-z_{r_2}\|=\|z+\tfrac{r_1}{\rs}(z-\cs)-
(z+\tfrac{r_2}{\rs}(z-\cs))\|=r_2-r_1.
$$
Therefore, by part (i) of Lemma \reff{SQ-PR},
$$
K_{r_1}(z)\subset K_{r_2}(z)~~~
\text{whenever}~~~0<r_1\le r_2
$$
proving that the family of squares $\{K_r(z):r>0\}$ is {\it ordered with respect to inclusion}. This motivates us to introduces the following
\begin{definition} {\em Let  $z\in\TB$. By $K(z)$ we denote {\it the maximal (with respect to inclusion) element} of the family of squares
$$
\Kc(z):=\{K_r(z):r>0, K_r(z)\subset\Omega\}.
$$
\par We let $c_z$ and $r_z$ denote the center and the ``radius'' of  $K(z)$ respectively.}
\end{definition}
\par Thus $K(z)$ is the square of the {\it maximal diameter} belonging to the family of squares $\Kc(z)$. It can be represented in the form
$$
K(z)=S(c_z,r_z),~~~ z\in\TB
$$
where
\bel{CF-Z}
c_z=z+\tfrac{r_z}{R}(z-\cs).
\ee
\par See Figure 5.
\begin{figure}[h]
\center{\includegraphics[scale=0.9]{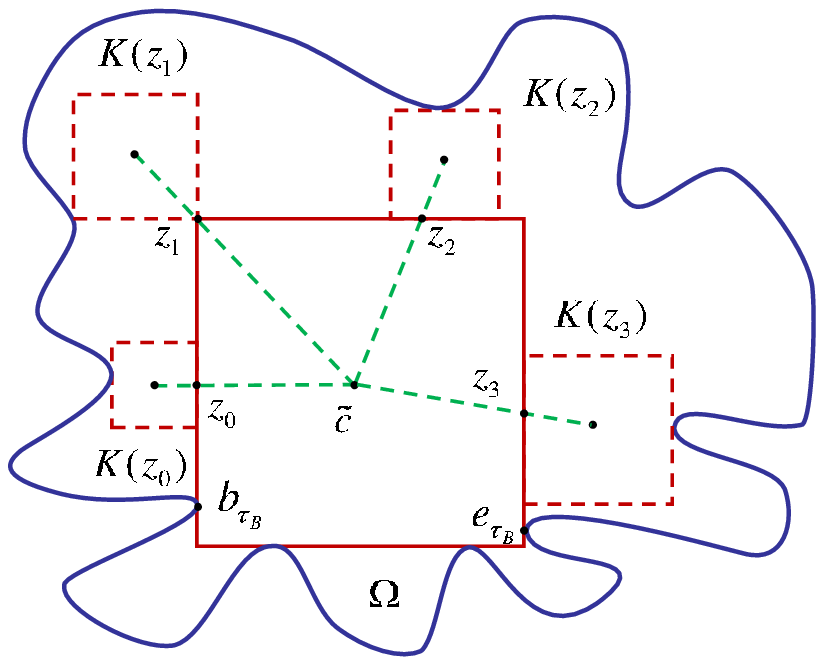}}
\caption{Examples of squares $K(z_i)$, $i=0,...3.$}
\end{figure}
\newpage
\par Let us describe several simple properties of the squares $K(z), z\in\TB$.
\begin{lemma}\lbl{PR-KZ} Let $z\in\TB$.\medskip
\par (a). The square $K(z)$ is well defined;
\medskip
\par (b). $K(z)$ and $\SB$ are touching squares such that $K(z)^{\cl}\cap\DS\cap\Omega\ne\emp$;
\medskip
\par (c). $K(z)\subset\Omega\sm\SB^{\cl}$ and $\dist(K(z),\DO\setminus\DS)=0$\,;
\par (d). The line segment $[\cs,c_z]$ lies in $\Omega$:
\bel{C-KSZ}
[\cs,c_z]\subset\Omega.
\ee
\par Furthermore,
\bel{Z-KSZ}
z\in K(z)^{\cl}\cap\SCL\cap\Omega\,;
\ee
\par (e). For every $u\in K(z)^{\cl}\cap\Omega$ there exists a path $\gamma$ which joins $u$ to $B$ in $\Omega$ such that $(\gamma\sm\{u\})\cap\SCL=\emp$.
\par In particular, this implies that $K(z)$ and $B$ belong to the same connected component of\, $\Omega\sm\SCL$ (i.e., the component $\GB$).
\end{lemma}
\par {\it Proof.} Since $\Omega$ is a bounded domain and $K(z)$ is the square of the maximal diameter from the family $\Kc(z)$, this square is well defined. This proves (a).
\par In turn, property (b) follows from \rf{T-SK}, and property (c) from the maximality of the square $K(z)$. Property (d) follows from the fact that the point $z\in\TB\subset\Omega$ and $\SB\cup K(z)\subset\Omega$.\medskip
\par Prove (e). Since $z\in\TB$, by Definition \reff{DF-GTB} and Lemma \reff{CM-X}, there exists a path $\gamma_z$ which joins $B$ to $z$ in $\Omega$ such that $\gamma_z\sm\{z\}\subset\GB$. Recall that $z\in\Omega$ so that for some $\ve>0$ small enough the $\ve$-neighborhood of $z$, the square $S(z,\ve)\subset\Omega$.
\par Clearly, $(\gamma_z\sm\{z\})\cap S(z,\ve)\ne\emp$
and $K(z)\cap S(z,\ve)\ne\emp$ so that there exist points $a\in\gamma\sm\{z\}$, and $b\in K(z)$ which belong to $S(z,\ve)$.
\par Let $\gamma_1$ be the arc of $\gamma$ from $B$ to $a$. Clearly, $S(z,\ve)\sm\SCL$ is an open connected set so that there exists a path $\gamma_2$ in $S(z,\ve)\sm\SCL$ joining $a$ to $b$. Finally, let $\gamma_3:=[b,u]$.
\par Let $\gamma:=\gamma_1\cup\gamma_2\cup\gamma_3$. Then $\gamma$ is a path which joins $u$ to $B$ in $\Omega$. Since
$$
\gamma_1\cap\SCL=\gamma_2\cap\SCL=\emp,
$$
$b\in K(z)$ and $u\in K(z)^{\cl}\cap\Omega$, the path $\gamma\sm\{u\}$ does not intersect $\SCL$.
\par Prove the second statement of part (e). Since $K(z)\cap\SCL=\emp$, we conclude that {\it every point} $u\in K(z)$ can be joined to $B$ by a path $\gamma\subset\Omega$ such that $\gamma\cap\SCL=\emp$.
Clearly, this implies that $K(z)$ and $B$ belong to the same connected component of\, $\Omega\sm\SCL$.
\par The lemma is proved.\bx\smallskip
\begin{lemma}\lbl{C-KZL} Let $z,y\in\TB, z\ne y$. Suppose that $z$ and $y$ lie on a side $[a,b]$ of the square $\SB$. 
\par (i). If $z,y\in(a,b)$, then
$$
|r_y-r_z|\le \frac{(\rs+r_y+r_z)\,\|y-z\|}{\dist(\{z,y\},\{a,b\})}\,;
$$
\medskip
\par (ii). If $z\in\{a,b\}$ and $y\in (a,b)$, then
$$
r_y\le r_z+\frac{(\rs+r_z)\,\|y-z\|}{\|y-h\|}
$$
where $h:=\{a,b\}\setminus\{z\}$.
\end{lemma}
\par {\it Proof.} Without loss of generality we may assume that $\cs=(0,-\rs)$, $a=(-\rs,0)$ and $b=(\rs,0)$.
See Figure 6.
\begin{figure}[h]
\center{\includegraphics[scale=0.7]{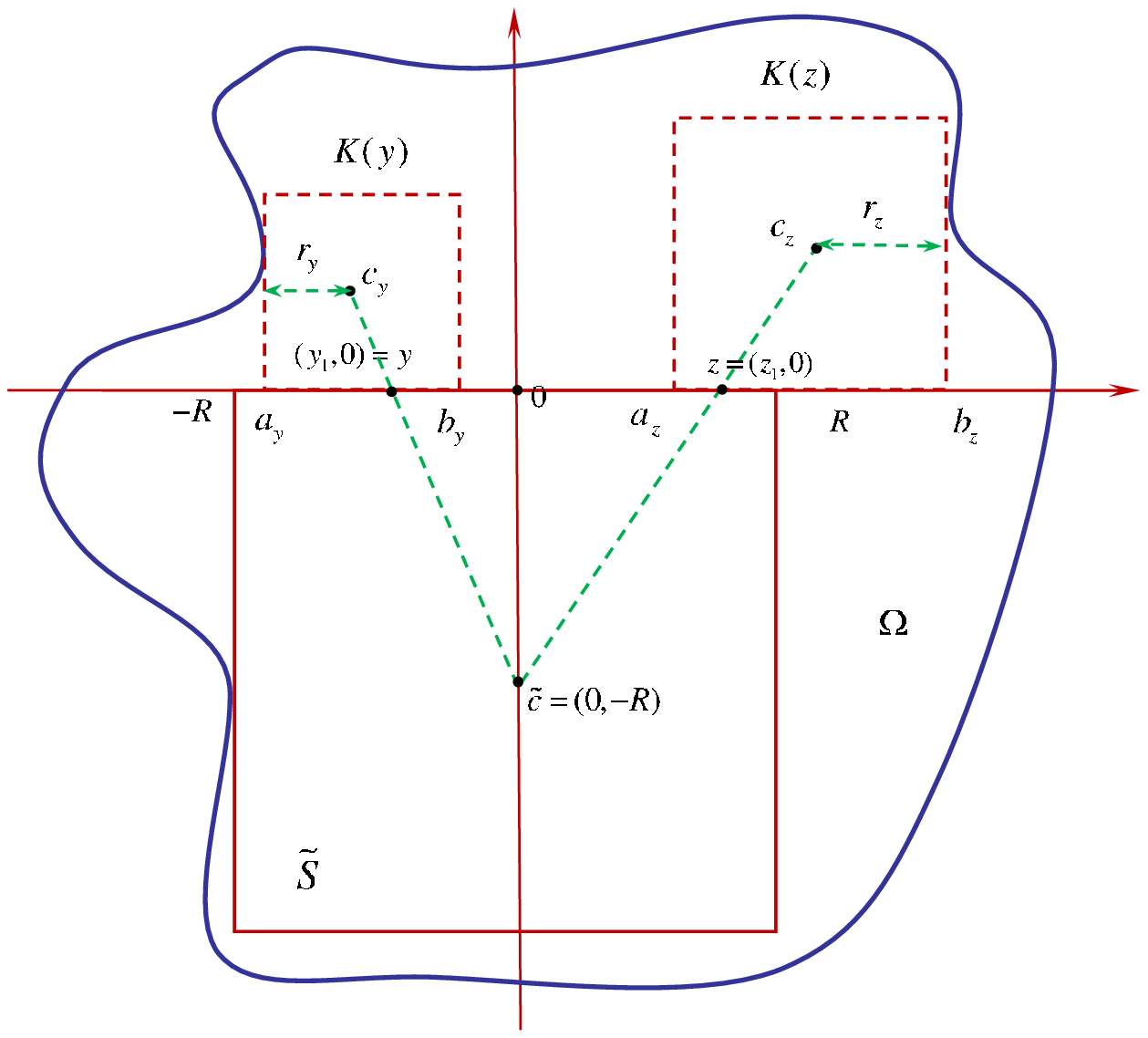}}
\caption{~}
\end{figure}
\par Since $y,z\in[a,b]\subset Ox$, we have $y=(y_1,0)$ and $z=(z_1,0)$ where $|y_1|\le\rs$ and $|z_1|\le\rs$. Since $K(z)$ and $\SB$ are touching squares, intersection of $K(z)^{\cl}$ with the axis $Ox$
is a closed line segment which coincides with a side of $K(z)$. Let $(a_z,0)\in Ox$ and $(b_z,0)\in Ox$ be the ends of this side so that
\bel{AZ-BZ}
K(z)^{\cl}\cap Ox=[(a_z,0),(b_z,0)].
\ee
In the same way we define points $(a_y,0),(b_y,0)\in Ox$; thus
$$
K(y)^{\cl}\cap Ox=[(a_y,0),(b_y,0)].
$$
\par Let us give explicit formulae for these points. By \rf{CF-Z},
$$
c_z=z+\tfrac{r_z}{\rs}(z-\cs)=
\left(z_1\left(1+\tfrac{r_z}{\rs}\right),r_z\right)
$$
so that
\bel{ABZ}
a_z=z_1\left(1+\tfrac{r_z}{\rs}\right)-r_z~~~\text{and}~~~
b_z=z_1\left(1+\tfrac{r_z}{\rs}\right)+r_z.
\ee
In the same way we obtain formulae for $a_y$ and $b_y$:
\bel{ABY}
a_y=y_1\left(1+\tfrac{r_y}{\rs}\right)-r_y~~~\text{and}~~~
b_y=y_1\left(1+\tfrac{r_y}{\rs}\right)+r_y.
\ee
\par Prove that either
\bel{O-1}
a_y\le a_z~~~\text{and}~~~b_y\le b_z
\ee
or
\bel{O-2}
a_z\le a_y~~~\text{and}~~~b_z\le b_y\,.
\ee
In fact, assume that both \rf{O-1} and \rf{O-2} do not hold. Then either
\bel{NO-1}
a_y<a_z~~~\text{and}~~~b_z<b_y
\ee
or
\bel{NO-2}
a_z<a_y~~~\text{and}~~~b_y<b_z\,.
\ee
\par Prove that \rf{NO-1} contradicts the maximality of the square $K(z)$. In fact, if \rf{NO-1} holds, then $K(z)^{\cl}\subset K(y)$ so that $K(z)^{\cl}\subset\Omega$. But this inclusion contradicts the equality
$$
\dist(K(z),\DO\setminus\DS)=0,
$$
see part (c) of Lemma \reff{PR-KZ}. In the same way we show that \rf{NO-2} is not true proving that either \rf{O-1} or  \rf{O-2} holds.\smallskip
\par We are in a position to prove part (i) of the lemma. Suppose that $y,z\in(a,b)$ and the option \rf{O-1} holds.
By \rf{ABZ} and \rf{ABY}, inequality $a_y\le a_z$ is equivalent to the inequality
$$
y_1\left(1+\tfrac{r_y}{\rs}\right)-r_y\le z_1\left(1+\tfrac{r_z}{\rs}\right)-r_z.
$$
Hence
$$
r_y-r_z\ge \frac{(\rs+r_z)(y_1-z_1)}{\rs-y_1}\,.
$$
In turn, inequality $b_y\le b_z$ implies that
\bel{ZY-2}
r_y-r_z\le \frac{(\rs+r_z)(z_1-y_1)}{\rs+y_1}\,.
\ee
\par In the same way we prove that \rf{O-2} implies the following:
$$
r_z-r_y\ge \frac{(\rs+r_y)(z_1-y_1)}{\rs-z_1}~~\text{and}~~
r_z-r_y\le \frac{(\rs+r_y)(y_1-z_1)}{\rs+y_1}\,.
$$
\par Summarizing these estimates, we obtain
$$
|r_y-r_z|
\le
|y_1-z_1|\,\max\left\{\frac{\rs+r_z}{\rs+y_1},
\frac{\rs+r_y}{\rs+z_1},\frac{\rs+r_z}{\rs-y_1},
\frac{\rs+r_y}{\rs-z_1}\right\}\,.
$$
But $|y_1-z_1|=\|y-z\|$ so that
$$
|r_y-r_z|
\le
\frac{\|y-z\|\,(\rs+r_y+r_z)}
{\min\{\rs+y_1,\rs-y_1,\rs+z_1,\rs-z_1\}}
=\frac{\|y-z\|\,(\rs+r_y+r_z)}
{\dist\{\{y,z\},\{a,b\}\}}
$$
proving part (i) of the lemma.
\medskip
\par Prove (ii). Let $z=b$ and let $y\in(a,b)$ so that $z_1=\rs$ and $-\rs<y_1<\rs$. By \rf{ABZ} and \rf{ABY},
$$
a_y=y_1\left(1+\tfrac{r_y}{\rs}\right)-r_y
=y_1+r_y\left(\tfrac{y_1}{\rs}-1\right)<y_1
$$
and
$$
a_z=z_1\left(1+\tfrac{r_z}{\rs}\right)-r_z=
\rs\left(1+\tfrac{r_z}{\rs}\right)-r_z=\rs
$$
so that $a_y< a_z$. Therefore, by \rf{O-1}, $b_y\le b_z$. Hence, by \rf{ZY-2},
$$
r_y-r_z\le \frac{(\rs+r_z)(z_1-y_1)}{\rs+y_1}=
\frac{(\rs+r_z)\|y-z\|}{\rs+y_1}\,.
$$
Since $a=(-R,0)$ and $y=(y_1,0)$, we have $\rs+y_1=\|y-a\|$ proving part (ii) of the lemma in the case under consideration. In the same fashion we prove (ii) whenever $z=a$.
\par The proof of the lemma is complete.\bx
\begin{lemma}\lbl{E-KZ} Let $z\in\TB$ and let $\ve>0$. There exists $\delta>0$ such that for every $y\in\TB$, $\|y-z\|<\delta$, the following inclusion
\bel{I-CN}
K(y)\subset[K(z)]_\ve
\ee
holds. Recall that the symbol $[\,\,\cdot\,\,]_\ve$ denotes the $\ve$-neighborhood of a set.
\end{lemma}
\par {\it Proof.} Clearly, $[K(z)]_\ve$ is a square with center $c_z$ and ``radius'' $r_z+\ve$, i.e.,
$$
[K(z)]_\ve=S(c_z,r_z+\ve).
$$
\par By part (i) of Lemma \reff{SQ-PR}, inclusion \rf{I-CN} is equivalent to the inequality
\bel{F-IN}
\|c_y-c_z\|+r_y\le r_z+\ve.
\ee
\par Let us consider two cases.\medskip
\par {\it The first case:} $z$ is not a vertex of the square $\SB$, i.e.,
$$
\tau:=\dist(z,V_{\TB})>0.
$$
Here $V_{\TB}$ is the family of vertices of $\SB$ which belong to $\TB$. In particular, every point $y\in\TB$ such that $\|y-z\|<\tau/2$ belongs to the same side of $\SB$ as the point $z$. Furthermore,
\bel{DY-2}
\dist(y,V_{\TB})\ge\tau/2>0.
\ee
\par By part (i) of Lemma \reff{C-KZL},
$$
|r_y-r_z|\le\frac{(\rs+r_y+r_z)\|y-z\|}{\dist(\{y,z\},V_{\TB})}
$$
so that, by \rf{DY-2},
$$
|r_y-r_z|\le (2/\tau)(\rs+r_y+r_z)\|y-z\|=\gamma_1\|y-z\|
$$
where $\gamma_1:=2(\rs+r_y+r_z)/\tau$.
\par By \rf{CF-Z},
$$
c_z=z+\tfrac{r_z}{\rs}(z-\cs)~~~
\text{and}~~~c_y=y+\tfrac{r_y}{\rs}(y-\cs)
$$
so that
$$
\|c_y-c_z\|\le \|y-z\|+\tfrac{|r_y-r_z|}{\rs}\|z-\cs\|
+\tfrac{r_y}{\rs}\|y-z\|\,.
$$
Since $\|z-\cs\|=\rs$, we obtain
\bel{H1}
\|c_y-c_z\|\le
\left(1+\tfrac{r_y}{\rs}\right)\|y-z\|+|r_y-r_z|\,.
\ee
Hence,
$$
\|c_y-c_z\|\le
\left(1+\tfrac{r_y}{\rs}+\gamma_1\right)\|y-z\|=
\gamma_2\|y-z\|
$$
with $\gamma_2:=1+\tfrac{r_y}{\rs}+\gamma_1$.
\par Now we are in a position to estimate the left-hand side of \rf{F-IN}:
$$
\|c_y-c_z\|+r_y\le \|c_y-c_z\|+|r_y-r_z|+r_z\le
\gamma_2\|y-z\|+\gamma_1\|y-z\|+r_z
=(\gamma_1+\gamma_2)\|y-z\|+r_z.
$$
\par This proves that whenever $\|y-z\|<\delta$ with $\delta:=\min\{\tau/2,\ve/(\gamma_1+\gamma_2)\}$ the inequality \rf{F-IN} holds.\bigskip
\par {\it The second case: $z$ is a vertex of $\SB$.} Let $y\in\TB$, $\|y-z\|<\rs/2$. Hence, $\|y-a\|>\rs/2$ for every vertex $a$ of $\SB$, $a\ne z$. Then, by part (ii) of Lemma \reff{C-KZL},
\bel{H2}
r_y\le r_z+\frac{(\rs+r_z)\,\|y-z\|}{(\rs/2)}
=r_z+2(1+r_z/\rs)\,\|y-z\|\,.
\ee
\par Prove inequality \rf{F-IN}. If $r_z\ge r_y$, then, by \rf{H1},
\bel{H2-T}
\|c_y-c_z\|+r_y\le \left(1+\tfrac{r_y}{\rs}\right)\|y-z\|+|r_z-r_y|+r_y=
\left(1+\tfrac{r_y}{\rs}\right)\|y-z\|+r_z\,.
\ee
\par If $r_z<r_y$, then, by \rf{H1} and \rf{H2},
\be
\|c_y-c_z\|+r_y&\le& \left(1+\tfrac{r_y}{\rs}\right)\|y-z\|+(r_y-r_z)+r_y\nn\\
&\le&
\left(1+\tfrac{r_y}{\rs}\right)\|y-z\|+2(r_y-r_z)+r_z\nn\\
&\le&
\left(1+\tfrac{r_y}{\rs}\right)\|y-z\|+
4\left(1+\tfrac{r_z}{\rs}\right)\|y-z\|+r_z\nn
\ee
so that
\bel{H3}
\|c_y-c_z\|+r_y\le 5\left(1+\tfrac{r_y}{\rs}+
\tfrac{r_z}{\rs}\right)\|y-z\|+r_z\,.
\ee
\par Combining this estimate with \rf{H2-T}, we conclude that inequality \rf{H3} is true {\it for all} choices of $y$. This shows that inequality \rf{F-IN} is satisfied provided $\|y-z\|<\delta$ where
$\delta:=\min\{\rs/2,\,\ve/5(1+(r_y+r_z)/\rs)$.
The proof of the lemma is complete.\bx
\begin{lemma}\lbl{T-KZ} Let $K$ be a square such that  $K\subset \GB$,
\bel{K-MG}
\KCL\cap \TB\ne\emp~~~\text{and}~~~\KCL\cap\DO\ne\emp.
\ee
\par Suppose that $B\in \GB\sm \KCL$. Then there exists at most one connected component $\TV=\TV(K)$ of the set $\TB\sm K^{\cl}$ which has the following property:
\bel{B-PRC}
\exists~ y\in\TV~~\text{and a path}~ \gamma_y~\text{joining}~y~\text{to}~B~\text{such that}~~\gamma_y\sm\{y\}\subset \GB\sm \KCL.
\ee
\par See Figure 7.
\begin{figure}[h]
\center{\includegraphics[scale=0.89]{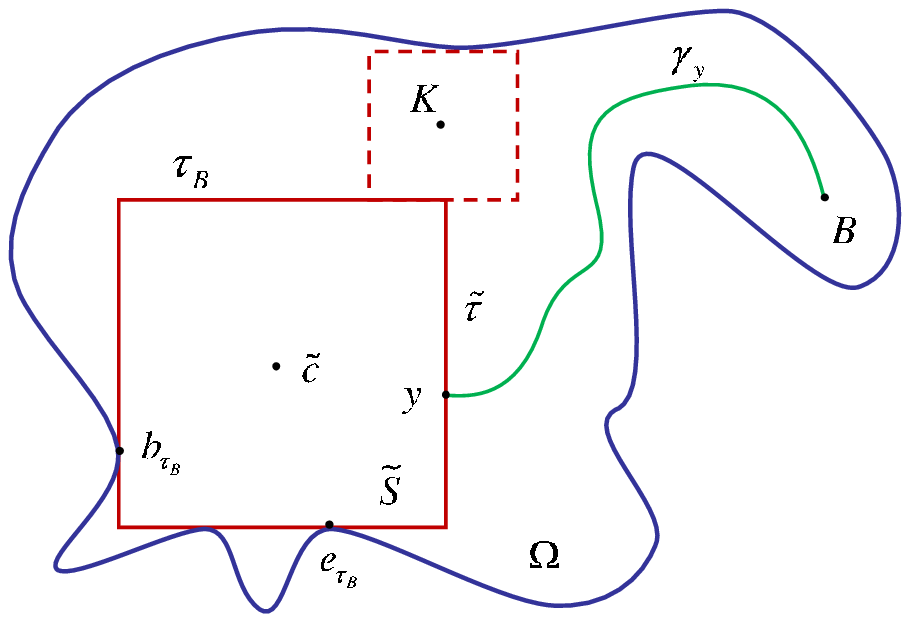}}
\caption{The path $\gamma_y$ connects $y$ to $B$ in $\GB\sm \KCL$.}
\end{figure}
\medskip
\par Furthermore, every point $x\in\TV$ has this property, i.e., it can be joined to $B$ by a path $\gamma_x$ such that $\gamma_x\sm\{x\}\subset \GB\sm \KCL$.  
\end{lemma}

\bigskip
\par {\it Proof.} Since
$$
K\subset \GB\subset \RT\sm\SCL~~~\text{and}~~~\KCL\cap \TB\ne\emp,
$$
we have $\SB\cap K=\emp$ and $\SCL\cap\KCL\ne\emp$, so that $\SB$ and $K$ are {\it touching squares}. Clearly, for each $p\in\SCL\cap\KCL$ we have
$$
[c_K,p)\subset K\subset \GB
$$
so that, by Lemma \reff{CM-X}, $p\in\TB$. Thus
\bel{T-K}
\TB\cap\KCL=\DS\cap \KCL=\SCL\cap\KCL,
\ee
so that, by part (iii) of Lemma \reff{SQ-PR},  $\TB\cap\KCL$ is either a line segment or a point. In particular, $\TB\sm \KCL$ has at most two connected components. Prove that $\TB\sm \KCL$ has at most one connected component $\TV$ satisfying \rf{B-PRC}.
\par Suppose that there exist two distinct connected components $\Tc'$ and $\Tc''$ of $\TB\sm\KCL$, points $y'\in\Tc'$ and $y''\in\Tc''$, paths $\Gamma'$ and $\Gamma''$ joining $B$ to $y'$ and $y''$ respectively such that
$$
\Gamma'\sm\{y'\}, \Gamma''\sm\{y''\}\subset \GB\sm\KCL.
$$
See Figure 8.\smallskip
\begin{figure}[h]
\center{\includegraphics[scale=1]{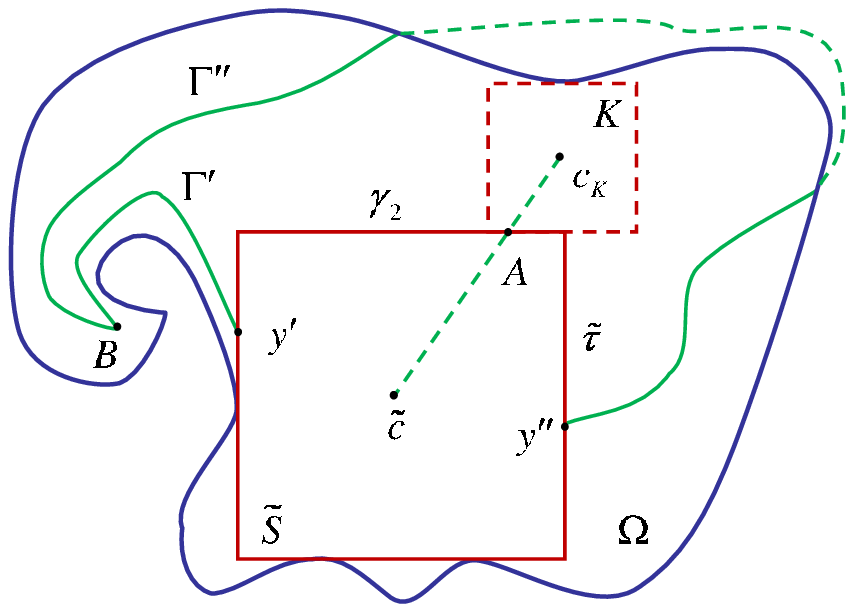}}
\caption{Paths $\Gamma'$ and $\Gamma''$ join $B$ to $y'$ and $y''$ in $\GB\sm\KCL$.}
\end{figure}
\medskip
\par We may assume that $y'\prec y''$. Since $\Tc'$ and $\Tc''$ are distinct connected components of $\TB\sm\KCL$,
we have
$$
\DS\cap \DKZ\subset \oi{y'}{y''}.
$$
\par By part (ii) of Lemma \reff{S-PH}, there exist a simple path $\gamma_1\subset \Gamma'\cup\Gamma''$ which joins $y'$ to $y''$ such that
\bel{G1-IN}
\gamma_1\sm\{y',y''\}\subset \GB\sm \KCL.
\ee
\par Let $\gamma_2:=\ci{y'}{y''}$ and let
$$
\tgm:=\gamma_1\cup\gamma_2.
$$
We know that $\gamma_2=\ci{y'}{y''}\subset \DS\cap \Omega$ and, by \rf{G1-IN}.
$\gamma_1\sm\{y',y''\}\subset \GB\subset\Omega$. This shows that $\gamma_2\cap (\gamma_1\sm\{y',y''\})=\emp$ so that  the path $\tgm$ is a {\it simple polygon in} $\Omega$. Hence, by part (i) of Statement \reff{ST-J}, the set $\RT\setminus\tgm$ consists of exactly two connected components - the interior $D_{int}$ (which is a bounded set), and the exterior component $D_{ext}$ (which is an unbounded set). Furthermore, $\tgm=\partial D_{int}=\partial D_{ext}$. Since $\Omega$ is a simply connected domain and $\tgm\subset\Omega$ is a simple polygon, by part (ii) of  Statement \reff{ST-J}, $D_{int}\subset\Omega.$
\par We also notice that $\tgm$ is a compact subset of $\Omega$ so that
\bel{DS-G}
\dist(\tgm,\DO)>0.
\ee
\par Prove that the centers of squares $\SB$ and $K$, the points $\cs$ and $c_K$, belong to distinct connected components of $\RT\sm\tgm$.
\par Since $\SB$ and $K$ are touching squares, by part (iii) of Lemma \reff{SQ-PR},
$$
[\cs,c_K]\cap\SCL\cap\KCL=\{A\}
$$
for some $A\in\RT$, see \rf{IN-S}. Hence, by \rf{T-K},
$A\in\TB\cap[\cs,c_K]$. On the other hand, $A$ is the unique point of intersection of $\DS$ and $[\cs,c_K]$. Since $\TB\subset\DS$, we conclude that $\{A\}=\TB\cap[\cs,c_K]$.
\par Furthermore, since $K\subset \GB$ and $\SB\cap\tgm=\emp$,
$$
\{A\}=\tgm\cap[\cs,c_K].
$$
\par We also notice that, by Definition \reff{S-C}, $[\cs,c_K]$ {\it strictly crosses} the polygon $\tgm$, so that, by Statement \reff{S-CR}, $\cs$ and $c_K$ belong to  distinct connected components of $\RT\sm\tgm$.
\par Since $\tgm\cap K=\emp$, for every $x\in K$ the line segment $[x,c_K]$ does not intersect $\tgm$ so that $K$ lie in the same connected component of $\RT\sm\tgm$ as $c_K$. The same is true for the square $\SB$ and $\cs$. This proves that the squares $\SB$ and $K$ lie in distinct connected components of $\RT\sm\tgm$.
\par Thus either $K\subset D_{int}$ or $\SB\subset D_{int}$. Recall that $D_{int}\subset \Omega$ and $\partial D_{int}=\tgm$ so that, by \rf{DS-G},
\bel{D-DO}
\dist(\partial D_{int},\DO)>0.
\ee
\par This inequality immediately leads us to a contradiction. In fact, if $K\subset D_{int}$, then $\KCL\subset (D_{int})^{\cl}$ so that, by \rf{D-DO},
$\dist(\KCL,\DO)>0$.
But, by the lemma's hypothesis, $\KCL\cap\DO\ne\emp$, see \rf{K-MG}, a contradiction.
\par On the other hand, if $\SB\subset D_{int}$, then the
same consideration shows that $\dist(\SCL,\DO)>0$
which contradicts to the assumption that $\SCL\cap\DO\ne\emp$.
\par It remains to show that {\it every} point $x\in\TV$ can be joined to $B$ by a path $\gamma_x$ such that
$$
\gamma_x\sm\{x\}\subset \GB\sm \KCL.
$$
We prove this statement using precisely the same arguments as used in the proof of Statement A from Lemma \reff{CM-X}. We leave the details to the interested reader.
\par The proof of the lemma is complete.\bx
\bigskip
\par {\bf 2.3. The final step of the proof of ``The Square Separation Theorem''.}
\addtocontents{toc}{~~~~2.3. The final step of the proof of ``The Square Separation Theorem''.\hfill \thepage\par\VST}
At this step we make the following
\begin{assumption}\lbl{MA-KZ} For every $z\in\TB$ the following conditions are satisfied:\medskip
\par (i). $B\notin K(z)^{\cl}$;
\medskip
\par (ii). There exist a point $z'\in\TB$ and a path $\gamma$ joining $z'$ to $B$ in $\Omega$ such that $$\gamma\sm\{z'\}\subset\GB\sm K(z)^{\cl}.$$
\end{assumption}
\par We will show that this assumption leads us to a contradiction which immediately implies the statement of  Theorem \reff{N-SQ}.
\par Assumption \reff{MA-KZ} and Lemma \reff{T-KZ} motivate the following
\begin{definition}\lbl{A-T} {\em Let $z\in \TB$. By $\TBZ$ we denote a connected component of $\TB\sm K(z)^{\cl}$ having the following property: for every point $y\in\TBZ$ there exists a path $\gamma$ which connects $y$ to $B$ in $\Omega$ such that
$$
\gamma\sm\{y\}\subset\GB\sm K(z)^{\cl}.
$$
We refer to $\TBZ$ as a {\it $B$-accessible component} of the set $\TB\sm K(z)^{\cl}$ (with respect to $z$).}
\end{definition}
\par By Assumption \reff{MA-KZ} and Lemma \reff{T-KZ}, the $B$-accessible component $\TBZ$ is well defined and {\it non-empty} for each $z\in\TB$.
\par Thus for every $z\in\TB$ the set $\TB\sm K(z)^{\cl}$  contains at least one and at most two connected components. One of them is the $B$-accessible component $\TBZ$ consisting of {\it all} points of $\TB$ connected to $B$ by paths which lie in $\GB\sm K(z)^{\cl}$. Another connected component (if it exists) consists of ``$B$-inaccessible'' points, i.e., those points $y\in\TB$ for which any path connecting $y$ to $B$ in $\GB$ crosses $K(z)^{\cl}$. See Figure 9.\newpage
\begin{figure}[h]
\center{\includegraphics[scale=1]{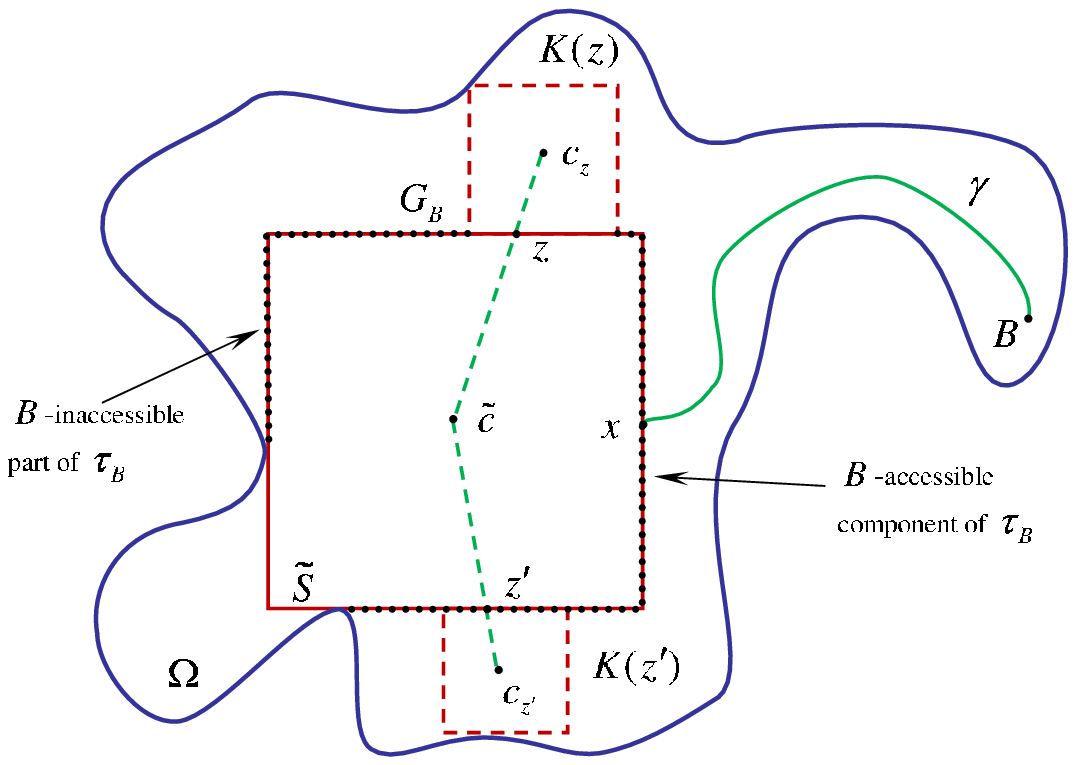}}
\caption{``$B$-accessible'' and ``$B$-inaccessible'' subsets of $\TB$.}
\end{figure}
\medskip
\par The next definition enables us to specify the position of the $B$-accessible component $\TBZ$ with respect to the interval $\DS\cap K(z)^{\cl}$.
\begin{definition} {\em By $\TBP$ we denote a set consisting of all points $z\in\TB$ such that
$$
x\prec y~~\text{for every}~~x\in\TB\cap K(z)^{\cl}~~
\text{and every}~~y\in\TBZ.
$$
\par Correspondingly, $\TBM$ is a subset of $\TB$ consisting of all points $z$ such that
$$
y\prec x~~\text{for every}~~x\in\TB\cap K(z)^{\cl}~~
\text{and every}~~y\in\TBZ.
$$}
\end{definition}
\par In particular, the point $z$ on Figure 9 belongs to $\TBP$ while the point $z'$ on this picture belongs $\TBM$. Note that, by Lemma \reff{T-KZ},
\bel{IN-E}
\TBP\cap\TBM=\emp.
\ee
In turn, by Assumption \reff{MA-KZ},
\bel{P-TN}
\TBP\cup\TBM=\TB
\ee
so that {\it $\TBP$ and $\TBM$ is a partition of $\TB$}.
\par Our goal at this step of the proof is to show that representation \rf{P-TN} leads to a contradiction. Our proof of this fact relies on the following two lemmas
which state that $\TBP$ and $\TBM$ are open subsets of $\TB$, and, under Assumption \reff{MA-KZ}, these sets are non-empty.
\begin{lemma}\lbl{A-KZ} The sets $\TBP$ and $\TBM$ are open subsets of $\TB$ in the topology induced by the Euclidean metric on $\TB$. In other words, for each $z\in\TBP$ there exists $\ve>0$ such that every point $y\in\TB$, $\|y-z\|<\ve$, belongs to $\TBP$ (and the same statement is true for $\TBM$).
\end{lemma}
\par {\it Proof.} Let $z\in\TBP$. As we have noted above, the set $\TBZ$ of all $B$-accessible points is non-empty
so that there exists a point $z_1\in\TBZ$. Recall that $z_1\in \TB\sm K(z)^{\cl}$. By Definition \reff{A-T}, there exists a path $\gamma_1$ which connects $z_1$ to $B$ in $\Omega$ such that $\gamma_1\sm{z_1}\subset\GB\sm K(z)^{\cl}$. Furthermore, since $z\in\TBP$, we have
$$
x\prec z_1~~\text{for every}~~x\in\TB\cap K(z)^{\cl}.
$$
\par Let $\ve_1:=\dist(K(z)^{\cl},\gamma_1)$. Since
$\gamma_1\sm{z_1}\subset\GB\sm K(z)^{\cl}$, the path $\gamma_1$ and $K(z)^{\cl}$ have no common points, so that $\ve_1>0$. Since $z_1\in\gamma_1$, we have $z_1\notin[K(z)]_{\ve_1}$ so that
\bel{Z-1}
p\prec z_1~~\text{for every}~~p\in\TB\cap[K(z)]_{\ve_1}.
\ee
\par By Lemma \reff{E-KZ}, there exists $\delta>0$ such that for every $y\in\TB$, $\|y-z\|<\delta$, we have
$K(y)\subset [K(z)]_{\ve_1}$. Hence,
$$
K(y)\cap\gamma_1=\emp~~~\text{for every}~~~y\in\TB,~ \|y-z\|<\delta.
$$
See Figure 10.\medskip
\begin{figure}[h]
\center{\includegraphics[scale=1]{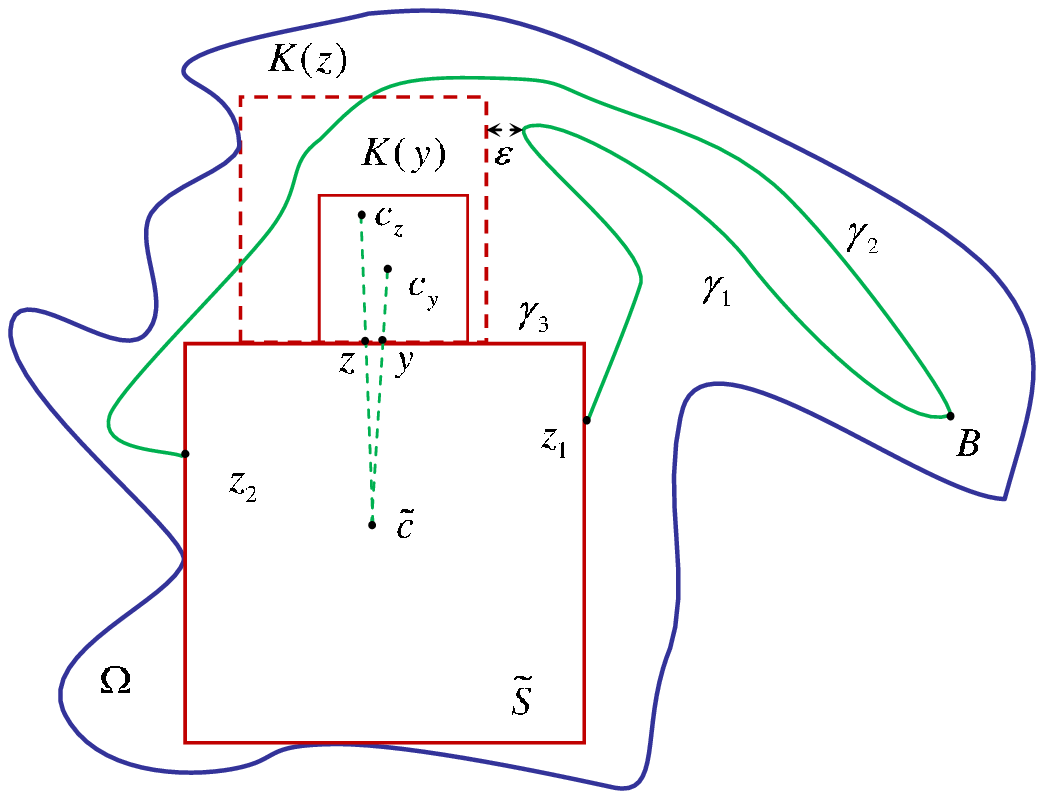}}
\caption{The path $\gamma_1$ joins $z_1$ to $B$ in $\GB\sm K(z)^{\cl}$.}
\end{figure}
\medskip
\par Prove that
$$
y\in\TBP~~~\text{for every}~~~y\in\TB,~ \|y-z\|<\delta.
$$
\par Suppose that there exists $y\in\TB$ such that $\|y-z\|<\delta$ but $y\notin\TBP$. Since the square $K(Y)\subset [K(z)]_{\ve_1}$, by \rf{Z-1},
$$
p\prec z_1~~\text{for every}~~p\in\TB\cap K(y)^{\cl}.
$$
\par By \rf{IN-E} and \rf{P-TN}, $y\in\TBM$ so that there exists a point $z_2\in\TB\sm K(y)^{\cl}$ such that
$$
z_2\prec x~~\text{for every}~~x\in\TB\cap K(y)^{\cl}.
$$
Furthermore, there exists a path $\gamma_2$ joining $z_2$ to $B$ in $\Omega$ such that $\gamma_2\sm\{z_2\}\subset \GB\sm K(y)^{\cl}$. See Figure 10.
\par Thus the point $B$ can be joined  by paths $\gamma_i$ in $\Omega$ to the points $z_i, i=1,2$ which belong to {\it distinct} connected components of $\TB\sm K(y)^{\cl}$. These paths have the following property: $\gamma_i\sm\{z_i\}\subset \GB\sm K(y)^{\cl}$, $i=1,2$. Furthermore, the square $K=K(y)$ satisfies conditions \rf{K-MG} of Lemma \reff{T-KZ}.
\par However, by this lemma, $B$ can be joined to at most {\it one} connected component of the set $\TB\sm K(y)^{\cl}$ by a path of such a kind, a contradiction. This contradiction proves that each point $y\in\TB$ in the $\delta$-neighborhood of $z$  belongs to $\TBP$.
\par In the same way we prove a similar statement for the set $\TBM$.
\par The lemma is completely proved.\bx
\begin{lemma}\lbl{N-EM} Under Assumption \reff{MA-KZ} both $\TBP$ and $\TBM$ are non-empty subsets of $\TB$.
\end{lemma}
\par {\it Proof.} Let us prove that $\TBM\ne\emp$.
\par Suppose that $\TBM=\emp$. Since $\TBP$ and $\TBM$ are a partition of $\TB$, we conclude that $\TB=\TBP$. This equality implies the following\medskip
\par {\it Statement B.} For every $z\in\TB$ there exists a point $y\in\TB\sm K(z)^{\cl}$ such that:\medskip
\par (i). $x\prec y$ for every $x\in\TB\sm K(z)^{\cl}$;
\par (ii). there exists a path $\gamma_y$ connecting $y$ to $B$ in $\Omega$ such that $\gamma_y\sm\{y\}\subset \GB\sm \KZC$.\medskip
\par See Figure 11. \medskip
\begin{figure}[h]
\center{\includegraphics[scale=0.75]{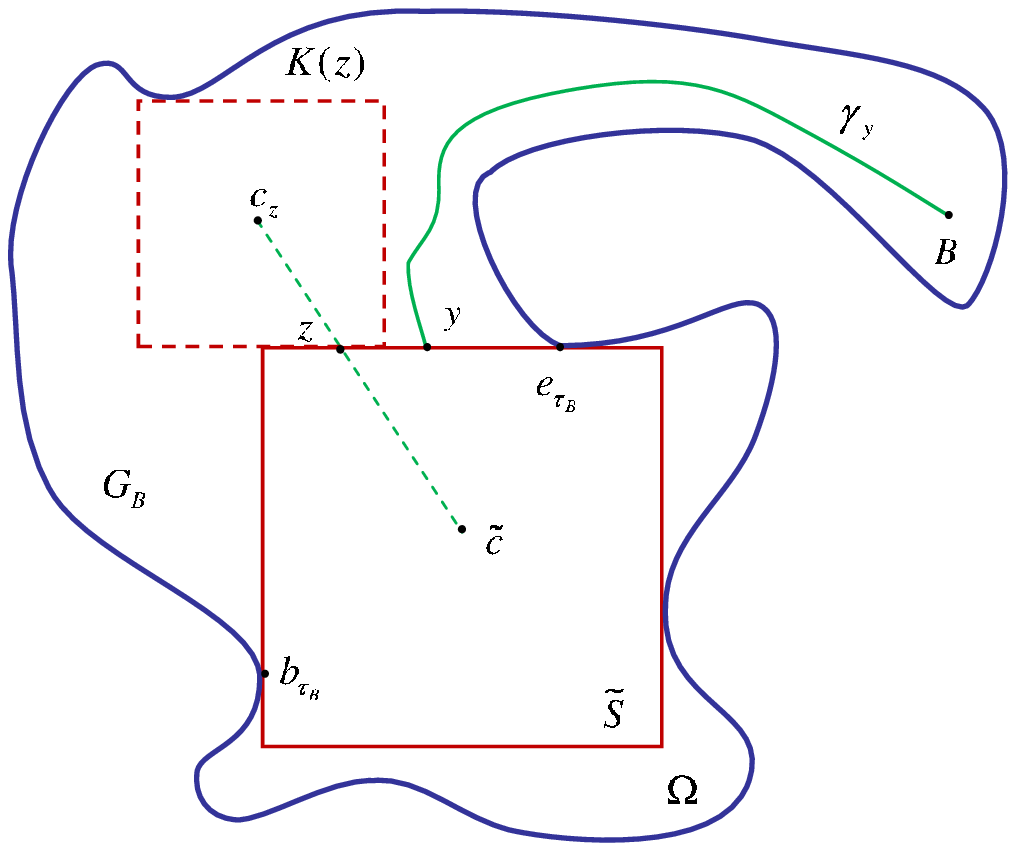}}
\caption{The path $\gamma_y$ connects $y$ to $B$ in
$\GB\sm \KZC$.}
\end{figure}
\medskip
\par Prove that Statement B leads to a contradiction whenever $z$ tends to the point $\ETB$ along $\TB$.
As in Lemma \reff{C-KZL}, without loss of generality we may assume that $\cs=(0,-\rs)$ where $\rs$ is the ``radius'' of $\SB$. Furthermore, $z\in[a,b)$ and $\ETB\in(a,b]$ where $a=(-\rs,0)$ and $b=(\rs,0)$. Thus $[a,b]$ is a side of $\SB$ lying on the real axes.
\par Let $z=(z_1,0)$ and $\ETB=(h,0)$ where $-\rs\le z_1<h$ and $-\rs<h\le\rs$.
\par We use the same notation as in Lemma \reff{C-KZL}. In particular, as in formulas \rf{AZ-BZ} and \rf{ABZ},
\bel{K-ABZ}
K(z)^{\cl}\cap Ox=[a_z,b_z]
\ee
where
$$
a_z=z_1\left(1+\tfrac{r_z}{\rs}\right)-r_z~~~\text{and}~~~
b_z=z_1\left(1+\tfrac{r_z}{\rs}\right)+r_z.
$$
See Figure 12.
\begin{figure}[h]
\center{\includegraphics[scale=1]{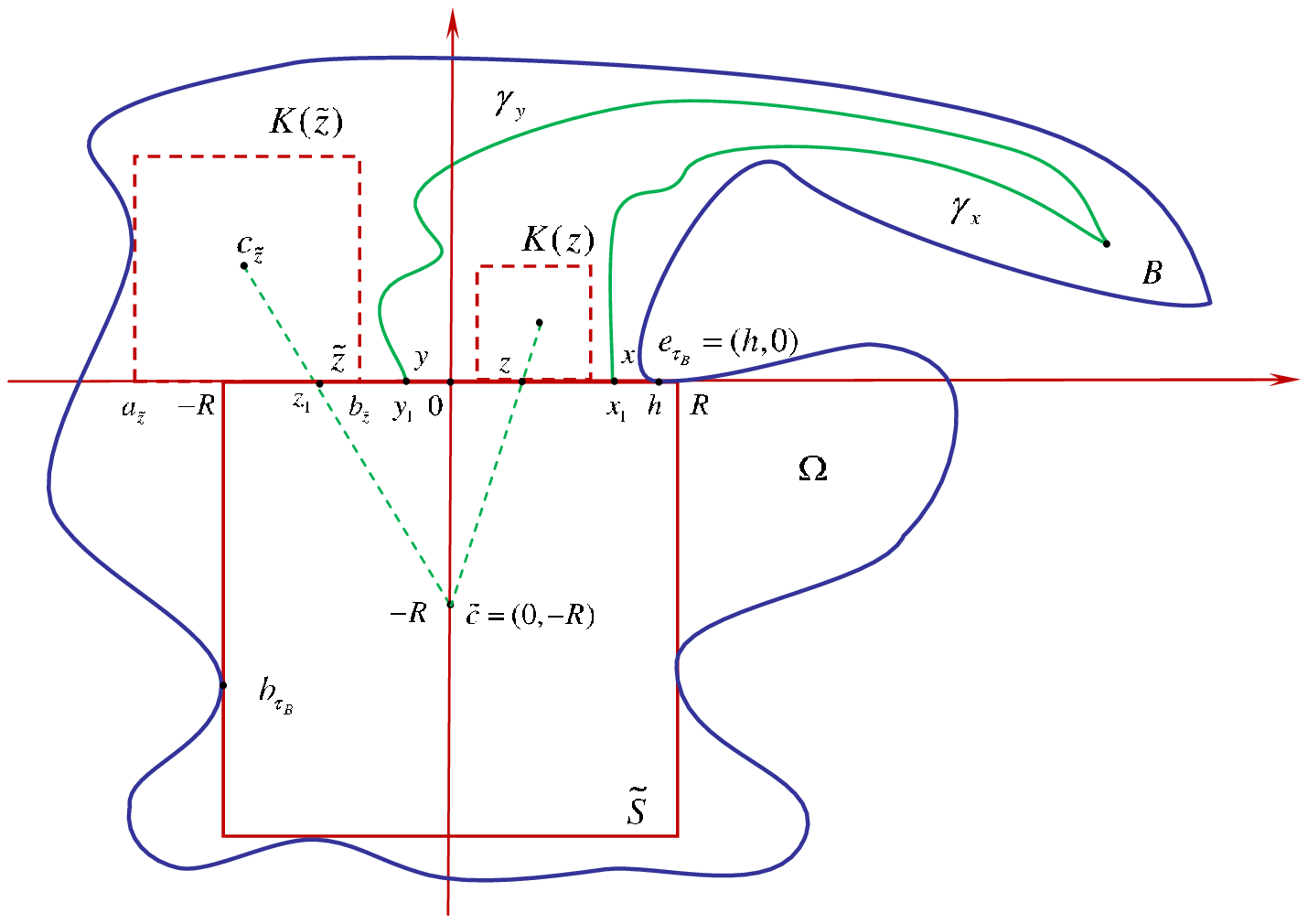}}
\caption{~}
\end{figure}
\medskip
\par Let $z\in\TB$ and $z\to \ETB$, i.e, $z_1\to h$. Consider two cases.\medskip
\par {\it The first case.} Let us assume that
\bel{L-S}
\limsup_{z\to\ETB,\, z\in\TB}r_z=L>0.
\ee
\par Prove that in this case there exists $\bz=(\bz_1,0)\in[a_z,\ETB)$ such that
\bel{E-L}
\ETB\in [a_{\bz},b_{\bz}].
\ee
\par Note that, since $a_{\bz}\le \bz_1<h$, property
\rf{E-L} is equivalent to the inequality $h\le b_{\bz}$.
\par Simple calculations show that if 
\bel{ZP-1}
r_{\bz}\ge L/2~~~\text{and}~~~ \|\bz-\ETB\|\le \tfrac14\min\{1,L/\rs\}\,\|\ETB-a\|,
\ee
then \rf{E-L} holds. In fact, since $r_{\bz}\ge L/2$, we obtain
$$
r_{\bz}\left(1+\tfrac{\bz_1}{\rs}\right)\ge \tfrac{L}{2}\left(1+\tfrac{\bz_1}{\rs}\right)
=\tfrac{L}{2\rs}\,\|\bz-a\|.
$$
Furthermore, by \rf{ZP-1},
$$
\|\bz-a\|\ge\|\ETB-a\|-\|\bz-\ETB\|\ge\|\ETB-a\|-
\tfrac14\|\ETB-a\|=\tfrac34\|\ETB-a\|
$$
proving that
$$
r_{\bz}\left(1+\tfrac{\bz_1}{\rs}\right)\ge
\tfrac{3L}{8\rs}\,\|\ETB-a\|.
$$
Hence, by \rf{ZP-1},
$$
r_{\bz}\left(1+\tfrac{\bz_1}{\rs}\right)\ge \|\bz-\ETB\|=
h-\bz_1
$$
so that
$$
b_{\bz}=\bz_1\left(1+\tfrac{r_{\bz}}{\rs}\right)+r_{\bz}\ge h
$$
proving \rf{E-L}.
\par Of course, condition \rf{L-S} guarantees the existence of a point $\bz\in\TB$ satisfying requirements \rf{ZP-1}.
\par Combining \rf{E-L} with \rf{K-ABZ} we conclude that $K(\bz)^{\cl}\ni\ETB$ so that the point $y$ satisfying conditions of part (i) of  Statement B {\it does not exist}. This contradiction shows that equality \rf{L-S} does not hold.\medskip
\par {\it The second  case.}
\bel{L-S0}
\lim_{z\to\ETB,\, z\in\TB}r_z=0.
\ee
\par Let $\tz=(\tz_1,0)$, $-\rs\le \tz_1<h$, and let
$$
K(\tz)^{\cl}\cap Ox=[a_{\tz},b_{\tz}]
$$
By Statement B, there exist a point $y=(y_1,0),$  $b_{\tz}<y_1<h,$ and a path $\gamma_y$ which joins $y$ to $B$ in $\Omega$ such that $\gamma_y\sm\{y\}\subset\GB\sm K(y)^{\cl}$. See Figure 12.
\par Let $\ve:=\dist(\gamma_y,\DO)$. Since $\gamma_y$ is a compact subset of $\Omega$, the number $\ve$ is positive. Note that the point $\ETB=(h,0)\in\DO$ so that
\bel{EP-1}
\dist(\gamma_y,\ETB)\ge \ve.
\ee
\par By \rf{L-S0}, there exist $\delta\in(0,\ve/4)$ such that
$$
r_z<\ve/8~~\text{for every}~~~z\in\TB,\|z-\ETB\|<\delta.
$$
See Figure 12.\bigskip
\par Fix such a point $z=(z_1,0)$ satisfying these conditions. Then
$$
\dist(K(z),\ETB)\le\|z-\ETB\|
+\|z-c_z\|+r_z<\delta+2r_z<\delta+\ve/4<\ve/2
$$
so that, by \rf{EP-1},
$$
\gamma_y\cap K(z)^{\cl}=\emp.
$$
\par On the other hand, by part (ii) of Statement B, there exists a point $x=(x_1,0)$ such that\medskip
\par (a). $z'\prec x$ for all $z'\in K(z)^{\cl}\cap\TB$;
\par (b). there exists a path$\gamma_x$ connecting $x$ to $B$ in $\Omega$ such that $\gamma_x\sm\{x\}\subset \GB\sm K(z)^{\cl}$.\medskip
\par Thus both connected components of $\TB\sm K(z)^{\cl}$ are $B$-accessible which contradicts Lemma \reff{T-KZ}.  This contradiction shows that Statement B is wrong in both cases proving that $\TBM\ne\emp$.
\par In the same way we show that the points of $\TB$ which are close enough to the point $b_{\TB}$ belong to $\TBP$ proving that $\TBP\ne\emp$.
\par The proof of the lemma is complete.\bx
\bigskip
\par We are in a position to finish the proof of Theorem \reff{N-SQ}.
\par {\it Proof of Theorem \reff{N-SQ}.} Under Assumption \reff{MA-KZ} the sets $\TBP$ and $\TBM$ are a partition of $\TB$. Clearly, $\TB$ is a connected topological space in induced Euclidean topology. But $\TBP$ and $\TBM$ are {\it non-empty and open subsets} of $\TB$ in this topology, see Lemma \reff{A-KZ} and Lemma \reff{N-EM}. This contradicts the connectedness of $\TB$.
\par Thus Assumption \reff{MA-KZ} is not true which easily implies the statement of Theorem \reff{N-SQ}. In fact, if there exists $z\in\TB$ such that $B\in K(z)$, then we put
$Q:=K(z)$. Since $z\in K(z)^{\cl}\cap\SCL\cap\Omega$, see \rf{Z-KSZ}, condition (i) of the theorem is satisfied. Furthermore, the first option of part (ii) of this theorem (i.e., the requirement $B\in K(z)$) holds, and the proof
in this case is complete.
\par Suppose that $B\notin K(z)$ for every $z\in\TB$. Since Assumption \reff{MA-KZ} is not true, there exists $z\in\TB$  such that part (ii) of Assumption \reff{MA-KZ} does not hold. This means that
\bel{GG-1}
\forall z'\in\TB, \forall~\text{path}~\gamma ~\text{joining}~z'~\text{to}~B,~ \gamma\sm\{z'\}\subset\GB,~\text{we have}~
\gamma\cap K(z)^{\cl}\ne\emp.
\ee
\par We again put $Q:= K(z)$. Since part (i) of Theorem \reff{N-SQ} is satisfied and, by the assumption, $B\notin Q=K(z)$, it remains to prove the statement \rf{Q-R}. This statement is equivalent to the following:
\bel{F-T}
\text{For every}~a\in\SB~\text{and every path}~\gamma~\text{joining}~a~\text{to}~B~\text{in}~\Omega~
\text{we have}~
\gamma\cap K(z)^{\cl}\ne\emp.
\ee
\par Prove this fact by representing $\gamma$ in a parametric form, i.e., as a graph of a continuous mapping $\Gamma:[0,1]\to \Omega$ such that $\Gamma(0)=a$ and $\Gamma(1)=B$. Let $a':=\Gamma(t_{\max})$ where
$$
t_{\max}:=\max\{t\in[0,1]:\gamma(t)\in\DS\}.
$$
Since $a\in\SB$ and $B\notin\SCL$, the point $a'$ is well defined. By $\gamma'$ we denote the arc of $\gamma$ from $a'$ to $B$. By definition of $t_{\max}$,
\bel{Y-G}
\gamma'\sm\{a'\}\cap \SCL=\emp,
\ee
so that, by Lemma \reff{CM-X} and Definition \reff{DF-GTB}, $a'\in\TB$ and $\gamma'\sm\{a'\}\subset\GB$. Then, by \rf{GG-1}, $\gamma'\cap K(z)^{\cl}\ne\emp$ proving \rf{F-T}.
\par The proof of ``The Square Separation Theorem'' \reff{N-SQ} is complete.\bx\medskip
\begin{remark} \lbl{L-SO} {\em Note that we are able to prove the following slight improvement of the statement \rf{F-T}:
\bel{F-TNEW}
\forall~a\in\SCL\cap\Omega~\text{and}~\forall~
\text{path}~\gamma~\text{joining}~a~\text{to}~B~\text{in}~\Omega~
\text{we have}~
\gamma\cap K(z)^{\cl}\ne\emp.
\ee
\par In fact, let $a\in\DS\cap\Omega$. If $\gamma\cap\SB\ne\emp$, then the proof of \rf{F-TNEW} is reduced to the previous case of $a\in\SB$ proven below. If $\gamma\cap\SB=\emp$, then we can put $a'=a$ in \rf{Y-G} so that this equality will be satisfied.
\par This enables us to modify the statement \rf{Q-R} of Theorem \reff{N-SQ} as follows:
$$
\SCL\cap\Omega~~\text{and}~~B~~\text{lie in different connected components of}~~\Omega\setminus Q^{\cl}.\rbf
$$}
\end{remark}
\par We finish the section with two remarks which present certain additional useful properties of the square $Q$ from formulation of Theorem \reff{N-SQ}.
\begin{remark} \lbl{L-SO2} {\em We notice that the square $Q$ from Theorem \reff{N-SQ} coincides with a square $K(z)$ for some $z\in\TB$. Applying part (d) and part (e) of Lemma \reff{PR-KZ} to $K(z)=Q$ we conclude that $Q$ has the following properties:\medskip
\par (i). The line segment $[\cs,c_Q]\subset \Omega$;
\smallskip
\par (ii). For every point $u\in Q^{\cl}\cap\Omega$ there exists a path $\gamma$ which joins $u$ to $B$ in $\Omega$ such that $(\gamma\sm\{u\})\cap\SCL=\emp$.\rbx}
\end{remark}
\par Our next remark relates to a certain improvement of part (ii) of ``The Square Separation Theorem'' \reff{N-SQ}, see Remark \reff{IMP-D} below. This improvement is based on the following
\begin{lemma} \lbl{SI-1} Let $K$ be a square and let $x,y\in\Omega\sm K^{\cl}$. Suppose there exists a polygonal path $\gamma$ which joins $x$ to $y$ in $\Omega$ such that $\gamma\cap K=\emp$.
\par Then there exists a polygonal path $\tgm$ joining $x$ to $y$ in $\Omega$ such that $\tgm\cap K^{\cl}=\emp$.
\end{lemma}
\par {\it Proof.} We will obtain the path $\tgm$ by a slight modification of $\gamma$ around the set $H:=\gamma\cap\,\partial K$. Since $\gamma$ is a polygonal path in $\Omega$, the set $H$ can be represented as a union of a {\it finite} number of pairwise disjoint subarcs of $\gamma$ lying on $\partial K$. In other words,
$$
H=\gamma\cap K^{\cl}=\bigcup_{i=1}^m\gamma_i
$$
where each $\gamma_i$ is either a subarc of $\gamma$ or a point of $\gamma$, and $\gamma_i\cap\gamma_j=\emp$, $1\le i,j\le m,$ $i\ne j$.
\par Let us represent $\gamma$ as a graph of a continuous mapping $\Gamma:[0,1]\to\Omega$ such that $\Gamma(0)=x$ and $\Gamma(1)=y$. Then each $\gamma_i$ is the graph of the mapping $\Gamma:[a_i,b_i]\to\Omega$ where
$$
0\le a_i\le b_i\le 1.
$$
Since the arcs $\gamma_i$ are disjoint, the line segments $[a_i,b_i], i=1,...,m,$ are disjoint as well.
\par Let $A_i:=\Gamma(a_i)$ and $B_i:=\Gamma(b_i)$ be the beginning and the end of the arc $\gamma_i$ respectively.
Let
$$
\ve:=\min_{1\le i,j\le m,\,i\ne j}
\{\dist(\gamma,\DO),\dist(\gamma_i, \gamma_j)\}.
$$
Then $[\gamma_i]_{\ve}\subset\Omega$, $1\le i\le m,$ and
$$
[\gamma_i]_{\ve}\cap [\gamma_j]_{\ve}=\emp,~~~i,j=1,...,m,~~i\ne j.
$$
(Recall that $[\,\cdot\,]_\ve$ denotes the $\ve$-neighborhood of a set.) Clearly, the set
$$
T_i:=[\gamma_i]_{\ve}\sm K^{\cl}
$$
is a {\it connected} open subset of $\Omega$.
\par Let $\gamma_i^{(p)}$ be the arc of $\gamma$ joining $x$ to $A_i$, and let $\gamma_i^{(f)}$ be the arc of $\gamma$ joining   $B_i$ to $y$. Since $\gamma$ is a continuous curve and $x,y\notin K^{\cl}$, there are exist points $\tA_i\in\gamma_i^{(p)}\cap T_i$ and $\tB_i\in\gamma_i^{(f)}\cap T_i$. Since $T_i$ is a connected subset of $\Omega$, there exists a polygonal path $\tgm_i$ joining $\tA_i$ to $\tB_i$ in $T_i$.
\par Now we replace the arc $\gamma_i$ by $\tgm_i$ for each $i\in\{1,...,m\}$. As a result we obtain a new polygonal path $\tgm$ which connects $x$ to $y$ in $\Omega$ and has no common points with $K^{\cl}$.\bx\medskip
\begin{remark}\lbl{IMP-D} {\em Lemma \reff{SI-1} and Remark \reff{L-SO} enable us to make further improvement of part (ii) of ``The Square Separation Theorem'' \reff{N-SQ}:\medskip
\par {\it (ii$'$). Either  $B\in Q^{\cl}$ or
\bel{Q-RN}
\SCL\cap\Omega~~\text{and}~~B~~\text{lie in different connected components of}~~\Omega\setminus Q.
\ee
}
\par Thus $\gamma\cap Q\ne\emp$ for every $z\in\SCL\cap\Omega$ and every path $\gamma$ which joins $z$ to $B$ in $\Omega$.
\rbx}
\end{remark}
\SECT{3. Proof of ``The Wide Path Theorem''}{3}
\setcounter{equation}{0}
\addtocontents{toc}{3. Proof of ``The Wide Path Theorem''.\hfill \thepage\par\VST}
\indent\par Basing on ``The Square Separation Theorem'' \reff{N-SQ} given $\x,\y\in\Omega$ we construct ``The Wide Path'' $\WPT^{(\x,\y)}$, see Theorem \reff{W-PATH}, as follows.
\par Let
$$
S_1:=S(\x,\dist(\x,\DO)).
$$
Thus $S_1$ is the maximal (with respect to inclusion) square in $\Omega$ centered at $\x$. If $\y\in S_1^{\cl}$, then we put $k=1$ and stop. If $\y\in\Omega\sm S_1^{\cl}$, we apply  Theorem \reff{N-SQ} to $\SB:=S_1$ and $B:=\y$. By this theorem, there exist a square $S_2\subset\Omega\sm S_1^{\cl}$ such that
$$
S_1^{\cl}\cap S_2^{\cl}\cap\Omega \ne\emptyset,
$$
and {\it either $\y\in S_2^{\cl}$} or
$$
{ S_2~~\text{\it{and}}~~\y~~\text{{\it lie in distinct connected components of}}~~\Omega\sm S_2^{\cl}.}
$$
If $\y\in S_2^{\cl}$, then we put $k=2$ and stop. If not, using ``The Square Separation Theorem'' we construct a square $S_3$, etc.
\par Continuing this procedure we obtain a sequence $\{S_1,S_2,...,S_m,...\}$ of squares (finite or infinite). Let $k$ be the number of its elements; thus $k=\infty$ whenever the sequence is infinite.
\par In the next lemma we present  main properties of the squares $S_i, i=1,2,...$. Let $c_i$ and $r_i$ be the center and ``radius'' of the square $S_i$ respectively, i.e.,
$$
S_i=S(c_i,r_i),~~~i=1,2,...~.
$$
\begin{lemma}\lbl{MAIN-WP} (a). $\x\in S_1$ and $\y\in S_{k}^{\cl}$ provided $k<\infty$. Furthermore, if $1<k<\infty$, then
$\dist(\y,S_{k-1})=\diam S_k$\,;\medskip
\par (b). $S_i\subset\Omega$ and $S_i^{\cl}\cap\DO\ne\emp$ for every $1\le i<k$\,;
\par (c). For all $i, 1\le i< k$, we have $S_i^{\cl}\cap S_{i+1}=\emp$, but
$S_i^{\cl}\cap S_{i+1}^{\cl}\cap\Omega \ne\emptyset$.
\par Furthermore,
\bel{CI-B}
[c_i,c_{i+1}]\subset\Omega,~~~1\le i< k\,;
\ee
\par (d). Let $1\le i<k-1$ and let $a\in S_i^{\cl}\cap \Omega$. Then $\gamma\cap S_{i+1}\ne\emp$ for any path $\gamma$ connecting $a$ to $\y$ in $\Omega$\,;
\par (e). For every $1\le i<k$ and every  $z\in S_{i+1}^{\cl}\cap\Omega$ there exists a path $\gamma$ joining $z$ to $\y$ in $\Omega$ such that $(\gamma\sm\{z\})\cap S_i^{\cl}=\emp$.
\end{lemma}
\par {\it Proof.} Parts (a) and (b) follow from the construction of the squares $\{S_i\}$ and the proof of ``The Square Separation Theorem'' \reff{N-SQ}; see part (b) of  Lemma \reff{PR-KZ}. Since the unique requirement to the square $S_k$ is that $S_k\ni \y$ and $S_k$ touches $S_{k-1}$, one can choose $S_k$ in such a way that
$\dist(\y,S_{k-1})=\diam S_k$.
\par Note that part (c) of the lemma directly follows from the construction of the squares $\{S_i\}$, part (i) of Theorem \reff{N-SQ} and \rf{C-KSZ}. In turn, part (d) and part (e) are consequences of \rf{Q-RN}, see Remark \reff{IMP-D}, and part (ii) of Remark \reff{L-SO2} respectively.\bx
\par In the next four lemmas we present additional properties of the squares $\{S_1,S_2,...\}$ which we need for the proofs of Theorems \reff{W-PATH} and \reff{WP-EXT}.
\begin{lemma}\lbl{AD-1} (i). Let $k>1$ and let $1\le i<k-1$. Let $a\in S_i^{\cl}\cap \Omega$ and let $\gamma$ be a path joining $a$ to $\y$ in $\Omega$. Then $\gamma\cap S_j\ne\emp$ for every $j, i<j<k$;
\par (ii). $S_i^{\cl}\cap S_j=\emp$ for all $i,j\ge 1$, $i\ne j$.
\end{lemma}
\par {\it Proof.} (i). See Figure 13 for an example of a path $\gamma$ joining in $\Omega$ a point $a\in S_2^{\cl}\cap\Omega$ to $\y$. 
\begin{figure}[h]
\center{\includegraphics[scale=1.0]{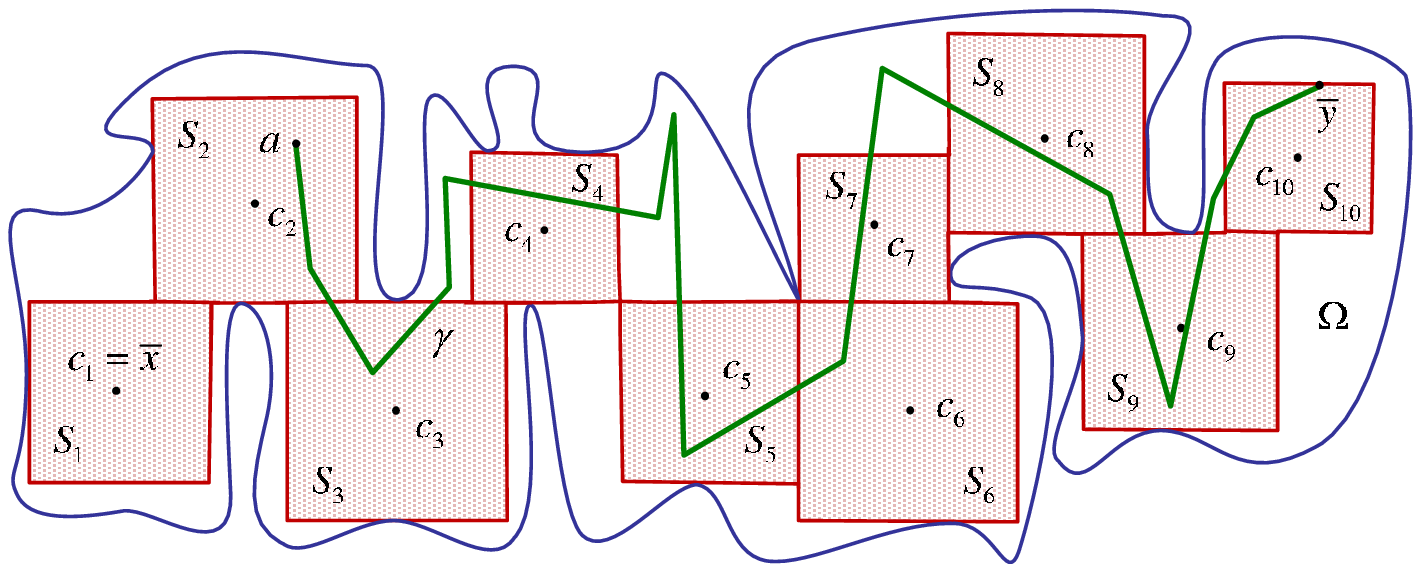}}
\caption{A path $\gamma$ connects a point $a\in S_2$ to $\y$ in $\Omega$.}
\end{figure}
\par We prove property (i) by induction on $j$. For $j=i+1$ it follows from part (d) of Lemma \reff{MAIN-WP}. Suppose that $\gamma\cap S_j\ne\emp$ for some $j>i+1$. Prove that $\gamma\cap S_{j+1}\ne\emp$ as well.
\par In fact, let $b\in\gamma\cap S_j$ and let $\gamma_b$ be the arc of $\gamma$ from $b$ to $\y$. Since $b\in S_j^{\cl}\cap\Omega$, by property (d) of Lemma \reff{MAIN-WP}, $\gamma_b\cap S_{j+1}\ne\emp$, proving
the statement (i) of the lemma.
\par (ii). Let $i<j$. Prove this statement by induction on $j$. By part (c) of Lemma \reff{MAIN-WP}, $S_i^{\cl}\cap S_{i+1}=\emp$.
\par Suppose that $S_i^{\cl}\cap S_{j}=\emp$ for some $j>i+1$, and prove that $S_i^{\cl}\cap S_{j+1}=\emp$ as well. Assume that it is not true, i.e., that there exists $z\in S_i^{\cl}\cap S_{j+1}$. Since $z\in S_{j+1}$, by part (e) of Lemma \reff{MAIN-WP}, there exists a path $\gamma_1$ joining $z$ to $\y$ in $\Omega$ such that
$\gamma_1\cap S_j^{\cl}=\emp.$
\par On the other hand, $z\in S_i^{\cl}\cap S_{j+1}\subset S_i^{\cl}\cap\Omega$ so that, by part (i) of the present lemma, $\gamma_1\cap S_j^{\cl}\ne\emp$, a contradiction which proves part (ii) for $i<j$.
\par Let $j<i$. As we have proved, in this case $S_j^{\cl}\cap S_i=\emp$ so that $S_j\cap S_i=\emp$ as well. Hence $S_i^{\cl}\cap S_j=\emp$, and the proof of the lemma is complete.\bx
\begin{lemma} $k<\infty$, i.e., $\{S_1,S_2,...\}$ is a finite family of squares.
\end{lemma}
\par {\it Proof.} Let $\gamma$ be a path connecting $\x$ to $\y$ in $\Omega$. Since $\x=c_{S_1}\in S_1$, by part (i) of Lemma \reff{AD-1}, $\gamma\cap S_i\ne\emp$ for every $1\le i<k$.
\par Note that the path $\gamma$ is a compact subset of $\Omega$ so that $\ve:=\dist(\gamma,\DO)>0$. Prove that for each square $S_i$, $i\ge 1$, we have $\diam S_i\ge \ve$.
\par In fact, let $a\in \gamma\cap S_i$. Then $\dist(a,\DO)\ge\dist(\gamma,\DO)=\ve$.
\par Recall that $S_i=S(c_i,r_i)$, $i=1,2,...\,$. By part (b) of Lemma \reff{MAIN-WP}, $S_i\subset\Omega$ and $S_i^{\cl}\cap\DO\ne\emp$, so that $r_i=\dist(c_i,\DO)$.
Hence,
$$
\ve\le \dist(a,\DO)\le \|a-c_i\|+\dist(c_i,\DO)\le r_i+r_i=\diam S_i.
$$
\par By part (ii) of Lemma \reff{AD-1}, the squares of the family $\Sc=\{S_i:1\le i<k\}$ are non-overlapping. Since the diameter of each square from $\Sc$ is at least $\ve$, the domain $\Omega$ contains at most $|\Omega|/\ve^2$ squares from this family. Since $\Omega$ is bounded, this number is finite, and the proof is complete.\bx
\medskip
\par The next lemma provides a certain improvement of Lemma \reff{AD-1}.
\begin{lemma}\lbl{AD-1-2} (i). Let $k>1$ and let\,
$1\le i<m-1\le k-1$. Let $a\in S_i^{\cl}\cap\Omega$,
$b\in S_m^{\cl}\cap\Omega$, and let
$\gamma$ be a path joining $a$ to $b$ in $\Omega$.
Then $\gamma\cap S_j\ne\emp$ for every $j,\,i<j< m$; \medskip
\par (ii). $S_i^{\cl}\cap S_j^{\cl}\cap\Omega=\emp$ for every $1\le i,j\le k$ such that $|i-j|>1$;
\medskip
\par (iii). Let $1\le i\le m\le k$ and let $a\in S_i^{\cl}\cap\Omega$, $b\in S_m^{\cl}\cap\Omega$. There exists a simple path $\gamma$ which joins $a$ to $b$ in $\Omega$ such that
\bel{G-UJ}
\gamma\subset \bigcup_{j=i}^m\,\left(S_j^{\cl}\,\ism\,\Omega\right).
\ee
\par Furthermore, $\gamma\cap S_j^{\cl}=\emp$ provided $j>m+1$ or $j<i-1$, and $(\gamma\sm\{a\})\cap S_{i-1}^{\cl}=\emp$ and  $(\gamma\sm\{b\})\cap S_{m+1}^{\cl}=\emp$.
\end{lemma}
\par {\it Proof.} (i). We prove the statement (i) by induction on $n=m-j$, $1\le n<m-i$. Let $n=1$, i.e., $j=m-1$. Prove that $\gamma\cap S_{m-1}\ne\emp$.
\par Since  $b\in S_m^{\cl}\cap\Omega$, by property (e) of Lemma \reff{MAIN-WP}, there exists a path $\gamma_1$ joining $b$ to $\y$ in $\Omega$ such that
\bel{G1-B}
(\gamma_1\sm\{b\})\cap S_{m-1}^{\cl}=\emp.
\ee
\par Let $\tgm:=\gamma\cup\gamma_1$. Then $\tgm$ is a path which connects $a$ to $\y$ in $\Omega$ so that, by part (i) of Lemma \reff{AD-1}, $\tgm\cap S_{m-1}\ne\emp$.
\par Recall that $b\in S_m^{\cl}$. Since $S_m^{\cl}\cap S_{m-1}\ne\emp$, see part (ii) of Lemma \reff{AD-1}, $b\notin S_{m-1}$. Combining this with \rf{G1-B} we conclude that $\gamma_1\cap S_{m-1}=\emp$. Since
$$
\tgm\cap S_{m-1}=(\gamma\cup\gamma_1)\cap S_{m-1}\ne\emp,
$$
we obtain that $\gamma\cap S_{m-1}\ne\emp$.\medskip
\par Now given $j=m-n+1$ suppose that $\gamma\cap S_j\ne\emp$. Prove that $\gamma\cap S_{j-1}\ne\emp$ as well.
\par We follow the same scheme as for the case $j=m-1$. Let $\tb\in\gamma\cap S_j$. Since
$\tb\in S_j\subset S_j^{\cl}\cap\Omega$, by part (e) of
Lemma \reff{MAIN-WP}, there exists a path $\gamma'$ which joins $\tb$ to $\y$ in $\Omega$ such that
\bel{GP-1}
(\gamma'\sm\{\tb\})\cap S_{j-1}^{\cl}=\emp.
\ee
\par Let $\gamma''$ be the arc of $\gamma$ from $a$ to $\tb$. Then the path $\bg:=\gamma'\cup\gamma''$ joins $a$ to $\y$ in $\Omega$ so that, by part (i) of Lemma \reff{AD-1}, $\bg\cap S_{j-1}\ne\emp$.
\par Since $\tb\in S_j$ and $S_j\cap S_{j-1}=\emp$, see part (ii) of Lemma \reff{AD-1}, we conclude that $b\notin S_{j-1}$. This and \rf{GP-1} imply that $\gamma'\cap S_{j-1}^{\cl}=\emp$. Since
$$
\bg\cap S_{j-1}=(\gamma'\cup\gamma'')\ne\emp,
$$
we conclude that $\gamma''\cap S_{j-1}^{\cl}\ne\emp$. But $\gamma''$ is a subarc of $\gamma$ so that $\gamma\cap S_{j-1}^{\cl}\ne\emp$ proving part (i) of the lemma. \medskip
\par (ii). Suppose that $S_i^{\cl}\cap S_j^{\cl}\cap\Omega\ne\emp$ for some $1\le i<j\le k$ such that $i+1<j$. Let $z\in S_i^{\cl}\cap S_j^{\cl}\cap\Omega$. Since $S_i\cap S_j=\emp$, the point $z\in \partial S_i\cap
\partial S_j$.
\par Let $\gamma:=[c_i,z]\cup[z,c_j]$. (Recall that $c_i$ is the center of $S_i$.) Clearly,
$c_i\in S_i^{\cl}\cap\Omega$ and $c_j\in S_j^{\cl}\cap\Omega$ so that, by part (i) of the present lemma,
\bel{GM-S}
\gamma\cap S_{i+1}\ne\emp.
\ee
\par However $\gamma\sm \{z\}\subset S_i\cup S_j$. Since $S_i,S_{i+1}$ and $S_j$ are pairwise disjoint,  $S_{i+1}\cap (S_i\cup S_j)=\emp$, so that, by \rf{GM-S},
$z\in S_{i+1}$. Since $z\in S_i^{\cl}$, this implies
$S_i^{\cl}\cap S_{i+1}\ne\emp$ which contradicts part (ii) of Lemma \reff{AD-1}.
\medskip
\par (iii). Recall that $S_j=S(c_j,r_j)$. Let $\gamma_1:=[a,c_i]$ and let $\gamma_3:=[c_m,b]$.
(Whenever $a=c_i$ or $b=c_m$ we ignore $\gamma_1$ or $\gamma_3$ respectively.)
By $\gamma_2$ we denote a polygonal path with vertices in $c_i,c_{i+1},...,c_{m-1},c_m$. Then a path $\gamma:=\gamma_1\cup\gamma_2\cup\gamma_3$ connects $a$ to $b$.
\par Clearly, $\gamma_1=[a,c_i]\subset S_i^{\cl}\cap\Omega$ and $\gamma_1\sm\{a\}\subset S_i$. Also
$\gamma_3=[c_m,b]\subset S_m^{\cl}\cap\Omega$ and $\gamma_3\sm\{b\}\subset S_m$.
\par On the other hand, by property \rf{CI-B}, see part (c) of Lemma \reff{MAIN-WP},
$$
[c_j,c_{j+1}]\subset (S_j^{\cl}\cup S_{j+1}^{\cl})\cap\Omega
$$
so that $\gamma_2\subset \cup\{S_{j}^{\cl}\cap\Omega:i\le j\le m\}$. These properties of the paths $\gamma_i$, $i=1,2,3,$ prove \rf{G-UJ}.
\par The second statement of part(iii) immediately follows from the fact that the squares $\{S_j\}$ are pairwise disjoint and $S_{j_1}^{\cl}\cap S_{j_2}^{\cl}\cap\Omega=\emp$ whenever $|j_1-j_2|>1$. See
part (ii) of Lemma \reff{AD-1} and part (ii) of the present lemma.
\par The proof of the lemma is complete.\bx
\begin{lemma}\lbl{AD-3} Let $1<i<k$. Then the set~$\cup\{S_j:i<j\le k\}$ and the point $\y$
belong to the same connected component of\, $\Omega\sm S_i^{\cl}$.\smallskip
\par In turn, the set\, $\cup\{S_j:1\le j<i\}$ and the point $\x$ belong to another connected component of\, $\Omega\sm S_i^{\cl}$.
\end{lemma}
\par {\it Proof.} Let $a\in\cup\{S_j:i<j\le k\}$
so that $a\in S_j$ for some $i+1\le j\le k$. Since $\y\in S_k^{\cl}$, by part (iii) of Lemma \reff{AD-1-2}, there exists a path $\gamma$ which connects $a$ to $\y$ in $\Omega$ such that $\gamma\sm\{a\}\cap S_i^{\cl}=\emp$.
But $S_j\cap S_i^{\cl}=\emp$, see part (ii) of Lemma \reff{AD-1}, so that $\gamma\sm\{a\}\cap S_i^{\cl}=\emp$.
This  proves that $a$ and $\y$ belong to the same connected component of $\Omega\sm S_i^{\cl}$.
\par In the same way we show that every point
$$
b\in \cup\{S_j:1\le j<i\}
$$
belong to the same connected component of $\Omega\sm S_i^{\cl}$ as the point $\x$.
\par It remains to note that, by part(i) of Lemma \reff{AD-1-2}, $\gamma\cap S_i^{\cl}\ne\emp$ for every path $\gamma$ joining $\x$ to $\y$ in $\Omega$ so that $\x$ and $\y$ belong to {\it distinct} connected components of $\Omega\sm S_i^{\cl}$.
\par The proof of the lemma is complete.\bx
\bigskip
\par {\it Proof of ``The Wide Path Theorem'' \reff{W-PATH}.} The proof immediately follows from lemmas proven in this section. In fact, part (i) and part (ii) of Theorem \reff{W-PATH} follow from part (a) and part (c) of Lemma \reff{MAIN-WP} respectively, and part (iii) follows from Lemma \reff{AD-3}.
\par ``The Wide Path Theorem'' \reff{W-PATH} is completely proved.\bx
\SECT{4. Sobolev extension properties of ``The Wide Path''}{4}
\setcounter{equation}{0}
\addtocontents{toc}{4. Sobolev extension properties of ``The Wide Path''.\hfill \thepage\par\VST}
\indent\par {\bf 4.1. ``The arc diameter condition'' and the structure of ``The Wide Path''.}
\addtocontents{toc}{~~~~4.1. ``The arc diameter condition'' and the structure of ``The Wide Path''.\hfill \thepage\par}
In this section we prove Theorem \reff{WP-EXT} which states that given $\x,\y\in\Omega$ any ``Wide Path'' $\WPT^{(\x,\y)}$ joining $\x$ to $\y$ in $\Omega$, see \rf{WP-DEF}, has the Sobolev extension property provided the domain $\Omega$ has.
\par We recall that, by the Sobolev imbedding theorem, see e.g., \cite{M}, p. 73, every function $f\in \LM(\Omega)$, $p >2$, can be redefined, if necessary, in a set of Lebesgue measure zero so that it belongs to the space $C^{m-1}(\Omega)$. Thus, for $p>2$, we can identify each element $f\in \LM(\Omega)$ with its unique $C^{m-1}$-representative on $\Omega$. This will allow us to restrict our attention to the case of Sobolev $C^{m-1}$-functions.\medskip
\par In this section and in Sections 5 and 6 we assume that $\Omega$ is a {\it simply connected} bounded domain in $\RT$ satisfying the hypothesis of Theorem \reff{MAIN-NEC}:\medskip
\par {\it There exists a constant $\CE\ge 1$ such that} 
\bel{S-EXT}
\forall f\in \LMPO~~ \exists\, F\in\LMT~\text{\it such that}~
F|_\Omega=f~\text{\,\it{and}\,}~ \|F\|_{\LMT}\le\CE\|f\|_{\LMPO}\,.
\ee
\par In other words, we assume that $e_{m,p}(\Omega)\le \CE$, see \rf{E-MP}.\smallskip
\par The following well known property of Sobolev extension domains proven by  Gol'dshtein and Vodop'janov \cite{GV1} shows that every domain $\Omega\subset\RT$ satisfying \rf{S-EXT} is ``almost quasiconvex''. Here we present a slight improvement of this property given in \cite{GR}, Chapter 6, Theorems 2.5 and 2.8.
\begin{theorem}\lbl{T-CM} Let $p>2$, $m\in \N$, and let $\Omega$ be a domain in $\RT$ satisfying condition \rf{S-EXT}. 
\par Then for every $a,b\in\Omega$ there exists a path $\gamma$ which connects $a$ to $b$ in $\Omega$ such that
$\diam \gamma\le \eta\|a-b\|$.
\par Here $\eta$ is a positive constant such that the following inequality $\eta\le C(m,p)\,\CE$ holds.
\end{theorem}
\par Following \cite{GV1} we refer to this property as ``the arc diameter condition''.\smallskip
\par Theorem \reff{T-CM} enables us to prove an additional geometrical property of the family of squares $\{S_1,...,S_k\}$ defined in the previous section.
\par Consider two subsequent squares from this family, say $S_i$ and $S_{i+1}$, $1\le i<k$, such that $\#(S_i^{\cl}\cap S_{i+1}^{\cl})>1$. Since $S_i$ and $S_{i+1}$ are touching squares, intersection of their closures is a line segment which we denote by $[u_i,v_i]$:
\bel{UV-I}
[u_i,v_i]:=S_i^{\cl}\cap S_{i+1}^{\cl}\,.
\ee
\par Note that in this case
\bel{SGM-I}
(S_i^{\cl}\cup S_{i+1}^{\cl})^\circ=S_i\cup S_{i+1}\cup(u_i,v_i).
\ee
\begin{lemma}\lbl{H-MS} Let $\Omega$ be a simply connected bounded domain in $\RT$ satisfying condition \rf{S-EXT}. Let $\x,\y\in\Omega$ and let $\Sc_\Omega(\x,\y)=\{S_1,...,S_k\}$ be the sequence of squares constructed in Theorem \reff{W-PATH}.\smallskip
\par (i). Let $1\le i<k$ and let $S_i,S_{i+1}$ be two consecutive squares from this family such that $\#(S_i^{\cl}\cap S_{i+1}^{\cl})>1$. Then $(u_i,v_i)\subset\Omega$ and
$$
(S_i^{\cl}\cup S_{i+1}^{\cl})^\circ\subset\Omega\,;
$$
\par (ii). $\#(S_i^{\cl}\cap S_{i+2}^{\cl})\le 1$ for all $i, 1\le i\le k-2$, and
\bel{SIJ-L}
S_i^{\cl}\cap S_{j}^{\cl}=\emp~~~\text{if}~~~|i-j|>2,~1\le i,j\le k\,;
\ee
\par (iii). If $\#(S_i^{\cl}\cap S_{i+2}^{\cl})=1$ for some $i, 1\le i\le k$, then $S_i^{\cl}\cap S_{i+1}^{\cl}\cap S_{i+2}^{\cl}=\{a_{i+1}\}$ is a singleton. The point $a_{i+1}\in\DO$. This point is a common vertex of the squares $S_i,S_{i+1}$ and $S_{i+2}$ and belongs to the boundary of the set $S^{\cl}_i\cup S^{\cl}_{i+1}\cup S^{\cl}_{i+2}$.
\par See Figure 14. See also the squares $S_6, S_7$ and $S_8$ on Figure 1.
\begin{figure}[h]
\center{\includegraphics[scale=0.75]{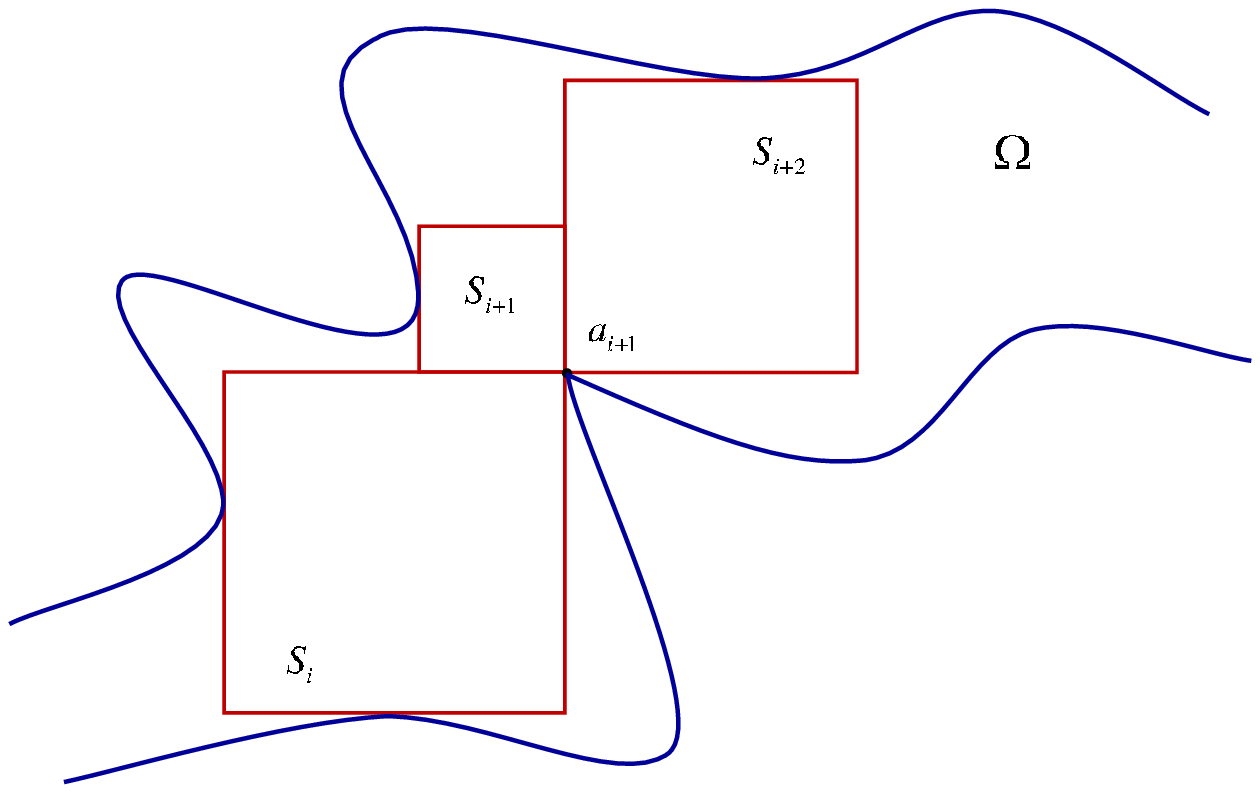}}
\caption{~}
\end{figure}
\end{lemma}
\par {\it Proof.} Let us prove part (i) of the lemma. Note that $S_i^{\cl}\cap S_{i+1}^{\cl}$ is a line segment because $S_i$ and $S_{i+1}$ are touching squares such that $\#(S_i^{\cl}\cap S_{i+1}^{\cl})>1$.
\par Prove that $(u_i,v_i)\subset\Omega$. In fact, $(u_i,v_i)\cap\Omega$ is an open set in the relative topology of the straight line passing through $u_i$ and $v_i$. By part (ii) of Theorem \reff{W-PATH}, this set is non-empty, so that $(u_i,v_i)\cap\Omega$ can be represented as a union of a finite or countable family $\Ic$ of pairwise disjoint open subintervals of $(u_i,v_i)$ with ends in $\DO$.
\par Let us show that this family contains precisely {\it one subinterval} of $(u_i,v_i)$, i.e., $\#\Ic=1$. Suppose that it is not true, i.e., that there exist two distinct line intervals from this family, say $I'=(x',y')$ and $I''=(x'',y'')$, $I'\ne I''$. Then $x',y',x'',y''\in\DO$ and $I'\cup I''\subset(u_i,v_i)\cap\Omega$. See Figure 15.
\medskip
\begin{figure}[h!]
\center{\includegraphics[scale=1]{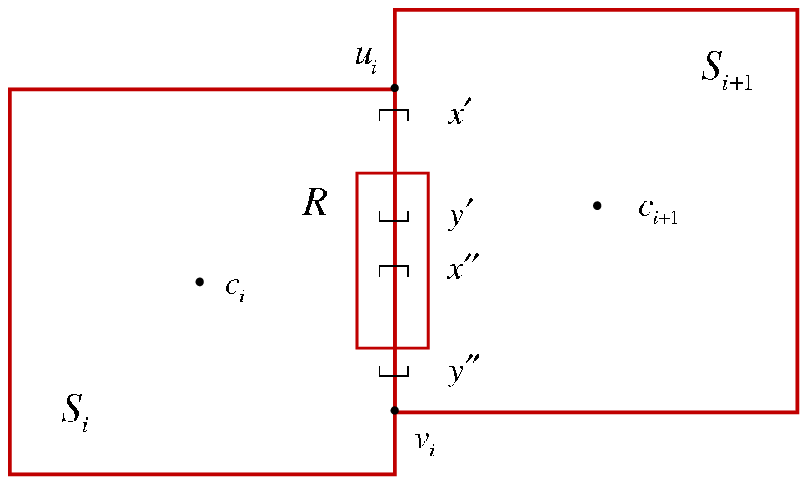}}
\caption{~}
\end{figure}
\medskip
\par We may assume that $y',x''\in(x',y'')$. Then there exists a rectangle $R$ with sides parallel to the coordinate axes and width small enough such that $y',x''\subset R\hspace{0.2mm}^{\circ}$ and $\partial R\subset\Omega$. See Figure 15. Since $\Omega$ is simply connected, $R\subset\Omega$ so that $y'\in\Omega$. But $y'\in\DO$, a contradiction.
\par Thus $\#\Ic=1$ so that
$$
(u_i,v_i)\cap\Omega=(z',z'')~~\text{for some}~~ z',z''\in[u_i,v_i].
$$
\par We may assume that
$$
\|z'-u_i\|<\|z''-u_i\|~~~\text{and}~~~ ~~~\|z''-v_i\|<\|z'-v_i\|.
$$
Prove that $z'=u_i$ and $z''=v_i$. Suppose that it is not true, and, for instance, $z'\ne u_i$. Then the line segment $[u_i,z']\subset\DO$.
\par Let $\tz:=(u_i+z')/2$. Then there exist sequences $\{a_j\}_{j=1}^\infty\subset S_i$ and  $\{b_j\}_{j=1}^\infty\subset S_{i+1}$ such that
$$
a_j\to \tz~~\text{and}~~b_j\to \tz~~\text{as}~~j\to\infty. $$
\par Since $[u_i,z']\subset\DO$, any path $\gamma$ joining $a_j$ to $b_j$ in $\Omega$ has the diameter at least $\|u_i-z'\|/8$ provided $a_j$ and $b_j$ are close enough to $\tz$. On the other hand, $\Omega$ satisfies condition \rf{S-EXT} so that, by Theorem \reff{T-CM}, the points $a_j$ and $b_j$  can be joined by a certain path $\gamma_j$ such that
$$
\diam \gamma_j\le\eta\|a_j-b_j\|.
$$
\par Hence,
$$
\|a_i-z'\|/8\le \diam \gamma_j\le\eta\|a_j-b_j\|\to 0 ~~\text{as}~~j\to\infty.
$$
Since $a_j\to \tz$, we have $\tz=z'$ so that $z'=u_i$, a contradiction. In the same fashion we prove that $z''=v_i$ so that $(u_i,v_i)=(z',z'')\subset\Omega$.
\par Finally, we obtain that
$$
(S_i^{\cl}\cup S_{i+1}^{\cl})^\circ=S_i\cup S_{i+1}\cup(u_i,v_i)\subset\Omega
$$
proving part (i) of the lemma.\medskip
\par Prove part (ii) and (iii). First prove that
\bel{F-DM}
\#(S_i^{\cl}\cap S_{j}^{\cl})\le 1~~~\text{provided}~~~|i-j|>1,~1\le i,j\le k.
\ee
\par Suppose that $1\le i<j\le k$ and $S_i^{\cl}\cap S_{j}^{\cl}\ne\emp$.
\par Since $S_i\cap S_{j}=\emp$, we have
$S_i^{\cl}\cap S_{j}^{\cl}=\partial S_i\cap \partial S_{j}$ so that
\bel{SCL-1}
S_i^{\cl}\cap S_{j}^{\cl}=[a,b]~~~\text{for some}~~~a,b\in\RT.
\ee
\par We know that
$S_i^{\cl}\cap S_{j}^{\cl}\cap\Omega=\emp$ whenever $|i-j|>1$, see part (ii) of Lemma \reff{AD-1-2}. Hence $[a,b]\subset\RT\sm\Omega$. On the other hand, by \rf{SCL-1}, $[a,b]\subset \Omega^{\cl}$ so that  $[a,b]\subset \DO$.
\par Let us assume that $\#(S_i^{\cl}\cap S_{j}^{\cl})>1$, i.e., that $a\ne b$. Let $z:=(a+b)/2$. Since $z\in S_i^{\cl}\cap S_{j}^{\cl}$, there exist sequences of points
\bel{ST-S}
\{s_n\}_{n=1}^\infty\subset S_i~~~\text{and}~~~
\{t_n\}_{n=1}^\infty\subset S_j~~~\text{such that}~~~s_n,t_n\to z~~~\text{as}~~~n\to\infty.
\ee
\par Since $\Omega$ satisfies condition \rf{S-EXT} and $s_n,t_n\to z$, by Theorem \reff{T-CM}, there exists a path $\gamma_n$ connecting $s_n$ to $t_n$ in $\Omega$ such that
\bel{G-FT}
\diam\gamma_n\le\eta\|s_n-t_n\|
\ee
provided $n>N$ where $N$ is big enough. We may also assume that $N$ is so big that
\bel{W-1}
\|z-s_n\|<\|a-b\|/(8\eta)~~~\text{and}~~~
\|z-t_n\|<\|a-b\|/(8\eta)~~~\text{for}~~~n>N.
\ee
\par Note that the straight line passing through $a$ and $b$ separates $s_n$ and $t_n$ and the path $\gamma_n$ does not cross the line segment $[a,b]$. Therefore
\bel{CR-2}
\diam\gamma_n\ge \tfrac12\|a-b\|-\tfrac18\|a-b\|=\tfrac38\|a-b\|.
\ee
On the other hand, by \rf{G-FT} and \rf{W-1},
$$
\diam\gamma_n\le\eta\|s_n-t_n\|\le \eta(\|z-s_n\|+\|z-t_n\|)\le 2\eta\,\frac{\|a-b\|}{8\eta}=\tfrac14 \|a-b\|.
$$
\par This inequality contradicts to inequality \rf{CR-2} proving \rf{F-DM}.\smallskip
\par Now suppose that
$$
S_i^{\cl}\cap S_{j}^{\cl}\ne\emp~~~\text{for some}~~~1\le i<j\le k,
$$
and prove that this condition is satisfied only for $j=i+2$.
\par In fact, by \rf{F-DM},
$S_i^{\cl}\cap S_{j}^{\cl}=\{a\}$
for some $a\in\DO\cap\partial S_i\cap \partial S_{j}$. Prove that
\bel{A-SL}
a\in S_\ell^{\cl}~~~\text{for every}~~~i\le\ell\le j.
\ee
\par As above, by $\{s_n\}_{n=1}^\infty$ and  $\{t_n\}_{n=1}^\infty$ we denote the sequences of points satisfying \rf{ST-S}, and by $\gamma_n$ we denote a path joining $s_n$ to $t_n$ in $\Omega$ such that \rf{G-FT} holds. Then, by part (i) of Lemma \reff{AD-1-2},
$$
\gamma_n\cap S_\ell\ne\emp~~~
\text{for every}~~~\ell,~ i\le \ell\le j.
$$
\par Let $b^{(n)}_\ell\in \gamma_n\cap S_\ell$, $i<\ell<j$. We also put $b^{(n)}_i:=s_n$ and
$b^{(n)}_j:=t_n$. Then, by \rf{G-FT},
$$
\|b^{(n)}_\ell-s_n\|\le \diam\gamma_n\le\eta\|s_n-t_n\|.
$$
Since $\|s_n-t_n\|\to 0$ and $s_n\to a$ as $n\to\infty$, we conclude that $b^{(n)}_\ell\to a$ for every $\ell, i\le\ell\le j$. Hence, $a\in S^{\cl}_\ell$ proving \rf{A-SL}.
\par Thus, by \rf{F-DM} and \rf{A-SL}, if $i+2\le j\le k$ and $S_i^{\cl}\cap S_{j}^{\cl}\ne\emp$, then there exists a point $a\in\partial S_i$ such that
\bel{FN-L}
S_i^{\cl}\cap S_{\ell}^{\cl}=\{a\}~~~\text{for all}~~~ \ell,~i\le\ell\le j.
\ee
\par Since $\{S_\ell:i\le \ell\le j\}$ are pairwise disjoint squares, this property easily implies the required restriction $j=i+2$. In fact, since $S_i\cap S_{\ell}=\emp$ and
$$
S_i^{\cl}\cap S_{\ell}^{\cl}=\{a\},
$$
the point $a$ is a {\it vertex} of the square  $S_\ell$ for every  $\ell,~i\le\ell\le j$. In particular, $a$ is a common vertex of $S_i$ and $S_{i+2}$. Note that, by \rf{FN-L}, $a\in S_{i+1}^{\cl}$. Since $S_i,S_{i+1},S_{i+2}$ are pairwise disjoint squares, this  implies that the point $a_{i+1}=a$ is a vertex of the square $S_{i+1}$ as well.
\par By part (ii) of Lemma \reff{AD-1-2}, $S^{\cl}_i\cap S^{\cl}_{i+2}\cap \Omega=\emp$ so that $a\in\DO$. It is also clear that $a$ is a boundary point of the set $S^{\cl}_i\cup S^{\cl}_{i+1}\cup S^{\cl}_{i+2}$.
See Figure 16.
\begin{figure}[h]
\center{\includegraphics[scale=0.8]{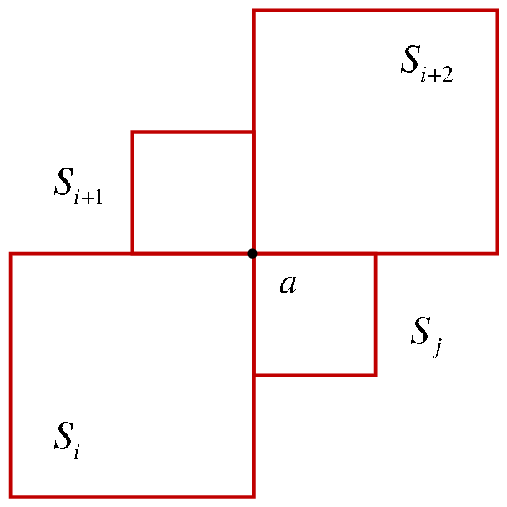}}
\caption{~}
\end{figure}
\par But, $a$ is also a vertex of the square $S_j$. Since  $S_i,S_{i+1},S_{i+2}$ and $S_j$ are pairwise disjoint squares, and $a$ is a common vertex of these squares, the intersection of $S_i^{\cl}$ and $S_j^{\cl}$ is a line segment (of positive length). See Figure 16.
\par Thus $\#(S_i^{\cl}\cap S_j^{\cl})>1$ whenever $j>i+2$ which contradicts \rf{F-DM}. Hence, $S_i^{\cl}\cap S_j^{\cl}=\emp$ provided $|i-j|>2$.
\par The proof of the lemma is complete.\medskip\bx
\par Part (iii) of Lemma \reff{H-MS} motivates us to introduce the following
\begin{definition}\lbl{R-SQ-P} {\em Let $1\le i\le k-1$ and let $a_{i+1}\in\DO$ be a {\it common vertex} of the squares $S^{\cl}_i,S^{\cl}_{i+1}$ and $S^{\cl}_{i+2}$, i.e.,
$$
\{a_{i+1}\}=S^{\cl}_i\cap S^{\cl}_{i+1}\cap S^{\cl}_{i+2}.
$$
\par We refer to the point $a_{i+1}$ as a {\it rotation point} of ``The Wide Path'' $\Wc=\WPT^{(\x,\y)}$, and to the square $S_{i+1}$ as a {\it rotation square} of $\Wc$.}
\end{definition}
\par See Figure 14. Another example is given on Figure 17. 
\begin{figure}[h]
\center{\includegraphics[scale=0.75]{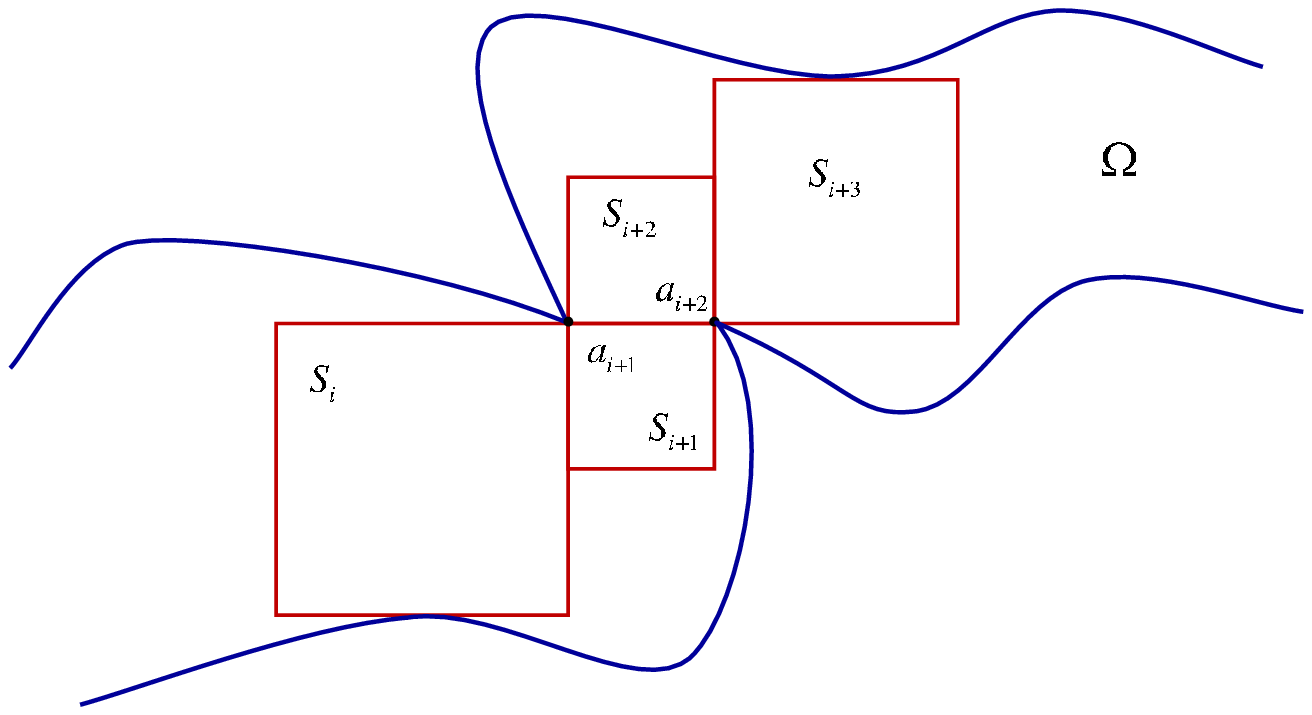}}
\caption{~}
\end{figure}
\par Here $\{a_{i+1}\}$ and $\{a_{i+2}\}$ are the rotation points corresponding to the rotation squares $S_{i+1}$ and  $S_{i+2}$. Note that rotation points and rotation squares play an important in construction of ``The Narrow path''. See Section 5.\bigskip
\par {\bf 4.2. Subhyperbolic properties of elementary squarish domains.}
\addtocontents{toc}{~~~~4.2. Subhyperbolic properties of elementary squarish domains.\hfill \thepage\par}
We will be needed several auxiliary results related to subhyperbolic properties of domains in $\RT$ consisting of a ``small number'' of open squares. We refer to such sets as ``elementary squarish domains''.
\begin{lemma}\lbl{Q-DA} Let $Q$ be a square in $\RT$ and let $a,b\in Q^{\cl}$. Then there exists a path $\gamma_{ab}$ joining $a$ to $b$ and consisting of at most two edges such that $\gamma_{ab}\sm\{a,b\}\subset Q$ and for every $\alpha\in(0,1]$ the following inequality
\bel{LA-Q}
\len_{\alpha,Q}(\gamma_{ab})\le \tfrac{3}{\alpha}\,\|a-b\|^\alpha
\ee
holds. See \rf{SH-LN}.
\end{lemma}
\par {\it Proof.} Let $a=(a_1,a_2)$, $b=(b_1,b_2)$
and let $Q=(u_1,u_2)\times(v_1,v_2)$. Suppose that $|a_2-b_2|\le |a_1-b_1|=\|a-b\|$. Since $a_2,b_2\in [v_1,v_2]$ and
$$
|v_1-v_2|=\diam Q\ge\|a-b\|\ge |a_2-b_2|
$$
there exists a line segment $[s_1,s_2]$ such that
$$
a_2,b_2\in [s_1,s_2]\subset[v_1,v_2]~~~\text{and}~~~
|s_1-s_2|=|a_1-b_1|=\|a-b\|.
$$
\par Let $Q_{ab}:=(a_1,b_1)\times (s_1,s_2)$. Then $Q_{ab}\subset Q$, $a,b\in\partial Q_{ab}$ and $\diam Q_{ab}=\|a-b\|$. Let $Q_{ab}=S(c,r)$, i.e., $c$ is the center of $Q_{ab}$ and $r=\tfrac12\|a-b\|$ is its ``radius'', and let
$$
\gamma_{ab}:=[a,c]\cup[c,b].
$$
Clearly, $\gamma_{ab}$ is a two edges path connecting $a$ to $b$ such that $\gamma_{ab}\sm\{a,b\}\subset Q$.
\par Prove inequality \rf{LA-Q}. By definition \rf{SH-LN},  
\be
\len_{\alpha,Q}(\gamma_{ab})&:=&
\intl_{\gamma_{ab}}
\dist(z,\partial Q)^{\alpha-1}\,ds(z)\nn\\
&\le&
\intl_{[a,c]}
\dist(z,\partial Q_{ab})^{\alpha-1}\,ds(z)+
\intl_{[c,b]}
\dist(z,\partial Q_{ab})^{\alpha-1}\,ds(z)\nn\\
&=&I_1+I_2.\nn
\ee
\par Note that $\dist(z,\partial Q_{ab})=\|z-a\|$ for every $z\in[a,c]$. Hence,
$$
I_1:=
\intl_{[a,c]}
\dist(z,\partial Q_{ab})^{\alpha-1}\,ds(z)=\intl_{[a,c]}
\|z-a\|^{\alpha-1}\,ds(z)=\intl_0^1(\|a-c\|\,t)^{\alpha-1}
\|a-c\|_2\,dt
$$
where $\|\cdot\|_2$ denotes the Euclidean norm in $\RT$.
\par Recall that $a\in\partial Q_{ab}$ so that
$$
\|a-c\|=r=\tfrac12\diam Q_{ab}=\tfrac12\,\|a-b\|.
$$
We obtain:
$$
I_1\le r^{\alpha-1}(\sqrt{2}r)\,\intl_0^1\,t^{\alpha-1}\,dt=
\tfrac{\sqrt{2}}{\alpha}\,r^{\alpha}.
$$
In the same way we prove that
$$
I_2:=
\intl_{[c,b]}
\dist(z,\partial Q_{ab})^{\alpha-1}\,ds(z)
\le \tfrac{\sqrt{2}}{\alpha}\,r^{\alpha}.
$$
Hence,
$$
\len_{\alpha,Q}(\gamma_{ab})\le I_1+I_2\le
\tfrac{2\sqrt{2}}{\alpha}\,r^{\alpha}=
\tfrac{2\sqrt{2}}{\alpha\,2^\alpha}\,\|a-b\|^{\alpha}
\le\tfrac{3}{\alpha}\,\|a-b\|^{\alpha}.
$$
\par The proof of the lemma is complete.\bx
\begin{lemma}\lbl{2-SQ} Let $Q_1$ and $Q_2$ be squares in $\RT$ such that $\#(Q_1^{\cl}\cap Q_2^{\cl})>1$, and let
$$
G:=(Q_1^{\cl}\cup Q_2^{\cl})^\circ.
$$
\par Then for every $a,b\in G^{\cl}$, $a\ne b$, there exists a path $\gamma_{ab}(G)$ consisting of at most four edges which joins $a$ to $b$ in $G$ such that $\gamma_{ab}(G)\setminus\{a,b\}\subset G$ and for every $\alpha\in(0,1]$
\bel{L-ABA}
\len_{\alpha,G}(\gamma_{ab}(G))\le \tfrac{12}{\alpha}\,\|a-b\|^\alpha.
\ee
\end{lemma}
\par {\it Proof.} Suppose that $a\in Q_1^{\cl}$ and
$b\in Q_2^{\cl}$. If $a,b\in Q_1^{\cl}$ or $a,b\in Q_2^{\cl}$, then the lemma directly follows from Lemma \reff{Q-DA}. Thus we can assume that $a\in Q_1^{\cl}\sm Q_2^{\cl}$ and $b\in Q_2^{\cl}\sm Q_1^{\cl}$.
\par Let $a=(a_1,a_2),b=(b_1,b_2)$ and let
$$
\Pi(a,b):=[a_1,b_1]\times[a_2,b_2].
$$
Thus $\Pi(a,b)$ is the smallest closed rectangle with sides parallel to the coordinate axes containing $a$ and $b$.
Then by Helly's intersection theorem for rectangles
$$
Q_1^{\cl}\cap Q_2^{\cl}\cap \Pi(a,b)\ne\emp.
$$
\par Let $\tw\in Q_1^{\cl}\cap Q_2^{\cl}\cap \Pi(a,b)$.
Since $\tw\in\Pi(a,b)$, we have
\bel{AB-TW}
\|a-\tw\|,\|b-\tw\|\le \|a-b\|.
\ee
\par Let $R:=G\cap Q_1^{\cl}\cap Q_2^{\cl}$. Then $R$ is either an open line interval or an open rectangle. In both cases $\tw\in R^{\cl}$ so that there  exists $w\in R$ such that $\|\tw-w\|<\|a-b\|$. By this inequality and \rf{AB-TW},
\bel{AB-W}
\|a-w\|,\|b-w\|\le 2\|a-b\|.
\ee
\par Since $w\in Q_1^{\cl}$, by Lemma \reff{Q-DA}, there exists a path $\gamma_1$ (consisting of at most two edges) which joins $a$ to $w$ such that $\gamma_1\sm\{a,w\}\subset Q_1$ and
$$
\len_{\alpha,Q_1}(\gamma_{1})\le \tfrac{3}{\alpha}\,\|a-w\|^\alpha.
$$
\par In a similar way we construct a path $\gamma_2$ (consisting of at most two edges) which connects $b$ to $w$ such that $\gamma_2\sm\{b,w\}\subset Q_2$ and
$$
\len_{\alpha,Q_2}(\gamma_{2})\le \tfrac{3}{\alpha}\,\|b-w\|^\alpha.
$$
\par Since $Q_i\subset G$, we have $\dist(z,\partial Q_i)\le\dist(z,\partial G)$ for every $z\in Q_i$, so that, by definition \reff{SH-LN}, $\len_{\alpha,G}(\gamma_{i})\le \len_{\alpha,Q_i}(\gamma_{i})$, $i=1,2$. Hence,
$$
\len_{\alpha,G}(\gamma_{1})\le \tfrac{3}{\alpha}\,\|a-w\|^\alpha~~~\text{and}~~~
\len_{\alpha,G}(\gamma_{2})\le \tfrac{3}{\alpha}\,\|b-w\|^\alpha.
$$
\par Let $\gamma_{ab}(G):=\gamma_1\cup\gamma_2$. Then $\gamma_{ab}(G)$ is a path consisting of at most four edges which connects $a$ to $b$ in $G$ such that
$$
\len_{\alpha,G}(\gamma_{ab}(G))=\len_{\alpha,G}(\gamma_{1})+ \len_{\alpha,G}(\gamma_{2})\le \tfrac{3}{\alpha}\,\|a-w\|^\alpha+
\tfrac{3}{\alpha}\,\|b-w\|^\alpha.
$$
This inequality and \rf{AB-W} imply \rf{L-ABA} proving the lemma.\bx
\begin{lemma}\lbl{EXT-Q} (i). Let $G\subset\RT$ be one of the following sets:\medskip
\par (a). $G=(Q_1^{\cl}\cup Q_2^{\cl})^\circ$ where $Q_1$ and $Q_2$ are disjoint squares such that $\#(Q_1^{\cl}\cap Q_2^{\cl})>1;$
\par (b). $G=Q_1\cup Q_2\cup Q_3$ where $Q_1$ and $Q_2$ are disjoint squares such that $Q_1^{\cl}\cap Q_2^{\cl}$ is a singleton, and $Q_3$ is a square centered at $Q_1^{\cl}\cap Q_2^{\cl}$.
\par Then $G$ is an $\alpha$-subhyperbolic domain for every $\alpha\in(0,1]$. See Definition \reff{ASHD}. Furthermore,
for every $a,b\in G^{\cl}$, $a\ne b$, there exists a path $\gamma_{ab}(G)$ which joins $a$ to $b$ in $G$ such that $\gamma_{ab}(G)\setminus\{a,b\}\subset G$ and for every $\alpha\in(0,1]$
\bel{L-AB-G1}
\len_{\alpha,G}(\gamma_{ab}(G))\le \tfrac{12}{\alpha}\,\|a-b\|^\alpha\,;
\ee
\par (ii). Every domain $G$ satisfying either condition (a) or condition (b) is a Sobolev $\LM$-extension domain with $e(\LM(G))\le C(m,p)$. See \rf{E-MP}.
\end{lemma}
\par {\it Proof.} If $G$ satisfies condition of part (a), then the statement (i) of the lemma directly follows from  Definition \reff{ASHD} and Lemma \reff{2-SQ}.
\par Let $G$ be a domain from part (b) of the lemma, and let $a,b\in G^{\cl}$. If $a,b\in Q_1^{\cl}\cup Q_3^{\cl}$ or $a,b\in Q_2^{\cl}\cup Q_3^{\cl}$, then, by Lemma \reff{2-SQ}, there exists a path $\gamma_{ab}(G)$ satisfying inequality \rf{L-AB-G1}.
\par Now suppose that $a\in Q_1^{\cl}\sm Q_3$ and $b\in Q_2^{\cl}\sm Q_3^{\cl}$. Let $c$ be the center of the square $Q_3$, i.e., $\{c\}=Q_1^{\cl}\cap Q_2^{\cl}$. Then, by Lemma \reff{Q-DA}, there exists a path $\gamma_{ac}$ joining $a$ to $c$ such that $\gamma_{ac}\sm\{a,c\}\subset Q_1$ and $\len_{\alpha,Q_1}(\gamma_{ac})\le \tfrac{3}{\alpha}\,\|a-c\|^\alpha$. In the same way we prove the existence of a path $\gamma_{cb}$ joining $c$ to $b$ such that $\gamma_{cb}\sm\{b,c\}\subset Q_2$ and $\len_{\alpha,Q_2}(\gamma_{cb})\le \tfrac{3}{\alpha}\,\|b-c\|^\alpha$.
\par Let $\gamma_{ab}(G):=\gamma_{ac}\cup\gamma_{cb}$. Since $Q_1,Q_2\subset G$,
$$
\len_{\alpha,G}(\gamma_{ab}(G))=\len_{\alpha,G}(\gamma_{ac})+ \len_{\alpha,G}(\gamma_{cb})\le \len_{\alpha,Q_1}(\gamma_{ac})+ \len_{\alpha,Q_2}(\gamma_{cb})
$$
so that
$$
\len_{\alpha,G}(\gamma_{ab}(G))\le  \tfrac{3}{\alpha}\,(\|a-c\|^\alpha+\|b-c\|^\alpha).
$$
\par Clearly, $c\in\Pi(a,b)$ where $\Pi(a,b):=[a_1,b_1]\times[a_2,b_2]$ provided
$a=(a_1,a_2)$ and $b=(b_1,b_2)$. Hence
$$
\|a-c\|,\|b-c\|\le \|a-b\|~~~\text{so that}~~~
\len_{\alpha,G}(\gamma_{ab}(G))\le  \tfrac{6}{\alpha}\,\|a-b\|^\alpha
$$
proving inequality \rf{L-AB-G1} and part (i) of the lemma. \medskip
\par It remains to note that part (ii) of the lemma directly follows from part (i) of the present lemma and Theorem \reff{SOB-EXT-RN}.
\par The proof of the lemma is complete.\bx
\bigskip
\par {\bf 4.3. Main geometrical properties of ``The Wide Path''.}
\addtocontents{toc}{~~~~4.3. Main geometrical properties of ``The Wide Path''.\hfill \thepage\par}
Let us give a precise definition of the family of sets $\{\hS_i:1\le i\le k\}$ which we have used in definition \rf{WP-DEF} of ``The Wide Path''. See Section 1.
\begin{definition}\lbl{DF-SH} {\em We put
\bel{SH-1}
\hS_i:=\emp~~\text{if}~~i=k~~\text{or}~~\#(S_i^{\cl}\cap S_{i+1}^{\cl})>1\,.
\ee
\par We also put
\bel{SH-ND}
\hS_i:=S(w_i,\hdl)~~~\text{if}~~~
\#(S_i^{\cl}\cap S_{i+1}^{\cl})=1
\ee
where
\bel{WI-DF}
\{w_i\}:=S_i^{\cl}\cap S_{i+1}^{\cl}
\ee
and 
\bel{DLT-FIN}
\hdl:=\tfrac18\min\{\hdl_1,\hdl_2,\hdl_3\}.
\ee
Here\, $\hdl_1:=\min\{\dist(w_m,\DO):m\in I\}$,\, $\hdl_2:=\min\{\diam S_m:1\le m\le k\}$, and
\bel{DLT-3}
\hdl_3:=\min\{\dist(w_m,S_j):m\in I,1\le j\le k,
j\ne m,m+1\}
\ee
where
$$
I:=\{m\in\{1,...,k-1\}:
\#(S_m^{\cl}\cap S_{m+1}^{\cl})=1\}.
$$}
\end{definition}
\medskip
\par Prove that $\hdl>0$, i.e., that the squares $\hS_i$ in \rf{SH-ND} are well defined. In fact, since $w_i\in[c_i,c_{i+1}]$, $i\in I$, by inclusion \rf{CI-B}, $w_i\in\Omega$. (Recall that $c_i$ denotes the center of the square $S_i$.) Hence, $\hdl_1>0$. It is also clear that $\hdl_2>0$. By part(ii) of Lemma \reff{AD-1-2}, $w_i\notin S_j^{\cl}$ whenever $j\ne i,i+1$, so that $\hdl_3>0$ as well. Hence, $\hdl>0$.\smallskip
\par Our proof of the Sobolev extension property of ``The Wide Path''
$$
\Wc:=\WPT^{(\x,\y)}
$$
relies on a series of results which describe a geometrical structure of $\Wc$ and its complement
$\Hc:=\Omega\sm\Wc.$ Let us recall that
\bel{DF-W2}
\Wc=\left(\,\bigcup_{i=1}^k \left(\,S_i^{\cl}\,\usm\, \hS_i\,\right)\right)^{\circ}\,.
\ee
\par In the next lemma we present several useful properties of the sets  $\hS_i$ which directly follow from Definition \reff{DF-SH}.
\begin{lemma}\lbl{H-IJ} Let $\x,\y\in\Omega$ and let $S_i,S_{i+1}\in \Sc_\Omega(\x,\y)$ be two squares such that
$S_i^{\cl}\cap S_{i+1}^{\cl}$ is a singleton. Then
$$
\diam \hS_i\le \tfrac14\min\{\diam S_i,\diam S_{i+1}\}.
$$
\par Furthermore, the sets of the family $\{2\hS_i:1=1,...,k\}$ are pairwise disjoint subsets of $\Omega$ satisfying the following condition:
\bel{SH-JN}
(2\hS_i^{\cl})\cap S_j^{\cl}=\emp~~\text{for every}~~ 1\le i,j\le k,~j\ne i,i+1.
\ee
\par In particular,
$$
\diam\hS_i\le  2\dist(\hS_i,S_j)~~~\text{for all}~~~
1\le i,j\le k,~j\ne i,i+1,
$$
and
$$
\diam\hS_i+\diam\hS_j\le 4\dist(\hS_i,\hS_j)~~~\text{for all}~~~1\le i,j\le k,~j\ne i.
$$
\end{lemma}
\begin{proposition}\lbl{W-G1} ``The Wide Path'' $\Wc:=\WPT^{(\x,\y)}$ is an open connected subset of $\Omega$ which has the following representation:
\bel{R-WP}
\Wc=\bigcup_{i=1}^{k-1}\left(\,S_i^{\cl}\,\usm \, S_{i+1}^{\cl}\,\usm \,\hS_i\,\right)^{\circ}\,.
\ee
\end{proposition}
\par {\it Proof.} Let
$$
\tW:=\bigcup_{i=1}^{k-1}\left(\,S_i^{\cl}\,\usm\, S_{i+1}^{\cl}\,\usm\, \hS_i\,\right)^{\circ}\,.
$$
Clearly, $\Wc\supset (S_i^{\cl}\cup S_{i+1}^{\cl}\cup \hS_i)^{\circ}$ for every $i=1,...,k-1,$ so that $\Wc\supset\tW$.
\par Prove that $\Wc\subset\tW$. Let $a\in \Wc$. Then, by \rf{DF-W2}, there exists $\delta>0$ such that
\bel{B-SJ}
S(a,\delta)\subset\, \bigcup_{j=1}^k\, \left(\,S_j^{\cl}\,\usm\, \hS_j\,\right)\,.
\ee
\par Let us consider the following two cases.\medskip
\par {\it The first case.} {\it There exists $i\in\{1,...,k-1\}$ such that $a\in 2\hS_i^{\cl}$.}
\par By Lemma \reff{H-IJ}, $(2\hS_i)\cap (2\hS_j)=\emp$ for every $j\ne i, 1\le j\le k$. Furthermore, by \rf{SH-JN},
\bel{J-IN}
(2\hS_i^{\cl})\cap S_j^{\cl}\ne\emp~~~\text{if and only if}~~~j=i~~~\text{or}~~~j=i+1.
\ee
\par Hence, $a\notin 2\hS_j^{\cl}$ provided $j\ne i$ so that
$$
\eta_1:=\tfrac12\,\min\{\dist(a,\hS_j^{\cl}):1\le j\le k,j\ne i\,\}\,>0.
$$
Thus the square $S(a,\eta_1)$ does not cross any square $\hS_j$ whenever $j\ne i$.
\par Also note that, by \rf{J-IN},
\bel{ETA2}
\eta_2:=\tfrac12\,\min\{\dist(a,S_j^{\cl}):1\le j\le k,j\ne i,i+1\,\}\,>0.
\ee
\par Let $\tdl:=\min\{\delta,\eta_1,\eta_2\}$. Then the $\tdl$-neighborhood of $a$, the square $S(a,\tdl)$, may contain only points from the squares $S_i^{\cl}, S_{i+1}^{\cl}$ and $\hS_i$. Hence, by \rf{B-SJ},
$$
a\in\left(S_i^{\cl}\,\usm \, S_{i+1}^{\cl}\,\usm \,\hS_i\,\right)^{\circ}\subset \tW.
$$
\par {\it The second case.} {\it Let $a\in\Wc$ but $a\notin 2\hS_j$ for every $j,1\le j\le k$.} In particular,
$$
a\notin \hS_j^{\cl}~~\text{for every}~~1\le j\le k.
$$
Since
$$
a\in\Wc\subset \bigcup_{i=1}^k\,(S_i^{\cl}\,\usm\, \hS_i)
$$
we conclude that there exists $i\in\{1,...,k\}$ such that $a\in S_i^{\cl}\cap\Omega$.
\par By part (ii) and (iii) of Lemma \reff{H-MS}, we may choose the index $i$ in such a way that either
\bel{FC-A1}
a\in S_i^{\cl}\cup S_{i+1}^{\cl}~~\text{and}~~a\notin S_j^{\cl}~~\text{for every}~~j\ne i,i+1,
\ee
or
\bel{FC2}
\{a\}=\{a_{i+1}\}=S_i^{\cl}\cap S_{i+1}^{\cl}\cap S_{i+2}^{\cl}~\text{and}~a~\text{is a common vertex of}~ S_i, S_{i+1},S_{i+2}\,.
\ee
Furthermore, in this case
\bel{FC-AH2}
a\notin S_j^{\cl}~~\text{for every}~~j\ne i,i+1,i+2,
\ee
and $a$ is a boundary point of the set $S^{\cl}_i\cup S^{\cl}_{i+1}\cup S^{\cl}_{i+2}$. In other words, $a=a_{i+1}$ is a {\it rotation point} of ``The Wide Path'' $\Wc$, and the square $S_{i+1}$ is its {\it rotation square} associated with $a_{i+1}$. See Definition \reff{R-SQ-P}.\medskip
\par We begin with the first case described by \rf{FC-A1}. In this case the quantity $\eta_2$ defined by \rf{ETA2} is positive. Note that the following quantity
$$
\rho_1:=\tfrac12\,\min\{\dist(a,\hS_j^{\cl}):
1\le j\le k\,\}
$$
is positive as well.
\par Let $\rho:=\min\{\delta,\eta_2,\rho_1\}$. Clearly, $\rho>0$. Then the $\rho$-neighborhood of $a$, the square $S(a,\rho)$, does not intersect $S_j^{\cl}$ for all $j\ne i,i+1,1\le j\le k,$ and does not intersect $\hS_j^{\cl}$ for all $1\le j\le k$. Hence, by \rf{B-SJ},
$S(a,\rho)\subset S_i^{\cl}\cup S_{i+1}^{\cl}$ proving
that $a\in  (S_i^{\cl}\cup S_{i+1}^{\cl})^\circ\subset \tW$.\medskip
\par Consider the second case determined by \rf{FC2}. Again in this case $\rho_1>0$. Let
$$
\tau_1:=\tfrac12\,\min\{\dist(a,S_j^{\cl}):1\le j\le k,j\ne i,i+1,i+2\,\}.
$$
Then, by \rf{FC-AH2}, $\tau_1>0$, so that the quantity $\tau:=\min\{\delta,\rho_1,\tau_1\}>0$ as well.
\par Then, by \rf{B-SJ} and by the choice of $\tau$, we have
$$
S(a,\tau)\subset V_i:=S_i^{\cl}\cup S_{i+1}^{\cl}\cup S_{i+2}^{\cl}.
$$
Thus $a$ belongs to the interior of the set $V_i$. On the other hand, $a$ is a boundary point of this set, a contradiction. This contradiction shows that the second case described by \rf{FC2} is impossible proving that $a\in\tW$ for all $a\in \Wc$.
\medskip
\par It remains to show that $\Wc$ is a {\it connected} set. First consider points $c_i$ and $c_j$, $1\le i<j\le k$, the centers of the squares $S_i$ and $S_j$ respectively. Let
$$
\gamma_{ij}:=\bigcup\limits_{m=i}^{j-1}\,[c_m,c_{m+1}].
$$
Prove that $\gamma_{ij}\subset\Wc$. In fact, if $\hS\ne\emp$, i.e., $\#(S_{m}^{\cl}\cap S_{m+1}^{\cl})=1$, then clearly
$$
[c_m,c_{m+1}]\subset(S_{m}^{\cl}\cup S_{m+1}^{\cl}\cup \hS_{m})^\circ=S_{m}\cup S_{m+1}\cup\hS_{m}.
$$
\par Suppose that $\hS_m=\emp$, i.e., that $S_{m}^{\cl}\cap S_{m+1}^{\cl}$ is a line segment
$$
[u_m,v_m]:=S_m^{\cl}\cap S_{m+1}^{\cl}=\partial S_m\cap\partial S_{m+1}.
$$
See \rf{UV-I}.
\par By part (i) of Lemma \reff{H-MS}, $(u_m,v_m)\subset(S_m^{\cl}\cap S_{m+1}^{\cl})^\circ$. On the other hand, by \rf{CI-B},
$[c_m,c_{m+1}]\subset\Omega$. Let
$$
z_m:=[c_m,c_{m+1}]\cap\partial S_m\cap\partial S_{m+1}.
$$
Then $z_m\in\Omega$ so that $z_m\in(u_m,v_m)$. Hence
$$
[c_m,c_{m+1}]\subset S_m\cup S_{m+1}\cup(u_m,v_m)=(S_{m}^{\cl}\cup S_{m+1}^{\cl})^\circ
$$
proving that $[c_m,c_{m+1}]\subset\Wc$. This proves that
$\gamma_{ij}\subset\Wc$ as well.\smallskip
\par Let now $a,b\in\Wc$. Then, by \rf{R-WP}, there exist $i,j\in\{1,...,k\}$ such that
$$
a\in A_i:=(S_{i}^{\cl}\cup S_{i+1}^{\cl}\cup \hS_{i})^\circ~~~\text{and}~~~
b\in A_j:=(S_{j}^{\cl}\cup S_{j+1}^{\cl}\cup \hS_{j})^\circ.
$$
\par By \rf{R-WP}, $A_i,A_j\subset\Wc$. Furthermore, it is clear that $A_i$ and $A_j$ are connected sets containing $c_i$ and $c_j$ respectively. Therefore there exist a path $\gamma_a$ connecting $a$ to $c_i$ in $A_i$, and a path $\gamma_b$ connecting $b$ to $c_j$ in $A_j$. Then the path
$\gamma:=\gamma_a\cup\gamma_{ij}\cup\gamma_j$ joins $a$ to $b$ in $\Wc$.
\par The proposition is completely proved.\bx\medskip
\par Proposition \reff{W-G1} and Lemma \reff{H-MS} enable us to give the following representation of ``The Wide Path'' $\Wc:=\WPT^{(\x,\y)}$. To its formulation we recall that  $[u_i,v_i]=S_i^{\cl}\cap S_{i+1}^{\cl}$ whenever $\#(S_i^{\cl}\cap S_{i+1}^{\cl})>1$, $1\le i<k$. See \rf{UV-I}.
\par Let $1\le i< k$ and let
\bel{T-I}
T_i:=\left\{
\begin{array}{ll}
\hS_i,&~~\text{if}~~
\#(S_i^{\cl}\cap S_{i+1}^{\cl})=1,\\\\
(u_i,v_i),&~~\text{if}~~
\#(S_i^{\cl}\cap S_{i+1}^{\cl})>1.
\end{array}
\right.
\ee
We also put $T_k:=\emp$. We notice a useful formula for the interval $(u_i,v_i)$:
$$
(u_i,v_i)=S_i^{\cl}\cap S_{i+1}^{\cl}\cap\Omega~~~\text{provided}~~~\#(S_i^{\cl}\cap S_{i+1}^{\cl})>1.
$$
\par Now, by \rf{SGM-I} and by definition of $\hS_i$, see \rf{SH-ND},
\bel{B-GF}
\left(\,S_i^{\cl}\,\usm \, S_{i+1}^{\cl}\,\usm \,\hS_i\,\right)^{\circ}=S_i\,\usm \, S_{i+1}\,\usm \,T_i\,.
\ee
Combining this with \rf{R-WP} we obtain the following representation of ``The Wide Path'':
\bel{RW-E}
\Wc=\WPT^{(\x,\y)}=\usm\limits_{i=1}^k \left(S_i\,\usm\,T_i\right)\,.
\ee
C.f. \rf{WP-DEF}. We use this representation in the proof of the following important property of ``The Wide Path''.
\begin{lemma}\lbl{CR-WP} Let $a\in S_i$, $b\in S_{i+1}$, $1\le i<k$, and let $\gamma$ be a path joining $a$ to $b$ in $\WPT^{(\x,\y)}$. Then $\gamma\cap T_{i}\ne\emp$.
\end{lemma}
\par {\it Proof.} Assume that
\bel{A-B1}
\gamma\cap T_{i}=\emp.
\ee
Since $a\in S_i$ and $b\notin S_{i}^{\cl}$, there exists a point $h\in \partial S_i\cap\gamma$ such that the following condition is satisfied: {\it Let $\tgm$ be the subarc of the path $\gamma$ from $h$ to $b$. Then}
\bel{GW}
\tgm\sm\{h\}\subset \RT\sm S_{i}^{\cl}.
\ee
\par Since $h\in\tgm\cap S_{i}^{\cl}$, by \rf{A-B1}, $h\notin S_{i+1}^{\cl}$. On the other hand, $h\in\gamma\subset \WPT^{(\x,\y)}$ so that, by representation \rf{RW-E}, there exists $j, 1\le j\le k,$ $j\ne i$, such that $h\in S_j\cup T_j$. See \rf{T-I}.\medskip
\par Clearly, since $h\in\partial S_i$ and the squares of ``The Wide Path'' are touching, $h\notin S_j$ for every $j$, $1\le j\le k$. Hence $h\in T_j\cup S_{j}^{\cl}$ for some $j, 1\le j\le k$, $j\ne i$.
\par Prove that $j=i-1$. (In particular, it shows that $i\ge 2$.) If $\#(S_{j}^{\cl}\cap S_{j+1}^{\cl})>1$ for some $j\ne i$, $1\le j\le k,$ then
$$
T_j=(u_j,v_j)\subset  S_{j}^{\cl}\cap
S_{j+1}^{\cl}.
$$
See \rf{T-I} and \rf{UV-I}.
\par Hence $S_{j}^{\cl}\cap S_{j+1}^{\cl}\cap S_{i}^{\cl}\ni h$ so that
$S_{j}^{\cl}\cap S_{j+1}^{\cl}\cap S_{i}^{\cl}\ne \emp$.
Then, by part (ii) of Lemma \reff{AD-1-2}, $|j-i|\le 1$ and $|j+1-i|\le 1$. Since $j\ne i$, this implies that $j=i-1$.
\par Now let $\#(S_{j}^{\cl}\cap
S_{j+1}^{\cl})=1$ for some $j\ne i$, $1\le j\le k$, i.e., $T_j=\hS_j$, see \rf{T-I}. Then $\hS_j\cap S_{i}^{\cl}\ne\emp$ so that, by Lemma \reff{H-IJ}, see \rf{SH-JN}, either $j=i$ or $j=i-1$. But
we know that $j\ne i$ so that in this case $j=i-1$ as well.\smallskip
\par Thus
$$
h\in T_{i-1}\cap\partial S_i\cap\tgm
$$
where $\tgm$ is a path joining $h$ to $b$ in $\Omega$ which satisfies \rf{GW}.
\par Consider again two cases. If $\#(S_{i-1}^{\cl}\cap S_{i}^{\cl})>1$, i.e., if $T_{i-1}=(u_{i-1},v_{i-1})$, we have $T_{i-1}\subset S_{i-1}^{\cl}\cap S_{i}^{\cl}$ (see \rf{UV-I}), so that $h\in S_{i-1}^{\cl}\cap\Omega$. But $b\in S_{i+1}$ so that, by part (i) of Lemma \reff{AD-1-2}, $\tgm\cap S_i\ne\emp$ which contradicts \rf{GW}.
\par Consider the remaining case where $\#(S_{i-1}^{\cl}\cap S_{i}^{\cl})=1$, i.e., $T_{i-1}=\hS_{i-1}$. Choose a point
$\thw\in S_{i-1}\cap \hS_{i-1}$. It is clear that $\hS_{i-1}\sm S_{i}^{\cl}$ is a connected set so that we can join $h$ to $\thw$ by a path $\gamma_1$ which lies in $\hS_{i-1}\sm S_{i}^{\cl}$. Then the path $\gamma_2:=\gamma_1\cup \tgm$ connects in $\Omega$ the point $\thw\in S_{i-1}$ to the  point $b\in S_{i+1}$. Furthermore, $\gamma_2\cap S_i=\emp$. But this again contradicts part(i) of Lemma \reff{AD-1-2}.
\par The proof of the lemma is complete.\bx
\begin{proposition}\lbl{ST-CD} Let $\Hc=\Omega\sm\Wc$ and let $H$ be a connected component of $\Hc$. Suppose that there exist $i$ and $j$, $1\le i,j\le k$, such that
$$
H\cap S_i^{\cl}\cap \Omega\ne\emp~~~\text{and}~~~
H\cap S_j^{\cl}\cap \Omega\ne\emp.
$$
\par Then $|i-j|\le 1$.
\end{proposition}
\par {\it Proof.} Without loss of generality we may assume that $i\le j$. Suppose that $i+1<j$.
\par Let
$$
a\in H\cap S_i^{\cl}\cap \Omega~~\text{and let}~~
b\in H\cap S_j^{\cl}\cap \Omega.
$$
\par Since $a,b\in H$ and $H$ is a connected component of $\Hc$, there exists a path $\gamma$ connecting $a$ to $b$ in $\Hc$. We know that $\Hc\cap\Wc=\emp$ so that  $\gamma\cap\Wc=\emp$ as well. In particular, since  $S_{i+1}\subset \Wc$, see \rf{RW-E}, we conclude that  and $\gamma\cap S_{i+1}=\emp$. We also notice that $\gamma\subset \Hc\subset\Omega$.
\par On the other hand, $a\in S_i^{\cl}\cap \Omega$, $b\in S_j^{\cl}\cap \Omega$ and $i<j-1$, so that, by part (i) of Lemma \reff{AD-1-2}, $\gamma\cap S_{i+1}\ne\emp$, a contradiction.
\par This contradiction shows that our assumption that $i+1<j$ is not true, and the proof of the lemma is complete.\bx
\begin{proposition}\lbl{O-CMS} Let $H$ be a connected component of $\Hc=\Omega\sm \Wc$. Then \medskip
\par (i) either  there exists $i\in\{1,2,...,k\}$ such that
\bel{FC-1}
H\cap S_i^{\cl}\ne\emp~~~\text{and}~~~H\cap S_j^{\cl}=\emp~~\text{for every}~~1\le j\le k,\,j\ne i,
\ee
\par (ii) or there exists $i\in\{1,2,...,k-1\}$ such that
\bel{FC-2}
H\cap S_i^{\cl}\ne\emp,\,
H\cap S_{i+1}^{\cl}\ne\emp~~\text{and}
~~H\cap S_j^{\cl}=\emp~~\text{for every}~1\le j\le k,\,j\ne i,i+1\,.
\ee
\par Furthermore, in case (i)
\bel{S-G1}
H\cup S_i~~\text{is a subdomain  of}~~\Omega\,.
\ee
In turn, in case (ii)
\bel{S-G2}
H\cup S_i\cup S_{i+1}\cup T_i~~\text{is a subdomain  of}~~\Omega\,.
\ee
\end{proposition}
\par {\it Proof.} An example of connected components of the set $\Hc=\Omega\sm \Wc$ is given on Figure 18. In this example each of the connected components $H_1,...,H_5$  touches exactly one square from the family of squares $\Sc=\{S_1,...,S_{10}\}$. Thus the components $H_i, i=1,...,5,$ satisfy condition (i) of the lemma. Other connected components of $\Hc$ satisfy condition (ii), i.e., each of these components touches exactly two squares from $\Sc$.
\begin{figure}[h]
\center{\includegraphics[scale=0.9]{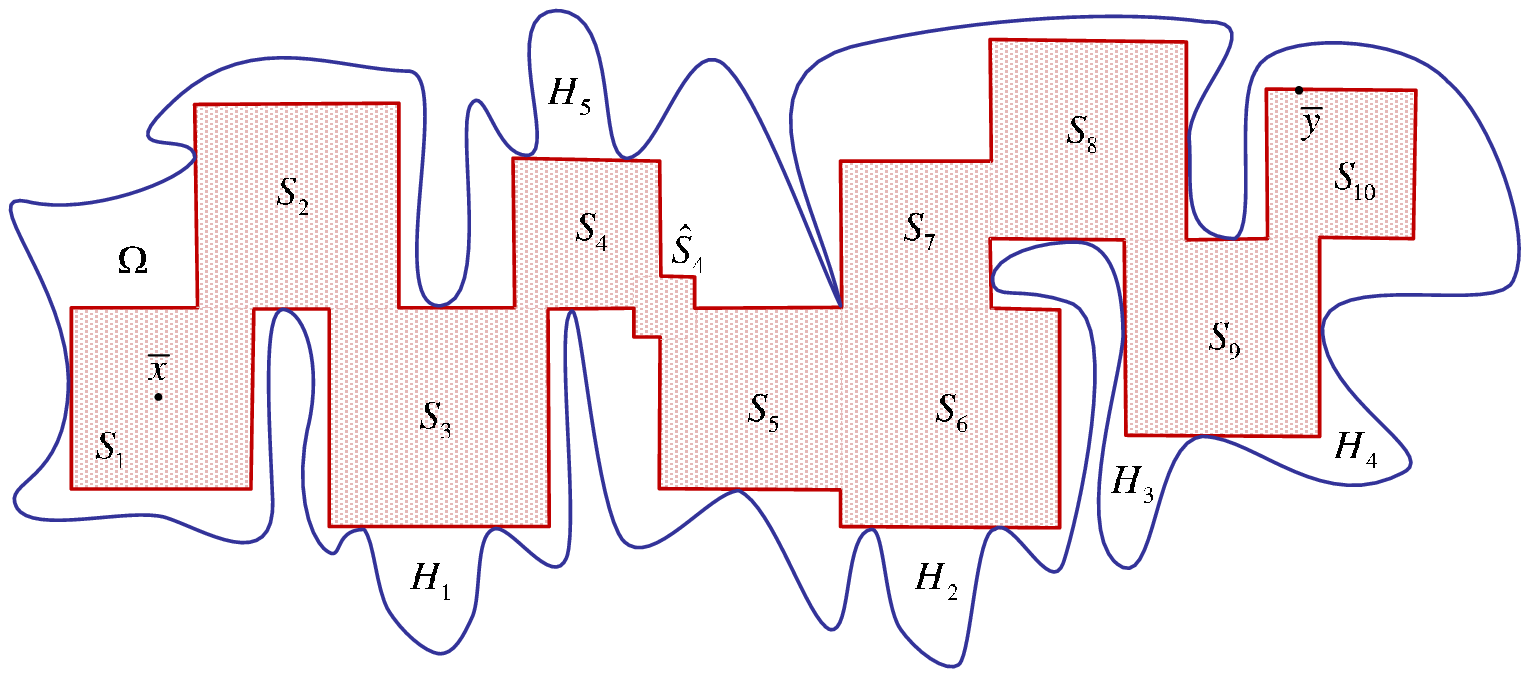}}
\caption{~}
\end{figure}
\medskip
\par We turn to the proof of the lemma. First let us prove that
\bel{DHA}
\DHA\cap H\ne\emp.
\ee
\par Fix a point $z_0\in H$. If $z_0\in\DHA$, then \rf{DHA} is proven. Suppose that $z_0\in H^\circ$. We know that $\x=c_1$, i.e., that $\x$ is the center of $S_1$. Hence, $\x\in\Wc$, see \rf{RW-E}. Let $\gamma$ be a path connecting $\x$ to $z_0$ in $\Omega$ so that $\gamma$ is a graph of a continuous mapping $\Gamma:[0,1]\to\Omega$ such that $\Gamma(0)=z_0$ and $\Gamma(1)=\x$.
\par Let $Y:=\{t\in[0,1]: \Gamma(t)\in \Wc\}$ and let $t':=\inf Y$. Since $\Gamma$ is a continuous mapping,  $z_0\in H^\circ$ and $\x\in \Wc^\circ=\Wc$, we conclude that $0<t'<1$.
\par Let $\tz:=\Gamma(t').$ Then, by definition of $t'$, the subarc of $\gamma$ from $z_0$ to $\tz$ lies in the set $\Hc=\Omega\sm\Wc$. Since $H$ is a {\it connected component of $\Hc$}, $\tz\in H$.
\par On the other hand, since $t'=\inf Y\notin Y$, there exists a sequence $\{t_m: m=1,2,...\}\subset Y$ which converges to $t'$ as $m\to\infty$. Let  $h_m:=\Gamma(t_m)$, $m=1,2,...$. Then $h_m\in\Wc$, $h_m\ne h_n$, if $m\ne n$ (because $\gamma$ is a {\it simple} path), and $h_m\to \tz$ as $m\to\infty$. Hence $\tz\in \DHA\cap H$ proving \rf{DHA}.\medskip
\par Prove the statements (i) and (ii). Since the parameter $k$ in representation \rf{RW-E} is finite, there exists $i\in\{1,...,k\}$ and an infinite subsequence $\{h_{m_j}:j=1,2,...\}$ of the sequence $\{h_{m}:m=1,2,...\}$ such that
$$
h_{m_j}\in S_i\cup T_i~~~\text{for all}~~~ j=1,2,...\,.
$$
Since $h_m\to \tz$ as $m\to\infty$, the subsequence $
h_{m_j}\to \tz$ as $j\to\infty$ proving that
$$
\tz\in S_i^{\cl}\cup T_i^{\cl}\,.
$$
\par Recall that the set $T_i$ is defined by \rf{T-I}. In particular, $T_k=\emp$ and $T_i\subset S_i^{\cl}$ whenever $\#(S_i^{\cl}\cap S_{i+1}^{\cl})>1$. Thus, in this case $\tz\in S_i^{\cl}$.\smallskip
\par Suppose that $1\le i<k$ and $\#(S_i^{\cl}\cap S_{i+1}^{\cl})=1$. In this case $\hS_i\ne\emp$ and is defined by the formula \rf{SH-ND}. Let us assume that $\tz\in\hS_i^{\cl}$ and prove that in this case
$$
H\cap S_i^{\cl}\ne\emp~~~\text{and}~~~
H\cap S_{i+1}^{\cl}\ne\emp.
$$
\par In fact, since $\tz\in\Hc=\Omega\sm\Wc$ and $S_i,S_{i+1},\hS_i\subset\Wc$, we have
$\tz\in\partial \hS_i\sm (S_i\cup S_{i+1})$.
\par We also recall that, by Lemma \reff{H-IJ}, see \rf{SH-JN}, $2\hS_i^{\cl}\cap S_j^{\cl}=\emp$ for every $j\in\{1,...,k\}$ such that $j\ne i,i+1$. This lemma also states that $(2\hS_i)\cap (2\hS_j)=\emp$ for every
$j\in\{1,...,k\}$, $j\ne i$. Hence, by representation \rf{R-WP} (or \rf{RW-E}), we have
$$
U_i:=(2\hS_i)\sm (S_i\cup S_{i+1}\cup\hS_i)\subset \Hc=\Omega\sm\Wc.
$$
\par Clearly, there exist a point $z_i\in S_i^{\cl}\cap U_i$ and a path $\gamma_1$ in $U_i$ which joins $\tz$ to
$z_i$. Hence, $z_i\in H\cap S_i^{\cl}$. Also there exist
a point $z_{i+1}\in S_{i+1}^{\cl}\cap U_i$ and a path $\gamma_2$ connecting $\tz$ to $z_{i+1}$ in $U_i$, so that  $z_{i+1}\in H\cap S_{i+1}^{\cl}$.
See Figure 19.
\begin{figure}[h]
\center{\includegraphics[scale=0.8]{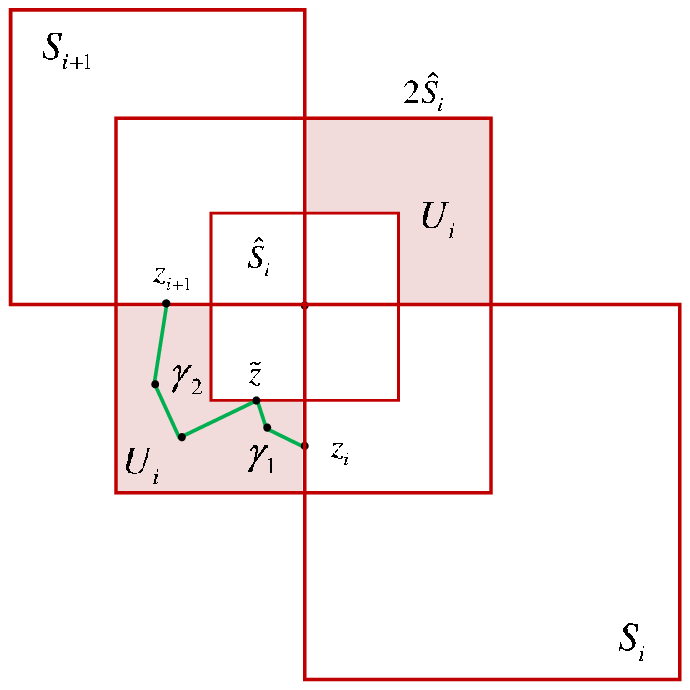}}
\caption{~}
\end{figure}
\medskip
\par Thus we have proved that either there exists
$i\in\{1,...,k\}$ such that $H\cap S_i^{\cl}\ne\emp$, or
there exists $i\in\{1,...,k-1\}$ such that $H\cap S_i^{\cl}\ne\emp$ and $H\cap S_{i+1}^{\cl}\ne\emp$. Then, by Proposition \reff{ST-CD}, all the conditions of part (i) and part (ii) are satisfied. See \rf{FC-1} and \rf{FC-2}.\smallskip
\par Prove \rf{S-G1}. Let $a\in H\cup S_i$. We have to find $\delta>0$ such that $S(a,\delta)\subset H\cup S_i$ provided conditions \rf{FC-1} hold.
\par Since $a\notin S_j^{\cl}$ for every $j\ne i$,
$$
\delta_1:=\tfrac12\dist
(a,\usm\limits_{j\ne i}S_j^{\cl})>0\,.
$$
\par As we have proved above,
$$
H\cap \hS_i^{\cl}\ne\emp~~\Longrightarrow~~H\cap S_i^{\cl}\ne\emp~~~\text{and}~~~
H\cap S_{i+1}^{\cl}\ne\emp.
$$
But, by \rf{FC-1}, $H\cap S_{i+1}^{\cl}=\emp$ so that
$$
H\cap \hS_j^{\cl}=\emp~~~\text{for every}~~~
j\in\{1,...,k\}.
$$
Hence
$$
\delta_2:=\tfrac12\dist
(a,\usm\limits_{j=1}^k \hS_j^{\cl})>0\,.
$$
\par Let $\delta_3:=\tfrac12\dist(a,\DO)$ and let
$$
\delta:=\min\{\delta_1,\delta_2,\delta_3\}\,.
$$
\par Then, by \rf{R-WP},
$$
S(a,\delta)\cap \Wc=S(a,\delta)\cap S_i \,.
$$
Hence,
$$
S(a,\delta)\cap \Hc=S(a,\delta)\cap (\Omega\sm \Wc)=
S(a,\delta)\sm S_i\,.
$$
Clearly, $S(a,\delta)\sm S_i$ is a {\it connected set} so that each $z\in S(a,\delta)\sm S_i$ can be joined to $a$ by a path $\gamma_z\subset S(a,\delta)\sm S_i\subset\Hc$. This implies that $z$ and $a$ {\it belong to the same connected component} of $\Hc$, i.e., that $z\in H$.
\par Hence $S(a,\delta)\sm S_i\subset H$ proving that $S(a,\delta)\subset S_i\cup H$.\medskip
\par Prove that $S_i\cup H$ is a connected set. We know that $H\cap S_i^{\cl}\ne\emp$ so that there exists $a\in H\cap S_i^{\cl}$.
\par Let $z\in H$. Since $H$ is a connected component of $\Hc$, this set is connected so that there exists a path $\gamma_z$ joining $z$ to $a$ in $H$. Then a path $\gamma=\gamma_z\cup [a,c_i]$ connects $z$ to $c_i$ in $S_i\cup H$. Thus each point $z\in S_i\cup H$ can be connected to $c_i$, the center of $S_i$, by a path in $S_i\cup H$ proving that this set is connected.\medskip
\par We turn to the proof of the statement \rf{S-G2}, the last statement of the proposition. Let $H$ be a connected component of $\Hc=\Omega\sm\Wc$ satisfying conditions \rf{FC-2}. Let
\bel{VI}
V_i:=S_i\cup S_{i+1}\cup T_i,
\ee
see \rf{T-I}, and let
\bel{GI}
a\in G_i:=H\cup V_i\,.
\ee
\par Prove the existence of $\ve>0$ such that $S(a,\ve)\subset G_i$. By \rf{FC-2},
$$
\ve_1:=\tfrac12\dist
(a,\usm\limits_{j\ne i,i+1}S_j^{\cl})>0\,.
$$
\par In the same way as we have proved \rf{FC-2}, we show that
$$
H\cap \hS_j^{\cl}=\emp~~~\text{for every}~~~
1\le j\le k,~j\ne i\,.
$$
Hence
$$
\ve_2:=\tfrac12\dist
(a,\usm\limits_{j\ne i}\hS_j^{\cl})>0\,.
$$
\par Finally, we put $\ve_3:=\tfrac12\dist(a,\DO)$,
$$
\ve_4:=\tfrac18\diam T_i,
$$
and
$$
\ve:=\min\{\ve_1,\ve_2,\ve_3,\ve_4\}\,.
$$
\par Then, by \rf{R-WP},
$$
S(a,\ve)\cap\Wc=S(a,\ve)
\cap(S_i^{\cl}\cup S_{i+1}^{\cl}\cup\hS_i)^\circ
$$
so that, by \rf{B-GF},
$$
S(a,\ve)\cap\Wc
=S(a,\ve)\cap(S_i\cup S_{i+1}\cup T_i)=S(a,\ve)\cap V_i\,.
$$
See \rf{VI}. Hence,
$$
S(a,\ve)\cap\Hc=S(a,\ve)\cap(\Omega\sm\Wc)
=S(a,\ve)\setminus V_i\,.
$$
\par It can be readily seen that, by definition of $\ve_4$, the set $S(a,\ve)\sm V_i$ is a connected set. Therefore every $z\in S(a,\ve)\sm V_i$ can be joined to $a$ by a path $\gamma_z\subset S(a,\ve)\sm V_i\subset \Hc$. Hence it follows that $z$ and $a$ belong to the same connected component of $\Hc$, i.e., that $z\in H$.
\par Thus we have proved that $S(a,\ve)\sm V_i\subset H$ so that $S(a,\ve)\subset V_i\cup H=G_i$. See \rf{GI}.
\par It remains to prove that the set $H\cup V_i$ is connected. The proof of this property is similar to that for the case \rf{FC-1}. As in that case we know that $H\cap S_i^{\cl}\ne\emp$ so that, using the same approach, we show that for every $z\in H$ there exists a path $\gamma\subset H\cup V_i$ joining $z$ to $c_i$. Clearly, $V_i$ is a connected set and $c_i\in V_i$. Hence $c_i$ can be connected by a path in $H\cup V_i$ to an {\it arbitrary} point $z\in H\cup V_i$ proving the connectedness of this set.
\par The proof of the proposition is complete.\bx
\bigskip
\par {\bf 4.4. Extensions of Sobolev functions defined on  ``The Wide Path''.}
\addtocontents{toc}{~~~~4.4. Extensions of Sobolev functions defined on ``The Wide Path''.\hfill \thepage\par\VST}
Proposition \reff{O-CMS} motivates us to introduce several important geometrical objects related to ``The Wide Path'' $\WPT^{(\x,\y)}$. Let
$$
\Cc:=\{H: H~~\text{is a connected component of}~~\Hc=\Omega\sm \Wc\}.
$$
Given $i\in\{1,...,k\}$ we define a subfamily $\Fc_i$ of $\Cc$ by
$$
\Fc_i:=\{H\in\Cc:H\cap S_i^{\cl}\ne\emp\,
~\text{and}
~H\cap S_j^{\cl}=\emp~\forall~1\le j\le k,\,j\ne i\}
\,.
$$
C.f., part(i) of Proposition \reff{O-CMS}. In turn, part (ii) of this proposition motivates us to introduce a subfamily $\Pc_i$ of $\Cc$ as follows: given $i\in\{1,...,k-1\}$ we put
$$
\Pc_i:=\{H\in\Cc:H\cap S_i^{\cl}\ne\emp,\,
H\cap S_{i+1}^{\cl}\ne\emp~\text{and}
~H\cap S_j^{\cl}=\emp~\forall~1\le j\le k,\,j\ne i,i+1\}\,.
$$
\par Note that, by Proposition \reff{O-CMS}, the family
\bel{PF-D}
\Fc\Pc:=\{\Fc_1,...,\Fc_k,\Pc_1,...,\Pc_{k-1}\}
\ee
provides a {\it a partition} of the family $\Cc$ of all connected components of the set $\Hc=\Omega\sm\Wc$. In other words, $\Fc\Pc$ consists of {\it pairwise disjoint sets} which cover the family $\Cc$, i.e.,
\bel{C-PT}
\Cc=\left(\usm\limits_{i=1}^k\Fc_i\right)\usm
\left(\usm\limits_{j=1}^{k-1}\Pc_i\right)\,.
\ee
\par The collection $\Fc\Pc$ enables us to introduces the following families of subsets of $\Omega$:
\bel{D-PHI}
\Phi_i:= \left(\,\usm\limits_{H\in\,\Fc_i} H\right)
\usm\, S_i, ~~~1\le i\le k\,.
\ee
and
\bel{D-PSI}
\Psi_i:= \left(\,\usm\limits_{H\in\,\Pc_i} H\right)
\usm\, S_i\,\usm\, S_{i+1}\,\usm\, T_i,~~~1\le i\le k-1\,.
\ee
\par Finally we put
$$
\Lambda:=\{\Phi_1,...,\Phi_k,\Psi_1,...,\Psi_{k-1}\}\,.
$$
\par The following proposition describes the main properties of the collection $\Lambda$. To its formulation given a family $\Ac=\{A_\alpha:\alpha\in I\}$ of sets in $\RT$ we let $M(\Ac)$ denote its {\it covering multiplicity}, i.e., the minimal positive integer $M$ such that every point $z\in\RT$ is covered by at most $M$ sets $A_\alpha$ from the family $\Ac$.
\begin{proposition}\lbl{LM-P} (i) The family $\Lambda$ consists of subdomains of $\Omega$ which cover $\Omega$ with co\-ve\-ring multiplicity $M(\Lambda)\le 3$;\smallskip
\par (ii) Let
$$
\Lambda_{\Hc}:=\{\Phi_1\sm\Wc,...,\Phi_k\sm\Wc,\Psi_1\sm\Wc,...,
\Psi_{k-1}\sm\Wc\}.
$$
Then the family $\Lambda_{\Hc}$ consists of pairwise disjoint sets;\smallskip
\par (iii)  For every domain $G\in\Lambda$ the set $G\cap\Wc$ is a Sobolev $\LM$-extension domain satisfying the following inequality
$$
e(\LM(G\cap\Wc))\le\, C(m,p).
$$
See \rf{E-MP}.
\end{proposition}
\par {\it Proof.} Prove (i). By Proposition \reff{O-CMS}, see \rf{S-G1}, for each connected component $H\in\Fc_i$, $1\le i\le k$, the set $H\cup S_i$ is open and connected. In turn, by  \rf{S-G2}, the set $H\cup V_i$ where
\bel{VI-D}
V_i:=S_i\cup S_{i+1}\cup T_i,~~~i=1,...k-1,
\ee
is open and connected provided $H\in\Pc_i$. Combining these facts with formulae \rf{D-PHI} and \rf{D-PSI}, we obtain that every set $G\in\Lambda$ is a union of domains which have a non-empty intersection. Hence $G$ is a {\it domain} as well.\smallskip
\par Recall that the family $\Fc\Pc$ defined by \rf{PF-D} is a partition of $\Cc$, see \rf{C-PT}. Combining this property with representation \rf{RW-E} of ``The Wide Path'' $\Wc$ we conclude that
$$
\Omega=\bigcup_{G\in \Lambda}\,G
$$
proving that $\Lambda$ is a {\it covering} of $\Omega$.\medskip
\par In a similar way we prove part (ii) of the proposition. In fact, by \rf{D-PHI} and \rf{D-PSI},
$$
\Phi_i\cap\Hc=\Phi_i\sm\Wc=\usm\limits_{H\in\,\Fc_i} H,
~~~1\le i\le k,
$$
and
$$
\Psi_i\cap\Hc=\Psi_i\sm\Wc=\usm\limits_{H\in\,\Pc_i} H, ~~~1\le i\le k-1.
$$
\par But the collection $\Fc\Pc$ is a partition of the family $\Cc$, see \rf{C-PT}, so that distinct members of the family $\Lambda_{\Hc}$ have no common points.\medskip
\par Prove that $M(\Lambda)\le 3$. Let $z\in\Hc=\Omega\sm\Wc$ and let $H\in\Cc$ be a connected component of $\Hc$ containing $z$. Since $\Fc\Pc$, see \rf{PF-D}, is a partition of the family $\Cc$ of {\it all} connected component of $\Hc$, there exists a {\it unique} domain $G\in\Lambda$ which contains $z$.
\par This also proves that $M(\Lambda)=\max\{1,M(\Lambda_{\Wc})\}$ where
$$
\Lambda_{\Wc}:=\{\Phi_1\cap\Wc,...,
\Phi_k\cap\Wc,\Psi_1\cap\Wc,...,
\Psi_{k-1}\cap\Wc\}.
$$
\par Note that, by definitions \rf{D-PHI} and \rf{D-PSI},
$$
\Lambda_{\Wc}=\{S_1,...,S_k,V_1,...,V_{k-1}\}
$$
where $V_i$ is defined by \rf{VI-D}.
\par It can be readily seen that $M(\Lambda_{\Wc})\le 3$. In fact, suppose that $z\in S_i$ for some $i\in\{1,...,k\}$. Then the point $z$ can also belong to $V_{i-1}=S_{i-1}\cup S_i\cup T_{i-1}$ and $V_{i}=S_{i}\cup S_{i+1}\cup T_{i}$. Other members of the family $\Lambda_{\Wc}$ do not contain $z$. (This follows from properties of the squares $\{S_j\}$ presented in Lemmas \reff{MAIN-WP}, \reff{AD-1} and \reff{AD-1-2}.) Thus in this case $z$ can be covered by at most $3$ members of the family $\Lambda_{\Wc}$.
\par Let $z\in T_i$ for certain $i\in\{1,...,k-1\}$, see \rf{T-I}. Clearly, in this case $z\in V_i$. By \rf{T-I}, if $\#(S_i^{\cl}\cap S_{i+1}^{\cl})>1$, i.e., if $T_i=(u_i,v_i)$, there are no exist other members of $\Lambda_{\Wc}$ which contain $z$. Whenever $\#(S_i^{\cl}\cap S_{i+1}^{\cl})=1$, i.e., $T_i=\hS_i$, only the squares $S_i$ and $S_{i+1}$ from the family $\Lambda_{\Wc}$ can contain $z$. (As in the previous case it directly follows from Lemmas \reff{MAIN-WP}, \reff{AD-1} and \reff{AD-1-2}.) Thus in this case again the point $z$ is covered by at most $3$ members of $\Lambda_{\Wc}$ proving that $M(\Lambda_{\Wc})\le 3$.
\par Hence $M(\Lambda)=\max\{1,M(\Lambda_{\Wc})\}\le 3$.
\medskip
\par Prove part (iii) of the proposition. Let $G\in\Lambda$. Then either $G=\Phi_i$ for some $i\in\{1,...,k\}$, or $G=\Psi_i$ for certain index $i\in\{1,...,k-1\}$. Hence either $G\cap\Wc=\Phi_i\cap\Wc=S_i$ or
$G\cap\Wc=\Psi_i\cap\Wc=V_i$. See \rf{VI-D}.
\par Then, by \rf{B-GF}, either
$V_i=(S_i^{\cl}\cup S_{i+1}^{\cl})^{\circ}$ or
$V_i=S_i\cup S_{i+1}\cup\hS_i$. Combining this description of $V_i$ with the statement of Lemma \reff{EXT-Q} we conclude that the set $G\cap\Wc$ is a Sobolev extension domain such that $e(\LM(G\cap\Wc))\le C(m,p)$.
\par The proposition is completely proved.\bx\bigskip
\par We turn to the proof of Theorem \reff{WP-EXT}. Clearly, this theorem immediately follows from definition \rf{E-MP} and the following result.
\begin{theorem} \lbl{WP-OM-EXT} Let $p>2$ and $m\in\N$. Let $\x,\y\in\Omega$ where $\Omega$ a simply connected bounded domain in $\RT$. Suppose that $\Omega$ is a Sobolev $L^m_p$-extension domain.
\par Let $\Wc=\WPT^{(\x,\y)}$ be a``Wide Path'' joining $\x$ to $\y$ in $\Omega$ and let $f\in \LM(\Wc)$. Then $f$ can be extended to a function $F\in\LMPO$ such that
$$
\|F\|_{\LMPO}\le C(m,p)\,\|f\|_{\LM(\Wc)}\,.
$$
\end{theorem}
\par For the proof of Theorem \reff{WP-OM-EXT} we are needed the following two auxiliary results.
\begin{proposition}\lbl{LOC-W}(\cite{MIC}, p. 128) If $\GW$ is a collection of non-empty open sets in $\RN$ whose union is $U$ and if $F\in L_{1,loc}(U)$ is such that for some multi-index $\alpha$ the $\alpha$-th weak derivative of $F$ exists on each member of $\GW$, then $F$ has the $\alpha$-th weak derivative on $U$.
\end{proposition}
\begin{proposition}\lbl{EXT-DM} Let $m\in\N$ and $1\le p<\infty$ and let $V$ be a domain in $\RT$.
\par Let $\Gc=\{G_i: i\in I\}$ be a family of domains in $\RT$ satisfying the following conditions:\smallskip
\par (i) $\Gc$ has finite covering multiplicity $M=M(\Gc)$;\smallskip
\par (ii) The sets of the family $\{G_i\sm V: i\in I\}$ are pairwise disjoint;\smallskip
\par (iii) For every $G\in \Gc$ the set $G\cap V$ is a non-empty Sobolev $\LM$-extension domain. Furthermore,
\bel{E-FN}
A:=\sup_{G\in\Gc}\,e(\LM(G\cap V))<\infty.
\ee
\par Let
\bel{U-V}
U:=V\,\usm\,\left\{\usm\limits_{G\in\Gc}G\right\}\,.
\ee
\par Then every function $f\in \LM(V)$ can be extended to a function $F\in \LM(U)$. Furthermore, $F$ depends on $f$ linearly and
$$
\|F\|_{\LM(U)}\le C\,M^{\frac1p}\,A\,\|f\|_{\LM(V)}
$$
where $C=C(m,p)$.
\end{proposition}
\par {\it Proof.} Let  $f\in \LM(V)$. We define the required extension $F$ of $f$ as follows. Let $G\in\Gc$. Then, by (iii), the set $G\cap\Wc$ is a Sobolev extension domain such that $e(\LM(G\cap V))\le A$, see \rf{E-FN}. Therefore there exists a function $F_G\in\LMT$
such that
$$
F_G|_{G\cap V}=f|_{G\cap V}
$$
and
\bel{FG-R}
\|F_G\|_{\LMT}\le A\,\|f|_{G\cap V}\|_{\LM(G\cap V)}\,.
\ee
\par By \rf{U-V} and by condition (ii), for each $z\in U\sm V$ there exists a unique domain $G^{(z)}\in\Gc$ such that $G^{(z)}\sm V\ni z$.
\par This property enables us to define the extension $F$ of $f$ by the following formula:
$$
F(z):=\left\{
\begin{array}{ll}
f(z),& z\in V,\\\\
F_{G^{(z)}}(z),& z\in U\sm V\,.
\end{array}
\right.
$$
\par Thus
\bel{EXT-FM}
F|_{G}=F_G|_{G}~~~\text{for every}~~G\in\Gc\,.
\ee
\par Prove that $F\in\LM(U)$. We know that the restriction of $F$ to $V$ and to any subdomain $G\in\Gc$ is a Sobolev function on $G$ so that {\it each weak derivative of $F$ of order at most $m$ exists on $G$}. Hence, by Proposition \reff{LOC-W}, {\it all partial distributional derivatives of $F$ of all orders up to $m$ exist on all of $U$}.
\smallskip
\par Now let us estimate the norm of $F$ in $\LM(U)$. We add the set $V$ to the family $\Gc$ and denote the new family by $\GW$. Clearly, by \rf{U-V}, the sets of the family $\GW$ cover the set $U$ so that
\be
\|F\|_{\LM(U)}^p&\le& C\,\sbig_{|\alpha|\le m}\, \intl_{U}|D^\alpha F|^p\,dz\le
C\,\sbig_{|\alpha|\le m}\,\,\sbig_{G\in\GW}\,
\intl_{G}|D^\alpha F|^p\,dz\nn\\
&=&
C\,\sbig_{|\alpha|\le m}\,\,\sbig_{G\in\GW}\,\intl_{G}
|D^\alpha F_G|^p\,dz
=\,C\,\sbig_{G\in\GW}\,\,\sbig_{|\alpha|\le m}\,\, \intl_{G}|D^\alpha F_G|^p\,dz\,.
\nn
\ee
Here $C=C(m,p)$. Hence, by \rf{FG-R},
$$
\|F\|_{\LM(U)}^p \le  \,C\,A^p\,\sbig_{G\in\GW}\,\,\sbig_{|\alpha|\le m}\,\,
\intl_{G\cap V}
|D^\alpha f|^p\,dz
=\,C\,A^p\,\sbig_{|\alpha|\le m}\,\,\sbig_{G\in\GW}\,\,
\intl_{G\cap V} |D^\alpha f|^p\,dz\,.
$$
\par By condition (i), covering multiplicity of the family $\{G\cap V: G\in \GW\}$ is bounded by $M+1$. Hence
$$
\|F\|_{\LM(U)}^p \le
\,C\,A^p\,(M+1)\,\sbig_{|\alpha|\le m}\,\,
\intl_{V} |D^\alpha f|^p\,dz\le \,C\,A^p\,M\,\|f\|_{\LM(V)}^p\,.
$$
\par It remains to note that, since $F_G$ depends on $f$ linearly, by \rf{EXT-FM}, {\it the function $F$ depends on $f$ linearly}. The proof of the proposition is complete.\bx
\bigskip
\par {\it Proof of Theorem \reff{WP-OM-EXT}.}\smallskip
\par Let $\x,\y\in\Omega$ and let $\Wc=\WPT^{(\x,\y)}$ be ``The Wide Path'' joining $\x$ to $\y$ in $\Omega$. We suppose that $\Omega$ is a Sobolev extension domain satisfying  condition \rf{S-EXT} for some $\CE\ge 1$.
Therefore, by Proposition \reff{LM-P}, there exists a finite family
$$
\Lambda:=\{\Phi_1,...,\Phi_k,\Psi_1,...,\Psi_{k-1}\}
$$
of subdomains of $\Omega$ satisfying conditions (i)-(iii) of this proposition. These conditions imply conditions (i)-(iii) of Proposition \reff{EXT-DM} provided 
\bel{DF-UVG}
U:=\Omega,~~V:=\WPT^{(\x,\y)}~~~\text{and}~~~\Gc:=\Lambda.
\ee
In these settings, by conditions (i) and (iii) of Proposition \reff{LM-P},
$$
M:=M(\Gc)=M(\Lambda)\le 3~~~\text{and}~~~ A:=\sup\{e(\LM(G\cap V)): G\in\Gc\}\le C(m,p).
$$
\par Now applying Proposition \reff{EXT-DM} to $U,V$ and $\Gc$ defined by \rf{DF-UVG} we prove that
for every $m\ge 1$, $p>2$, and every $f\in\LM(\Wc)$ there  exists a function $F\in\LMPO$ linearly depending on $f$ such that
$$
F|_{\Wc}=f~~~\text{and}~~~\|F\|_{\LMPO}\le C(m,p)\,\|f\|_{\LM(\Wc)}.
$$
\par The proofs of Theorem \reff{WP-OM-EXT} and Theorem \reff{WP-EXT} are complete.\bx
\bigskip
\par We finish the section with the following useful consequence of Theorem \reff{WP-EXT} and  Theorem \reff{T-CM}.
\begin{corollary}\lbl{ARC-WP} Let $\Omega$ be a simply connected bounded domain in $\RT$ satisfying condition \rf{S-EXT}. Then for every $\x,\y\in\Omega$ and every ``Wide Path'' $\WPT^{(\x,\y)}$ joining $\x$ to $\y$ in $\Omega$ the following condition is satisfied: for every $a,b\in\WPT^{(\x,\y)}$ there exists a path $\gamma$ connecting $a$ to $b$ in $\WPT^{(\x,\y)}$ such that
$$
\diam\gamma\le\eta_W\|a-b\|.
$$
Here $\eta_W$ is a positive constant satisfying the inequality $\eta_W\le C(m,p)\,\CE$ where $\CE$ is the parameter from condition \rf{S-EXT}.
\end{corollary}
\bigskip
\SECT{5. ``The Narrow Path''}{5}
\setcounter{equation}{0}
\addtocontents{toc}{5. ``The Narrow Path''.\hfill \thepage\par\VST}
\indent\par {\bf 5.1. ``The Narrow Path'' construction algorithm.}
\addtocontents{toc}{~~~~5.1. ``The Narrow Path'' construction algorithm.\hfill \thepage\par}
Let $\x,\y\in\Omega$ and let $\WPT^{(\x,\y)}$ be ``The Wide Path'' joining $\x$ to $\y$ in $\Omega$ which we have constructed in the preceding section. We also recall that the domain $\Omega$ satisfies condition \rf{S-EXT}.
\par In this section we construct a ``Narrow Path'' described in Section 1, and present its main geometrical and Sobolev extension properties.
\par We begin with the following important
\begin{lemma}\lbl{3-SQ-L} Let $\ve>0$. Let $K,K_1$ and $K_2$ be pairwise disjoint squares in $\RT$ such that $K^{\cl}\cap K_1^{\cl}\ne\emp$, $K^{\cl}\cap K_2^{\cl}\ne\emp$, and $\#(K_1^{\cl}\cap K_2^{\cl})\le 1$.
\par Then there exists a square $\tK\subset K$ such that
$\tK^{\cl}\cap K_1^{\cl}\ne\emp$, $\tK^{\cl}\cap K_2^{\cl}\ne\emp$ and
$$
\diam\tK\le 2\dist(K_1,K_2)~~\text{whenever}~~K_1^{\cl}\cap K_2^{\cl}=\emp,
$$
and
\bel{F-IN2}
\diam\tK=\ve~~~
\text{whenever}~~~K_1^{\cl}\cap K_2^{\cl}\ne\emp.
\ee
\par Furthermore, for every $j\in\{1,2\}$ the following is true:
\bel{N-IN1}
\text{if}~~~\#(K_j^{\cl}\cap K^{\cl})>1~~~\text{then}~~~
\#(K_j^{\cl}\cap \tK^{\cl})>1\,.
\ee
\end{lemma}
\par {\it Proof.} First prove the lemma whenever $K_1^{\cl}\cap K_2^{\cl}=\emp$.
\par We begin with the following statement: for every $a,b\in K^{\cl}$ there exists a square $K_{a,b}$ such that
\bel{Q-CL}
a,b\in K_{a,b}\subset K^{\cl}
\ee
and
\bel{D-KAB}
\diam K_{a,b}=\|a-b\|\,.
\ee
(Recall that we measure distances in the uniform metric.)
\par Let $K=(y',z')\times(y'',z'')$. Hence, $|y'-z'|=|y''-z''|=\diam K$. Let $a=(a_1,a_2)$ and $b=(b_1,b_2)$. We may assume that
\bel{A-B}
|a_1-b_1|\le |a_2-b_2|=\|a-b\|\,.
\ee
Since $a,b\in K$, we have $[a_1,b_1]\subset[y',z']$ and
\bel{A2-B2}
[a_2,b_2]\subset[y'',z''].
\ee
Since
$$
\|a-b\|\le \diam K=|y'-z'|,
$$
by \rf{A-B},
$$
|a_1-b_1|\le \|a-b\|\le |y'-z'|.
$$
Hence there exists a closed interval $[a_1',b_1']$ such that
$$
|a'_1-b'_1|=|a_2-b_2|=\|a-b\|
$$
and
\bel{I-AB}
[a_1,b_1]\subset[a'_1,b'_1]\subset [y',z'].
\ee
\par Let $K_{a,b}:=(a'_1,b'_1)\times(a_2,b_2)$. Then, by \rf{A2-B2} and \rf{I-AB}, inclusions \rf{Q-CL} hold. Furthermore, by \rf{A-B},
$$
\diam K_{a,b}=|a_2-b_2|=\|a-b\|
$$
proving that $K_{a,b}$ satisfies \rf{Q-CL} and \rf{D-KAB}.
\par Note that the requirements $K^{\cl}\cap K_1\ne\emp$ and $K^{\cl}\cap K_2\ne\emp$ imply the following equality:
\bel{DSI}
\dist(K_1,K_2)=\dist(K_1^{\cl}\cap K^{\cl},K_2^{\cl}\cap K^{\cl})\,.
\ee
A proof of this simple geometrical fact we leave to the reader as an easy exercise.
\par Let $[u_1,v_1]:=K_1^{\cl}\cap K^{\cl}$ and
$[u_2,v_2]:=K_2^{\cl}\cap K^{\cl}$. By \rf{DSI}, there exist points $a'\in [u_1,v_1]$ and $b'\in [u_2,v_2]$ such that
$$
\|a'-b'\|=\dist(K_1,K_2).
$$
Let $a:=a'$ whenever $u_1=v_1$, and let
\bel{C-AF}
a~~\text{be a point from}~~(u_1,v_1)~~\text{such that}~~
\|a'-a\|\le\tfrac12\dist(K_1,K_2)
\ee
whenever $u_1\ne v_1$. In a similar way we define a point $b$ by letting $b:=b'$ whenever $u_2=v_2$, and
\bel{C-BF}
b~~\text{be a point from}~~(u_2,v_2)~~\text{such that}~~
\|b'-b\|\le\tfrac12\dist(K_1,K_2)
\ee
provided $u_2\ne v_2$.
\par Let $\tK=K_{a,b}$ be the square satisfying \rf{Q-CL} and \rf{D-KAB}. Then $a,b\in\tK\subset K^{\cl}$ and
\be
\diam\tK&=&
\|a-b\|\le\|a-a'\|+\|a'-b'\|+\|b'-b\|\nn\\
&\le& \tfrac12\dist(K_1,K_2)+\dist(K_1,K_2)+\tfrac12\dist(K_1,K_2)
=2\dist(K_1,K_2)\,.\nn
\ee
\par Furthermore, by \rf{C-AF} and \rf{C-BF}, the square $\tK$ satisfies \rf{N-IN1}.\medskip
\par It remains to prove the statement of the lemma whenever $K_1^{\cl}\cap K_2^{\cl}$ is a singleton, see \rf{F-IN2}. Thus $\{a\}=K_1^{\cl}\cap K_2^{\cl}$ fore some $a\in\RT$. Since $K_1,K_2$ and $K$ are pairwise disjoint squares with sides parallel to the coordinate axes, the point $a$ is a common vertex of these squares.
See Figure 20.
\begin{figure}[h]
\center{\includegraphics[scale=1]{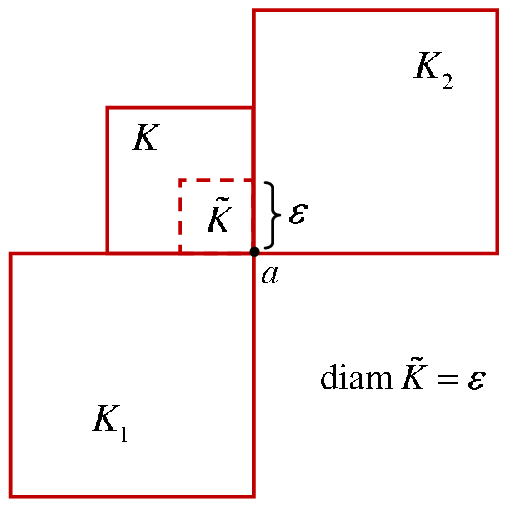}}
\caption{~}
\end{figure}
\medskip
\par This enables us to define the square $\tK:=K$ as follows: $\tK$ is a (unique) subsquare of $K$ with the vertex $a$ and $\diam\tK:=\ve$ as it shown on Figure 20. Clearly, $\tK$ satisfies conditions \rf{F-IN2} and \rf{N-IN1}.
\par The proof of the lemma is complete.\bx\smallskip
\par We are also needed the following auxiliary result.
\begin{lemma}\lbl{ADD-T} Let $1\le m\le k-3$. Let $S_{m+1}$ be a rotation square and let $a_{m+1}$ be a rotation point associated with the square $S_{m+1}$, see Definition \reff{R-SQ-P}. (We recall that in this case $S^{\cl}_m\cap S^{\cl}_{m+1}\cap S^{\cl}_{m+2}=\{a_{m+1}\}$.)
\par Let $H$ be a square such that $H\subset S_{m+1}$, the point $a_{m+1}$ is a vertices of $H$, and
\bel{U-DM}
\diam H\le\tfrac12\min\{\diam S_m,\diam S_{m+1},\diam S_{m+2}\}.
\ee
\par Then $H^{\cl}\cap S_{m+3}^{\cl}=\emp$.
\end{lemma}
\par {\it Proof.} First prove that
\bel{IM-Q3}
H^{\cl}\setminus\{a_{m+1}\}\subset S_m\cup S_{m+1}\cup S_{m+2}.
\ee
In fact, by \rf{U-DM},
$$
H^{\cl}\setminus\{a_{m+1}\}=H\cup ((H^{\cl}\cap S_m^{\cl})\setminus\{a_{m+1}\})\cup ((H^{\cl}\cap S_{m+2}^{\cl})\setminus\{a_{m+1}\}).
$$
\par We know that $H\subset S_{m+1}$. Recall that
$$
[u_m,v_m]=S_{m}^{\cl}\cap S_{m+1}^{\cl}~~~\text{and}~~~ [u_{m+1},v_{m+1}]=S_{m+1}^{\cl}\cap S_{m+2}^{\cl}\,.
$$
See \rf{UV-I}. (We can assume that $u_m=u_{m+1}=a_{m+1}$.) We also know that $T_m=(u_m,v_m)\subset\Omega$ and $T_{m+1}=(u_{m+1},v_{m+1})\subset\Omega$, see \rf{T-I} and \rf{RW-E}. These properties and inequality \rf{U-DM} imply the following:
$$
(H^{\cl}\cap S_m^{\cl})\setminus\{a_{m+1}\}\subset(u_m,v_m)\subset\Omega
~~\text{and}~~(H^{\cl}\cap S_{m+2}^{\cl})\setminus\{a_{m+1}\}\subset(u_{m+1},v_{m+1})
\subset\Omega,
$$
proving \rf{IM-Q3}.
\par By \rf{IM-Q3} an Lemma \reff{AD-1}, $(H^{\cl}\setminus\{a_{m+1}\})\cap S_{m+3}^{\cl}=\emp$. On the other hand, $a_{m+1}\in S_m^{\cl}$ and, by \rf{SIJ-L}, $S_m^{\cl}\cap S_{m+3}^{\cl}=\emp$. Hence $a_{m+1}\notin S_{m+3}^{\cl}$. Thus $H^{\cl}\cap S_{m+3}^{\cl}=\emp$, and the proof of the lemma is complete.\bx
\medskip
\par We turn to constructing ``The Narrow Path''. Let $\Sc_\Omega(\x,\y)=\{S_1, S_2,...,S_k\}$ be the family of squares constructed in ``The Wide Path Theorem'' \reff{W-PATH}.
\begin{proposition}\lbl{3-SQ} Let $k>2$. There exists a family
$$
\Qc_\Omega(\x,\y)=\{Q_1, Q_2,...,Q_k\}
$$
of pairwise disjoint squares such that:\medskip
\par (1). $Q_1=S_1$, $Q_k=S_k$, and $Q_i\subset S_i$ for every $i, 1\le i\le k$. Furthermore, $\x$ is the center of $Q_1$. In turn, $\y\in Q_k^{\cl}$ and $\dist(\y,Q_{k-1})=\diam Q_k$\,;
\smallskip
\par (2). $Q_i^{\cl}\cap Q_{i+1}^{\cl}\ne\emp$ for every $i\in\{1,...,k-1\}$, and $\#(Q_i^{\cl}\cap Q_{i+2}^{\cl})\le 1$ for every $i\in\{1,...,k-2\}$. Furthermore,
$$
\,Q_i^{\cl}\cap Q_{j}^{\cl}=\emp~~~\text{for all}~~~ i,j\in\{1,...,k\},~|i-j|>2\,;
$$
\par (3). If $\#(S_i^{\cl}\cap S_{i+1}^{\cl})>1$, then
$\#(Q_i^{\cl}\cap Q_{i+1}^{\cl})>1$. In turn, if $\#(S_i^{\cl}\cap S_{i+1}^{\cl})=1$, then
$\#(Q_i^{\cl}\cap Q_{i+1}^{\cl})=1$ as well;
\smallskip
\par (4). Let $1\le i\le k-2$. Then
\bel{DQ-1}
\diam Q_{i+1}\le 2\dist(Q_{i},Q_{i+2})~~~\text{if}~~~
Q_{i}^{\cl}\cap S_{i+2}^{\cl}=\emp,
\ee
and
\bel{DQ-2}
\diam Q_{i+1}\le \,
\tfrac14\min\{\diam Q_{i},\diam Q_{i+2}\}~~~\text{if}~~~
Q_{i}^{\cl}\cap S_{i+2}^{\cl}\ne\emp\,;
\ee
\par (5). If $Q_{i}^{\cl}\cap S_{i+2}^{\cl}\ne\emp$ then $Q_{i+1}^{\cl}\cap S_{i+3}^{\cl}=\emp$, $1\le i\le k-3$.
\smallskip
\end{proposition}
\par See Figure 2.\medskip
\par {\it Proof.} We obtain the family $\Qc_\Omega(\x,\y)$ as a result of a $k$ step inductive procedure based on Lemma \reff{3-SQ-L}. This procedure depends on a certain parameter $\ve_0>0$ which we define as follows. Let
$$
J:=\{m\in\{1,...,k-3\}: S_m^{\cl}\cap S_{m+1}^{\cl}\cap S_{m+2}^{\cl}\ne\emp\}.
$$
Thus for every $m\in J$ the square $S_{m+1}$ is a rotation square, see Definition \reff{R-SQ-P}. Let $a_{m+1}$ be the rotation point associated with $S_{m+1}$ so that
$$
\{a_{m+1}\}=S^{\cl}_m\cap S^{\cl}_{m+1}\cap S^{\cl}_{m+2}.
$$
We let $H_m$ denote a subsquare of $S_{m+1}$ such that $a_{m+1}$ is a vertices of $H_m$ and  
\bel{H-DM}
\diam H_m:=\tfrac14\min\{\diam S_m,\diam S_{m+1},\diam S_{m+2}\}.
\ee
\par Let
\bel{VE-AD}
\ve_0:=\tfrac14\min\{\dist(H_m,S_{m+3}): m\in J\}.
\ee
By Lemma \reff{ADD-T}, $\dist(H_m,S_{m+3})>0$ for every $m\in J$ so that $\ve_0>0$.\medskip
\par We are in a position to define the family of squares $\Qc_\Omega(\x,\y)$. At the first step of our inductive procedure we put $Q_1:=S_1$ and turn to the second step. We know that
$$
Q_1^{\cl}\cap S_2^{\cl}\ne\emp,~~S_2^{\cl}\cap S_3^{\cl}\ne\emp~~~\text{and}~~~\#(Q_1^{\cl}\cap S_3^{\cl})\le 1,
$$
see part (ii) of Lemma \reff{H-MS}. We put
$$
\ve:=\tfrac14\min\{\diam Q_1,\diam S_{2},\diam S_{3},\ve_0\}
$$
and apply Lemma \reff{3-SQ-L} to $\ve$ and pairwise disjoint squares $K_1:=Q_1,K:=S_2,$ and $K_3:=S_3$. By this lemma, there exists a square $\tK$ such that $\tK\subset S_2$,
$$
\tK^{\cl}\cap Q_1^{\cl}\ne\emp~~~\text{and}~~~
\tK^{\cl}\cap S_3^{\cl}\ne\emp.
$$
Furthermore,
$$
\diam \tK\le 2\dist(Q_1,S_3)~~~\text{if}~~~
Q_1^{\cl}\cap S_3^{\cl}=\emp,
$$
and
$$
\diam\tK=\ve~~~\text{if}~~~Q_1^{\cl}\cap S_3^{\cl}\ne\emp.
$$
In addition, if $\#(Q_1^{\cl}\cap S_2^{\cl})>1$, then
$\#(Q_1^{\cl}\cap \tK^{\cl})>1$, and if $\#(Q_1^{\cl}\cap S_2^{\cl})=1$, then $\#(Q_1^{\cl}\cap \tK^{\cl})=1$ as well. The same is true for the squares $S_3$ and $S_2$, i.e.,
$$
\text{if}~~~\#(S_2^{\cl}\cap S_3^{\cl})>1~~\text{then}~~
\#(\tK^{\cl}\cap S_3^{\cl})>1,
$$
and
$$
\text{if}~~~\#(S_2^{\cl}\cap S_3^{\cl})=1~~\text{then}~~
\#(\tK^{\cl}\cap S_3^{\cl})=1.
$$
\smallskip
\par We put $Q_2:=\tK$ and turn to the third step. We know that $Q_2^{\cl}\cap S_3^{\cl}\ne\emp$, $S_3^{\cl}\cap S_4^{\cl}\ne\emp$ and $\#(Q_2^{\cl}\cap S_4^{\cl})\le 1$ (because $Q_2\subset S_2$ and, by part (ii) of Lemma \reff{H-MS}, $\#(S_2^{\cl}\cap S_4^{\cl})\le 1$). This enables us to apply Lemma \reff{3-SQ-L} to
$$
\ve:=\tfrac14\min\{\diam Q_2,\diam S_{3},\diam S_{4},\ve_0\}
$$
and pairwise disjoint squares $K_1:=Q_2,K:=S_3$ and $K_3:=S_4$, and in this way to obtain a square $Q_3$, etc.
\par In a similar way we turn from the $m$-th step of this algorithm to its $(m+1)$-th step provided $1\le m<k-1$. After $m$ steps of this procedure we obtain a collection of  squares $\{Q_1,Q_2,...,Q_{m}\}$.
We know that $Q_{m}\subset S_{m}$,
$\,Q_{m}^{\cl}\cap S_{m+1}^{\cl}\ne\emp$, $\,S_{m+1}^{\cl}\cap S_{m+2}^{\cl}\ne\emp$ and $\,\#(Q_{m}^{\cl}\cap S_{m+2}^{\cl})\le 1$ (because $Q_{m}\subset S_{m}$ and, by part (ii) of Lemma \reff{H-MS}, $\#(S_{m}^{\cl}\cap S_{m+2}^{\cl})\le 1)$. We put
$$
\ve:=\tfrac14\min\{\diam Q_m,\diam S_{m+1},\diam S_{m+2},\ve_0\},
$$
$K_1:=Q_{m},K:=S_{m+1}$ and $K_2:=S_{m+2}$. Clearly, $K_1,K,K_2$ is a triple of pairwise disjoint squares satisfying the hypothesis of Lemma \ref{3-SQ-L}.
\par By this lemma, there exists a square $Q_{m+1}=\tK$ such that $Q_{m+1}\subset S_{m+1}$,
\bel{IN-QM}
Q_{m}^{\cl}\cap Q_{m+1}^{\cl}\ne\emp~~~\text{and}~~~
Q_{m+1}^{\cl}\cap S_{m+2}^{\cl}\ne\emp.
\ee
Furthermore,
\bel{Q-MIN}
\diam Q_{m+1}\le 2\dist(Q_{m},S_{m+2})~~~\text{whenever}~~~
Q_{m}^{\cl}\cap S_{m+2}^{\cl}=\emp.
\ee
\par See Figure 21.\bigskip
\begin{figure}[h]
\center{\includegraphics[scale=0.9]{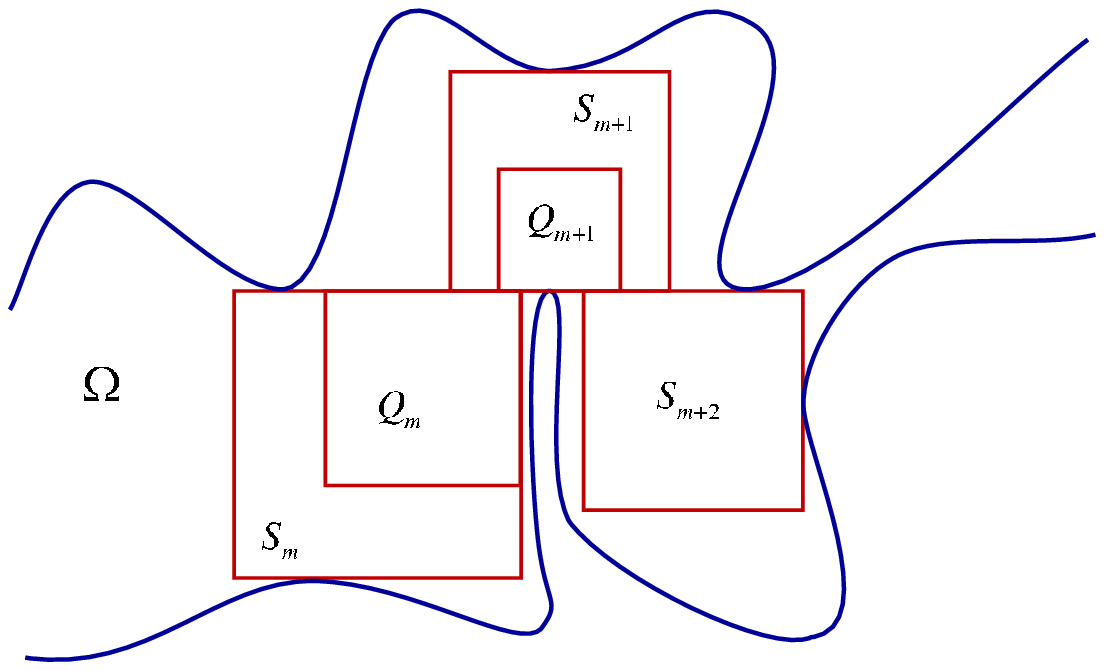}}
\caption{~}
\end{figure}
\medskip
\par Now let
\bel{Q-2R}
Q_{m}^{\cl}\cap S_{m+2}^{\cl}\ne\emp.
\ee
Since $Q_m\subset S_m$, we have
$S_{m}^{\cl}\cap S_{m+2}^{\cl}\ne\emp$, so that the square $S_{m+1}$ is a {\it rotation square} of ``The Wide Path'' $\Wc$, see Definition \reff{R-SQ-P}. We know that in this case the intersection $S_{m}\cap S_{m+1}\cap S_{m+2}$ is the {\it rotation point} $\{a_{m+1}\}$ associated with the rotation square $S_{m+1}$.
\par By Lemma \ref{3-SQ-L}, in this case
\bel{Q-O2}
\diam Q_{m+1}=\ve=\,\tfrac14
\min\{\diam Q_{m},\diam S_{m+1},\diam S_{m+2},\ve_0\}.
\ee

\par In addition, by \rf{N-IN1},
\bel{QM-IN1}
\text{if}~~~\#(Q_{m}^{\cl}\cap S_{m+1}^{\cl})>1~~\text{then}~~
\#(Q_{m}^{\cl}\cap Q_{m+1}^{\cl})>1,
\ee
and
\bel{QM-IN2}
\text{if}~~~\#(S_{m+1}^{\cl}\cap S_{m+2}^{\cl})>1~~\text{then}~~
\#(Q_{m+1}^{\cl}\cap S_{m+2}^{\cl})>1.
\ee
See Figure 22.
\begin{figure}[h]
\center{\includegraphics[scale=0.79]{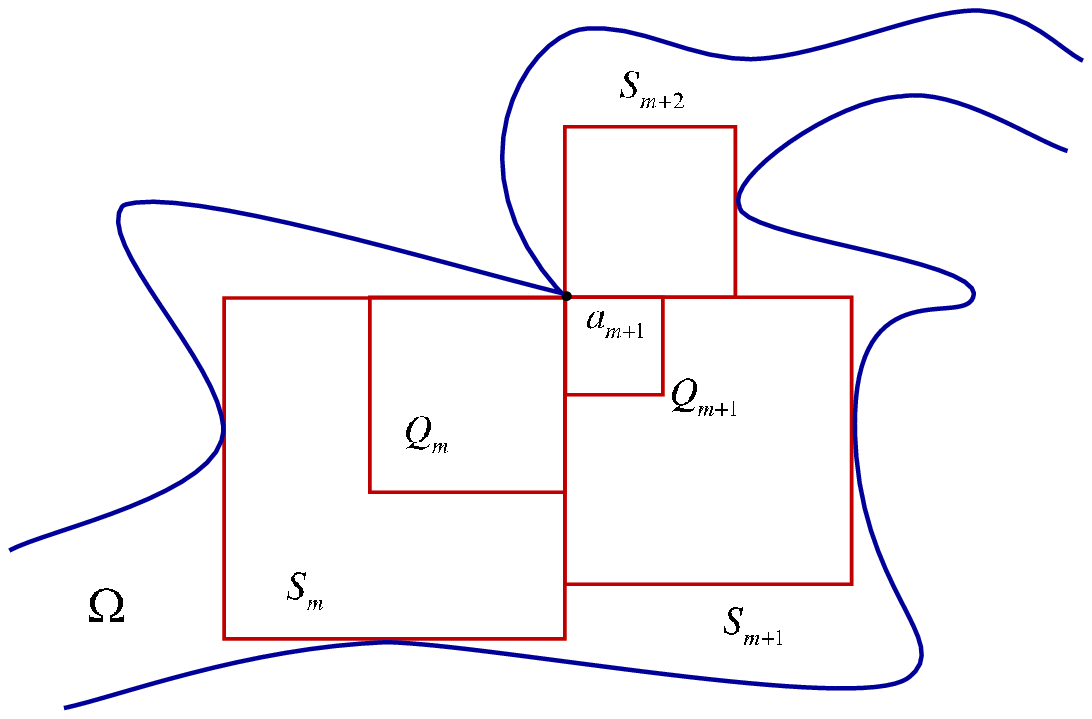}}
\caption{~}
\end{figure}
\par After $(k-1)$ steps of this algorithm we obtain a family of squares $\{Q_1,...,Q_{k-1}\}$. Finally, at the last step of this procedure we put $Q_k:=S_k$ and stop.\medskip
\par Let us prove that the obtained family $\{Q_1,...,Q_k\}$ of squares possesses properties (1)-(5) of the proposition.
\par Since $Q_1=S_1$, $Q_k=S_k$ and $Q_{m+1}\subset S_{m+1}$, the first part of property (1) holds. The second and third parts follow from part (a) of Lemma \reff{MAIN-WP}. Property (2) of the proposition follows from \rf{IN-QM} and part (ii) of Lemma \reff{H-MS}. Property (3) directly follows from properties \rf{QM-IN1} and \rf{QM-IN2} and the inclusion $Q_i\subset S_i$, $1\le i\le k$.

\par Prove property (4). Suppose that $Q_{m}^{\cl}\cap S_{m+2}^{\cl}=\emp$. Since $Q_{m+2}\subset S_{m+2}$, by \rf{Q-MIN},
$$
\diam Q_{m+1}\le 2\dist(Q_{m},S_{m+2})
\le 2\dist(Q_{m},Q_{m+2})
$$
proving \rf{DQ-1}.\smallskip
\par Now prove \rf{DQ-2} whenever $Q_{m}^{\cl}\cap S_{m+2}^{\cl}\ne\emp$, i.e., \rf{Q-2R} holds. In this case $Q_{m+1}$ is a subsquare of $S_{m+1}$, the rotation point $a_{m+1}$ is a vertices of $Q_{m+1}$ and its diameter is given by \rf{Q-O2}. See Figure 22.
\par Comparing $Q_{m+1}$ with the square $H_{m+1}$ defined at the beginning of the proof (see \rf{H-DM}) we conclude that $Q_{m+1}\subset H_{m+1}$. Note that, by \rf{VE-AD} and \rf{Q-O2},
$$
\diam Q_{m+1}\le \ve_0\le \tfrac14\dist(H_m,S_{m+3}).
$$
On the other hand, we know that $Q^{\cl}_{m+2}\cap Q^{\cl}_{m+1}\ne\emp$ and $Q^{\cl}_{m+2}\cap S^{\cl}_{m+3}\ne\emp$ so that
$$
\diam Q_{m+2}\ge \dist(Q_{m+1},S_{m+3})\ge \dist(H_{m+1},S_{m+3})\ge 4\ve_0.
$$
Hence $\diam Q_{m+1}\le \frac14\diam Q_{m+2}$.
\par In addition, by \rf{Q-O2}, $\diam Q_{m+1}\le \frac14\diam Q_{m}$ proving \rf{DQ-2} and property (4) of the proposition.\smallskip
\par Prove property (5). Suppose that \rf{Q-2R} holds so that $a_{m+1} $ is a rotation point associated with the rotation square $S_{m+1}$. See Figure 22.
\par Let $H:=Q_{m+1}$. Then $H\subset S_{m+1}$, the point $a_{m+1}$ is a vertices of $H$, and, by \rf{Q-O2}, inequality \rf{U-DM} of Lemma \reff{ADD-T} is satisfied.
By this lemma, $H^{\cl}\cap S_{m+3}^{\cl}=\emp$ proving  that $Q_{m+1}^{\cl}\cap S_{m+3}^{\cl}=\emp$. This implies property (5).
\par The proof of the proposition is complete.\bx
\bigskip
\par {\bf 5.2. Main geometrical properties of ``The Narrow Path''.}
\addtocontents{toc}{~~~~5.2. Main geometrical properties of ``The Narrow Path''. \hfill \thepage\par}
We recall that ``The Narrow Path'' $\NPT^{(\x,\y)}$ joining $\x$ to $\y$ in $\Omega$ is defined by formula \rf{NP}:
$$
\NPT^{(\x,\y)}:=\left(\,\bigcup_{i=1}^k \left(\,Q_i^{\cl}\,\usm\, \hS_i\,\right)\right)^{\circ}.
$$
Recall that $\{\hS_1,...,\hS_k\}$ is the family of sets (more specifically, squares or empty sets) introduced in  Definition \reff{DF-SH}.\smallskip
\par Let us present several useful geometrical properties of ``The Narrow Path'' which we will use later on in the study of the extension properties of $\NPT^{(\x,\y)}$ and differential properties of the ``\FG growing'' functions.
\begin{lemma}\lbl{H-QP} (i). $\hS_i=\emp$ whenever $i=k$ or $\#(Q_i^{\cl}\cap Q_{i+1}^{\cl})>1$, and
$$
\hS_i=S(w_i,\hdl)~~~\text{if}~~~
\#(Q_i^{\cl}\cap Q_{i+1}^{\cl})=1.
$$
Here $\{w_i\}=S_i^{\cl}\cap S_{i+1}^{\cl}$, see \rf{WI-DF}, and $\hdl$ is the number defined by \rf{DLT-FIN}.\medskip
\par (ii). $\{\hS_1,...,\hS_k\}$ is a family of pairwise disjoint subsets of $\Omega$ such that
\bel{QH-IN}
(\hS_i^{\cl})\cap Q_j^{\cl}=\emp~~\text{for every}~~ 1\le i,j\le k,\, j\ne i,i+1.
\ee
\par (iii). $\diam \hS_i\le \frac14\, \min\{\diam Q_i,\diam Q_{i+1}\}$ for every $i, 1\le i\le k-1$.
\end{lemma}
\par {\it Proof.} By part (3) of Proposition \reff{3-SQ},
$$
\#(Q_i^{\cl}\cap Q_{i+1}^{\cl})=1~~~\text{if and only if}~~~\#(S_i^{\cl}\cap S_{i+1}^{\cl})=1.
$$
This property, Definition \reff{DF-SH} (see\rf{SH-1} and \rf{SH-ND}) imply part (i) of the lemma.
\par Prove (ii). By Lemma \reff{H-IJ}, the sets of the family $\{2\hS_i:1=1,...,k\}$ are pairwise disjoint subsets of $\Omega$ so that the family $\{\hS_i:1=1,...,k\}$ consists of  pairwise disjoint subsets of $\Omega$ as well. This property, \rf{SH-JN} and the inclusion $Q_i\subset S_i$ immediately imply the statement \rf{QH-IN} proving (ii).
\par Prove (iii). Suppose that $\hS_i\ne\emp$, i.e., by part (i) of the present lemma,
$Q_i^{\cl}\cap Q_{i+1}^{\cl}=\{w_i\}$. (Recall that $w_i$ is the center of the square $\hS_i$.) Thus $w_i\in Q_i^{\cl}$.
\par If $i=1$ then, by part (1) of Proposition \reff{3-SQ}, $Q_1=S_1$. In turn, by \rf{DLT-FIN},
$$
\hdl\le \tfrac18\, \hdl_2=\min\{\diam S_j:1\le j\le k\}
$$
so that
\bel{DM-DQ}
\diam\hS_i=2\,\hdl\le \tfrac14 \diam S_j~~~\text{for every}~~~1\le i,j\le k.
\ee
In particular, $\diam\hS_1\le \tfrac14 \diam S_1=\tfrac14 \diam Q_1$.
\par Now let $i>1$. Since $Q_{i-1}^{\cl}\cap Q_i^{\cl}\ne\emp$ and $Q_{i-1}\subset S_{i-1}$, we have
$\dist(w_i,S_{i-1})\le \diam Q_i$. But, by \rf{DLT-3} and \rf{DLT-FIN},
$$
\diam \hS_i=2\hdl\le 2\cdot\tfrac18\,\hdl_3\le \tfrac14\dist(w_i,S_{i-1})\le \tfrac14\diam Q_i.
$$
\par In the same way we prove that $\diam \hS_i\le \tfrac14\diam Q_{i+1}$. In fact, let $i<k-1$. Since $Q_{i+1}^{\cl}\cap S_{i+2}^{\cl}\ne\emp$, we have $\dist(w_i,S_{i+1})\le \diam Q_{i+1}$. Hence,
$$
\diam \hS_i=2\hdl\le 2\cdot\tfrac18\,\hdl_3\le \tfrac14\dist(w_i,S_{i+1})\le \tfrac14\diam Q_{i+1}.
$$
\par If $i=k-1$ then, by part (1) of Proposition \reff{3-SQ}, $Q_{i+1}=Q_k=S_k$ so that, by \rf{DM-DQ},
$\diam\hS_{k-1}\le \tfrac14 \diam S_k=\tfrac14 \diam Q_{k}$ proving part (iii) and the lemma.\bx
\par The next proposition is an analog of Proposition \reff{W-G1} for ``The Narrow Path''. Its proof literally follows the scheme of the proof of Proposition \reff{W-G1};  we leave the details for the interested reader.
\begin{proposition}\lbl{NP-REP} ``The Narrow Path'' $\Nc:=\NPT^{(\x,\y)}$ is an open connected subset of the domain $\Omega$ which has the following representation:
\bel{R-NP}
\Nc=\bigcup_{i=1}^{k-1}\left(\,Q_i^{\cl}\,\usm \, Q_{i+1}^{\cl}\,\usm \,\hS_i\,\right)^{\circ}\,.
\ee
\end{proposition}
\par Let $Q_i$ and $Q_{i+1}$, $1\le i<k$, be two subsequent squares from ``The Narrow Path'' such that $\#(Q_i^{\cl}\cap Q_{i+1}^{\cl})>1$. Since $Q_i$ and $Q_{i+1}$ are touching squares, intersection of their closures is a line segment. We denote the ends of this segment by $s_i$ and $t_i$. Thus  
\bel{ST-DF}
[s_i,t_i]:=Q_i^{\cl}\cap Q_{i+1}^{\cl}~~~\text{whenever}~~~\#(Q_i^{\cl}\cap Q_{i+1}^{\cl})>1,
\ee
so that in this case
\bel{QST-I}
(Q_i^{\cl}\cup Q_{i+1}^{\cl})^\circ=Q_i\cup Q_{i+1}\cup(s_i,t_i).
\ee
\par Let $1\le i< k$ and let
\bel{Y-I}
Y_i:=\left\{
\begin{array}{ll}
\hS_i,&~~\text{if}~~
\#(Q_i^{\cl}\cap Q_{i+1}^{\cl})=1,\\\\
(s_i,t_i),&~~\text{if}~~
\#(Q_i^{\cl}\cap Q_{i+1}^{\cl})>1.
\end{array}
\right.
\ee
We also put $Y_k:=\emp$.
\par Then, by \rf{QST-I} and by definition of $\hS_i$, see \rf{SH-ND},
\bel{Q-BG}
\left(\,Q_i^{\cl}\,\usm \, Q_{i+1}^{\cl}\,\usm \,\hS_i\,\right)^{\circ}=Q_i\,\usm \, Q_{i+1}\,\usm \,Y_i\,.
\ee
Combining this equality with \rf{R-NP} we obtain the following representation of ``The Narrow Path'':
\bel{RN-E}
\Nc=\NPT^{(\x,\y)}=\usm\limits_{i=1}^k \left(Q_i\,\usm\,Y_i\right)\,.
\ee
\medskip
\par In the next two lemmas we present additional geometrical properties of ''The Narrow Path''.
\begin{lemma}\lbl{C-KZ} (i) Let $1\le i\le k-2$ and let $Q_i^{\cl}\cap S_{i+2}^{\cl}=\emp$. Then
\bel{DL-V}
\diam Q_{i+1}\le 4 \dist(Y_i,Y_{i+1}).
\ee
\par Furthermore,
\bel{DL-VX}
\diam Q_1\le 4 \dist(\x,Y_1)~~~\text{and}~~~
\diam Q_k\le 4 \dist(\y,Y_{k-1}).
\ee
\end{lemma}
\par {\it Proof.} (i) Suppose that $\#(Q_i^{\cl}\cap Q_{i+1}^{\cl})=1$ so that $Y_i=\hS_i$. See Definition \reff{DF-SH}. Consider two cases.
\par {\it The first case: $\#(Q_{i+1}^{\cl}\cap Q_{i+2}^{\cl})=1$.} In this case $Y_{i+1}=\hS_{i+1}$. Recall that the center of the square  $\hS_{i+1}$, the point $w_{i+1}$, is a common vertex of the squares
$Q_{i+1}^{\cl}$ and $Q_{i+2}^{\cl}$. Furthermore, since $\hS_{i}\cap\hS_{i+1}=\emp$, we have $w_i\ne w_{i+1}$.
\par Since $w_i$ is a vertex of $Q_{i+1}^{\cl}$ as well, we conclude that
\bel{C-Q}
\|w_i-w_{i+1}\|=\diam Q_{i+1}\,.
\ee
By part (iii) of Lemma \reff{H-QP},
$$
\diam \hS_{i},\,\diam \hS_{i+1}\le \tfrac14\diam Q_{i+1}\,.
$$
Combining this inequality with \rf{C-Q}, we obtain that
$$
\dist(Y_{i},Y_{i+1})=\dist(\hS_{i},\hS_{i+1})\ge \tfrac12\diam Q_{i+1}\,.
$$
\par In the same fashion, basing on property (1) of Lemma \reff{3-SQ}, we prove inequalities \rf{DL-VX}.
\smallskip
\par {\it The second case: $\#(Q_{i+1}^{\cl}\cap Q_{i+2}^{\cl})>1$.} In this case $Y_{i+1}=(s_{i+1},t_{i+1})\subset Q_{i+1}^{\cl}\cap Q_{i+2}^{\cl}$. Since $Q_i^{\cl}\cap S_{i+2}^{\cl}=\emp$, by \rf{DQ-1},
$$
\diam Q_{i+1}\le 2\dist(Q_i,Q_{i+2})\le
2\dist(w_i,Y_{i+1})\,.
$$
\par On the other hand, by part (iii) of Lemma \reff{H-QP}, $\diam \hS_i\le\tfrac14\diam Q_{i+1}$. Therefore, for each $z\in \hS_i$ we have
$$
\dist(Y_{i+1},z)\ge \dist(Y_{i+1},w_i)-\|z-w_i\|\ge
\tfrac12\diam Q_{i+1}-\tfrac14\diam Q_{i+1}=
\tfrac14\diam Q_{i+1}
$$
proving \rf{DL-V} in the case under consideration. \smallskip
\par Let now $\#(Q_i^{\cl}\cap Q_{i+1}^{\cl})>1$ and $\#(Q_{i+1}^{\cl}\cap Q_{i+2}^{\cl})>1$. In this case $Y_i=(s_i,t_i)\subset Q_i^{\cl}$ and $Y_{i+1}=(s_{i+1},t_{i+1})\subset Q_{i+2}^{\cl}$. Hence, by \rf{DQ-1},
$$
\dist(\hS_{i},\hS_{i+1})=\dist(Y_{i},Y_{i+1})\ge \dist(Q_i,Q_{i+2})\ge \tfrac12 \diam Q_{i+1}.
$$
\par Consider the remaining case where $\#(Q_i^{\cl}\cap Q_{i+1}^{\cl})>1$ and $\#(Q_{i+1}^{\cl}\cap Q_{i+2}^{\cl})=1$, i.e., $Y_i=(s_i,t_i)\subset Q_i^{\cl}$ and $Y_{i+1}=\hS_{i+1}$.
\par Recall that $\hS_{i+1}=S(w_{i+1},\hdl)$ where
$\{w_{i+1}\}=Q_{i+1}^{\cl}\cap Q_{i+2}^{\cl}$ and $\hdl$ is defined by \rf{DLT-FIN}. In particular, by $\hdl\le\tfrac18\, \hdl_3\le \tfrac18 \dist(w_{i+1},S_{i})$. But $w_{i+1}\in Q^{\cl}_{i+1}$ and $Q_{i+1}^{\cl}\cap S_{i}^{\cl}\ne\emp$ so that $\dist(w_{i+1},S_{i})\le\diam Q_{i+1}$. Hence $\hdl\le\tfrac18 \diam Q_{i+1}$.
\par This inequality implies the following:
$$
\dist(Y_{i},Y_{i+1})=\dist(Y_{i},\hS_{i+1})\ge \dist(Y_{i},w_{i+1})-\hdl\ge      \dist(Q_i,Q_{i+2})-\tfrac18 \diam Q_{i+1}.
$$
On the other hand, by \rf{DQ-1},
$\dist(Q_i,Q_{i+2})\ge \tfrac12 \diam Q_{i+1}$
which proves \rf{DL-V} and the lemma.\bx
\bigskip
\par {\bf 5.3. Sobolev extension properties of ``The Narrow Path''.}
\addtocontents{toc}{~~~~5.3. Sobolev extension properties of ``The Narrow Path''. \hfill \thepage\par\VST}
\begin{lemma}\lbl{NP-C1} Let $1\le i<k-1$ and let
$a\in S_i, b\in S_{i+2}$. Then there exists $z\in Q_{i+1}\cup \hS_i\cup \hS_{i+1}$ such that
$$
\|z-a\|\le 3\,\eta_W\|a-b\|.
$$
Here $\eta_W$ is the constant from Corollary \reff{ARC-WP}.
\end{lemma}
\par {\it Proof.} By Corollary \reff{ARC-WP}, there exists a path $\gamma$ joining $a$ to $b$ in $\WPT^{(\x,\y)}$ such that
\bel{DM-GM}
\diam\gamma\le\eta_W\|a-b\|.
\ee
\par In turn, by part (i) of Lemma \reff{AD-1-2}, $\gamma\cap S_{i+1}\ne\emp$ so that there exists a point $\tz\in \gamma\cap S_{i+1}$. However we can not guarantee that $\tz\in Q_{i+1}$.
\par By Lemma \reff{CR-WP},
\bel{G-TIN}
\gamma\cap T_{i}\ne\emp~~~\text{and}~~~
\gamma\cap T_{i+1}\ne\emp
\ee
where $T_i$ is the set defined by \rf{T-I}.
\par Suppose that $\#(S_{i}^{\cl}\cap S_{i+1}^{\cl})=1$. In this case, by part (3) of Proposition \reff{3-SQ}, we have $\#(Q_{i}^{\cl}\cap Q_{i+1}^{\cl})=1$ as well. Recall also that in this case $T_i=\hS_i$, see \rf{T-I}, so that $\gamma\cap \hS_i\ne\emp$.
\par In the same way we show that  $\gamma\cap \hS_{i+1}\ne\emp$ provided $\#(S_{i+1}^{\cl}\cap S_{i+2}^{\cl})=1$. Thus there exists $z\in\gamma\cap (\hS_{i}\cup \hS_{i+1})$ whenever
$$
\text{either}~~~\#(S_{i}^{\cl}\cap S_{i+1}^{\cl})=1~~~\text{or}~~~\#(S_{i+1}^{\cl}\cap S_{i+2}^{\cl})=1.
$$
Since $a,z\in \gamma,$ by \rf{DM-GM},
$$
\|z-a\|\le \diam \gamma\le \eta_W\|a-b\|.
$$
\par Thus we can assume that
$$
\#(S_{i}^{\cl}\cap S_{i+1}^{\cl})>1~~~\text{and}~~~\#(S_{i+1}^{\cl}\cap S_{i+2}^{\cl})>1
$$
so that $T_i=(u_i,v_i)$ and $T_{i+1}=(u_{i+1},v_{i+1})$. See \rf{T-I}. In particular, we have the following: $T_i^{\cl}=[u_i,v_i]=S_i^{\cl}\cap S_{i+1}^{\cl}$, see \rf{UV-I}.
\par By \rf{G-TIN} there exist points $a'\in T_i$ and $b'\in T_{i+1}$. Clearly,
$$
a'\in\gamma\cap S_{i}^{\cl}\cap S_{i+1}^{\cl}
~~~\text{and}~~~
b'\in \gamma\cap S_{i+1}^{\cl}\cap S_{i+2}^{\cl}\,.
$$
\par Note that $T_i^{\cl}$ and $T_{i+1}^{\cl}$ are closed line segments which lie on $\partial S_i$. Since the squares $S_i, S_{i+1}$ and $S_{i+2}$ are pairwise disjoint, intersection of $T_i^{\cl}$ and $T_{i+1}^{\cl}$ contains at most one point, i.e.,
\bel{T-INT}
\#(T_i^{\cl}\cap T_{i+1}^{\cl})\le 1.
\ee
\par We also notice that, by part (2) of Proposition \reff{3-SQ}, $Q_i^{\cl}\cap Q_{i+1}^{\cl}\ne\emp$. But
$Q_i\subset S_i$ and $Q_{i+1}\subset S_{i+1}$ so that
$Q_i^{\cl}\cap S_{i}^{\cl}\cap S_{i+1}^{\cl}\ne\emp$
proving that
\bel{P-Q1}
Q_{i+1}^{\cl}\cap T_{i}^{\cl}\ne\emp.
\ee
In the same fashion we prove that
\bel{P-Q2}
Q_{i+1}^{\cl}\cap T_{i+1}^{\cl}\ne\emp.
\ee
\par To finish the proof of the lemma we are needed the following simple geometrical\medskip
\par {\it Statement C.} {\it Let $S$ be a square in $\RT$ and let $T'\subset\partial S$ and $T''\subset\partial S$ be closed line segments such that $\#(T'\cap T'')\le 1$. Let $Q\subset S$ be a square such that
\bel{W1}
Q^{\cl}\cap T'\ne\emp~~~\text{and}~~~
Q^{\cl}\cap T''\ne\emp.
\ee
\par Then for every $\ta\in T'\sm Q^{\cl}$, $\tb\in T''\sm Q^{\cl}$ and $z\in Q^{\cl}$ the following inequality
$$
\|z-\ta\|\le \|\ta-\tb\|
$$
holds.}\medskip
\par We prove this statement with the help of projection on the coordinate axes. This enables us to reduce Statement C to the following trivial assertion: {\it Let $I_1$ and $I_2$ be closed intervals in $\R$ such that $\#(I_1\cap I_2)\le 1$. Let $I\subset \R$ be a closed interval such that $I_1\cap I\ne\emp$ and $I_2\cap I\ne\emp$. Then $|c-c_1|\le |c_1-c_2|$ provided $c_1\in I_1\sm I$, $c_2\in I_2\sm I$ and $c\in I$.}\medskip
\par Now we finish the proof of the lemma as follows.
First we notice that, by \rf{DM-GM},
\bel{EW1}
\|a'-b'\|\le \diam\gamma\le\eta_W\|a-b\|.
\ee
\par Let $S:=S_{i+1}$, $Q:=Q_{i+1}$, and let $T':=T_{i}^{\cl}$, $T'':=T_{i+1}^{\cl}$. Then \rf{P-Q1} and  \rf{P-Q2} imply \rf{W1}, and \rf{T-INT} implies inequality $\#(T'\cap T'')\le 1$. Hence, by Statement C, for every $z\in Q=Q_{i+1}$
$$
\|z-a'\|\le \|a'-b'\|~~~\text{provided}~~~a'\notin Q^{\cl}=Q_{i+1}^{\cl}~~~\text{and}~~~b'\notin Q^{\cl}=Q_{i+1}^{\cl}.
$$
Combining this inequality with \rf{EW1} we obtain:
\bel{EW2}
\|z-a'\|\le 2\,\eta_W\|a-b\|.
\ee
\par If $a'\in Q_{i+1}^{\cl}$, then we choose $z\in Q_{i+1}$ such that $\|z-a'\|\le \|a-b\|$. In turn, if $b'\in Q_{i+1}^{\cl}$, we can choose $z\in Q_{i+1}$ for which $\|z-b'\|\le \|a-b\|$. This inequality and \rf{EW1} imply the following:
$$
\|z-a'\|\le \|z-b'\|+\|a'-b'\|\le \|a-b\|+\eta_W\|a-b\|\le 2\,\eta_W\|a-b\|.
$$
(Of course, we assume that $\eta_W\ge 1$.)
\par These estimates show that there always exists a point $z\in Q_{i+1}$ satisfying inequality \rf{EW2}. Finally, by \rf{EW2} and \rf{DM-GM}, we obtain that
$$
\|z-a\|\le \|z-a'\|+\|a'-a\|\le 2\eta_W\|a-b\|+\diam \gamma\le 3\eta_W\|a-b\|
$$
proving the lemma.\bx
\medskip
\par Let us introduce two families of open subsets of $\Omega$, a family $\Gc$ and a family $\Hc$, which control Sobolev extension properties of ``The Narrow Path''. We define the members of these families as follows: Let
\bel{AI-RG}
A_i:=\left(Q_i^{\cl}\,\usm\, Q_{i+1}^{\cl}\,\usm\, \hS_i\right)^{\circ},~~~i=1,...,k-1,
\ee
and let
\bel{ST-RG}
G_i:=A_i\,\usm\, A_{i+1}~~~i=1,...,k-2.
\ee
\par Note that, by \rf{R-NP} and \rf{AI-RG},
\bel{AI-R}
\Nc=\NPT^{(\x,\y)}=\bigcup_{i=1}^{k-1}\,A_i
\ee
so that
\bel{GI-NP}
\Nc=\bigcup_{i=1}^{k-2}\,G_i\,.
\ee
\par We also put
$$
B_i:=\left(S_i^{\cl}\cup Q_{i+1}^{\cl}\,\usm\, \hS_i\right)^{\circ},~~~
C_i:=\left(Q_{i+1}^{\cl}\,\usm\, S_{i+2}^{\cl}\cup \hS_{i+1}\right)^{\circ},~~~i=1,...,k-2,
$$
and, finally,
\bel{HI-DF}
H_i:=B_i\,\usm\, C_{i}~~~i=1,...,k-2.
\ee
\par Note several useful representations of $G_i$ and $H_i$ which easily follow from their definitions and
part (ii) of Lemma \reff{H-MS}. In particular,
\bel{GI-FB}
G_i:=\left(\,Q_i^{\cl}\,\usm \, Q_{i+1}^{\cl}\,\usm \,Q_{i+2}^{\cl}\,\usm \,\hS_i\,\usm \,\hS_{i+1}\right)^{\circ},~~~1\le i\le k-2\,.
\ee
In turn,
\bel{HI-R2}
H_i:=\left(\,S_i^{\cl}\,\usm \, Q_{i+1}^{\cl}\,\usm \,S_{i+2}^{\cl}\,\usm \,\hS_i\,\usm \,\hS_{i+1}\right)^{\circ},~~~1\le i\le k-2\,.
\ee
\par We use representation \rf{AI-R} to prove the following important property of ``The Narrow Path''.
\begin{lemma}\lbl{NP-FC} ``The Narrow Path'' $\Nc:=\NPT^{(\x,\y)}$ is a simply connected domain.
\end{lemma}
\par {\it Proof.} The statement of the lemma easily follows from \rf{AI-R} and the following properties of simply connected domains: {\it Let $G$ and $G'$ be two simply connected domains in $\RT$ with simply connected intersection $G\cap G'$. Then $G\cup G'$ is simply connected as well.} See, e.g., \cite{Rem}, p. 175.\medskip 
\par Let $1\le m\le k-1$, and let
\bel{UI-R}
U_m:=\bigcup_{i=1}^{m}\,A_i.
\ee
Thus, by \rf{AI-R}, $U_{k-1}=\Nc$. Note that, by Proposition \reff{3-SQ} and Lemma \reff{H-QP},
\bel{IN-UI}
A_m\cap A_{m+1}=Q_{m+1}~~~\text{for all}~~~m=1,...,k-2.
\ee
\par Prove that {\it each set $U_m$, $m=1,...,k-1$, is a simply connected domain}. We do this by induction on $m$. First we note that each set
$
A_i=\left(Q_i^{\cl}\,\cup\, Q_{i+1}^{\cl}\,\cup\, \hS_i\right)^{\circ}
$
is a simply connected planar domain. (The reader  can easily see that $\partial A_i$ is a connected set which guarantees that $A_i$ is simply connected. See, e.g., \cite{Bea}, p. 81.)
\par In particular, $U_1=A_1$ is a simply connected domain. Suppose that $1\le m\le k-2$ and the set $U_m$ is simply connected. Then, by \rf{UI-R}, $U_{m+1}=U_m\cup A_{m+1}$. Furthermore, by Proposition \reff{3-SQ} and Lemma \reff{H-QP}, $A_i\cap A_j\ne\emp$ if and only if $|i-j|\le 1$ proving that $U_m\cap A_{m+1}=A_m\cap A_{m+1}$. Hence, by \rf{IN-UI}, $U_m\cap A_{m+1}=Q_{m+1}$.
\par Thus $U_m$ and $A_{m+1}$ are two simply connected planar domains with simply connected intersection. Therefore, by the above statement, their union $U_m\cup A_{m+1}=U_{m+1}$ is a simply connected domain as well.
\par The proof of the lemma is complete.\bx\medskip
\par In Section 6 we will be needed the following important geometrical property of ``The Narrow Path''.
\begin{lemma}\lbl{C-KZ-2} Let $\gamma$ be a path joining $\x$ to $\y$ in ``The Narrow Path'' $\Nc=\NPT^{(\x,\y)}$. There exist points $s_n,t_n\in\gamma$, $1\le n\le k$, such that:\medskip
\par (1). $s_{1}=\x$, $t_k=\y$,
$$
s_n\in \gamma\cap Y_{n-1}^{\cl}~~\text{for all}~~2\le n\le k, ~~\text{and}~~t_n\in \gamma\cap Y_{n}^{\cl}~~\text{for all}~~ 1\le n\le k-1\,.
$$
\par (2). Let $\gamma_n$ be a subarc of $\gamma$ with the ends in $s_n$ and $t_n$, $1\le n\le k$. Then $\gamma_n\subset Q_n^{\cl}$.\medskip
\par (3). The sets of the family
$\{\gamma_n\sm\{s_{n},t_{n}\}:1\le n\le k\}$ are pairwise disjoint.
\end{lemma}
\par {\it Proof.} Let $1\le n\le k$ and let
$$
\Nc_n:=\usm\limits_{i=1}^n (Q_i\,\usm\, Y_i).
$$
In particular, $\Nc=\Nc_k$, see \rf{RN-E}.
\par Let
$$
a\in \Nc\setminus\Nc_n=\usm\limits_{i=n+1}^k (Q_i\,\usm\, Y_i).
$$
Prove that the subarc $\gamma_{\x a}$ of the path $\gamma$ from $\x$ to $a$ intersects $Y_n^{\cl}$, i.e.,
\bel{G-AX1}
\gamma_{\x a}~\ism\, Y_n^{\cl}\ne\emp.
\ee
\par In fact, since $\x\in\Nc_n$ and $a\notin\Nc_n$, we have $\gamma_{\x a}\cap\partial\Nc_n\ne\emp$. On the other hand, the subarc $\gamma_{\x a}\subset\gamma\subset\Nc$. Hence
$$
\gamma_{\x a}~\ism\, \partial\Nc_n\subset\Nc~\ism\, \partial\Nc_n.
$$
Using Proposition \reff{3-SQ} and Lemma \reff{H-QP}, we conclude that $\Nc\cap\partial\Nc_n\subset Y_n^{\cl}$.
Hence,
$$
\gamma_{\x a}\cap\partial\Nc_n\subset Y_n^{\cl}
$$
proving \rf{G-AX1}.
\par In the same way we prove a similar statement: Let $1\le n\le k-1$ and let
$$
\widetilde{\Nc}_n:=\usm\limits_{i=n+1}^k (Q_i\,\usm\, Y_i).
$$
Then for every $b\in\Nc\setminus \widetilde{\Nc}_n=\cup_{i=1}^{n-1}(Q_i\cup Y_i)$ the subarc $\gamma_{b\y}$ of $\gamma$ from $b$ to $\y$ intersects $Y_n^{\cl}$, i.e.,
\bel{G-AX2}
\gamma_{b\y}~\ism\, Y_n^{\cl}\ne\emp.
\ee
\par Note that, by \rf{G-AX1}, for every $1\le n\le k-1$ we have $\gamma\cap Y_n^{\cl}\ne\emp$.
\par Now let us represent the path $\gamma$ in a parametric form, i.e., as a graph of a continuous mapping $\Gamma:[0,1]\to \Nc$ such that $\Gamma(0)=\x$ and $\Gamma(1)=\y$. Let
\bel{BN-D}
B_n:=\{u\in[0,1]:\Gamma(u)\in Y_n^{\cl}\},~~~1\le n\le k-1.
\ee
Then $B_n$ is a non-empty compact subset of $[0,1]$. Let
$$
v_n:=\min B_n~~~~\text{and}~~~V_n:=\max B_n.
$$
Let
$$
s_n:=\Gamma(V_{n-1}),~~~2\le n\le k,
$$
and let $t_n:=\Gamma(v_{n})$, $1\le n\le k-1$.
Then $s_n\in\gamma\cap Y_{n-1}^{\cl}$ and  $t_n\in\gamma\cap Y_{n}^{\cl}$.
\par Let $\gamma_n$ be the subarc of $\gamma$ from $s_n$ to $t_n$. Prove that
\bel{F-G}
\gamma_n\setminus \{s_n,t_n\}\subset Q_n.
\ee
Clearly, since the squares $\{Q_i\}$ are pairwise disjoint, this inclusion imply properties (2) and (3) of the lemma.
\par First let us show that
$$
(\gamma_n\setminus \{s_n,t_n\})\,\ism\,\Nc_{n-1}=\emp.
$$
\par In fact, suppose there exists a point $b$ such that  $b\in \Nc_{n-1}$ and $b\in\gamma_n\setminus \{s_n,t_n\}$. Then, by \rf{G-AX2}, $\gamma_{b\y}\,\cap\, Y_{n-1}^{\cl}\ne\emp$. Therefore there exists $\bu\in[0,1]$ such that $\Gamma(\bu)\in \gamma_{b\y}\,\cap\, Y_{n-1}^{\cl}$.
\par Recall that $b\in\gamma_n$ and $\gamma_n$ is the subarc of $\gamma$ which joins $s_n=\Gamma(V_{n-1})$ to $t_n=\Gamma(v_n)$. Since $\Gamma(\bu)\in \gamma_{b\y}$, we conclude that $\bu>V_{n-1}$. At the same time $\Gamma(\bu)\in Y_{n-1}^{\cl}$ so that, by \rf{BN-D}, $\bu\in B_{n-1}$. Hence, $\bu\le\max B_{n-1}=V_{n-1}$, a contradiction.
\par In the same way, using \rf{G-AX1}, we show that
$
(\gamma_n\setminus \{s_n,t_n\})\,\ism\,\widetilde{\Nc}_{n}=\emp.
$
\par Hence we conclude that
$$
(\gamma_n\setminus \{s_n,t_n\})\subset \Nc\setminus\,(\Nc_{n-1}\cap\widetilde{\Nc}_{n})\subset Q_n
$$
proving \rf{F-G} and the lemma.\bx\bigskip
\par The next lemma describes Sobolev extension properties of the sets from the families $\Gc:=\{G_i:1\le i\le k-2\}$ and $\Hc:=\{H_i:1\le i\le k-2\}$. See \rf{GI-FB}  and \rf{HI-R2}.
\begin{lemma}\lbl{M-PE} Let $m\ge 1$, $2<p<\infty$, and let $\Omega$ be a domain satisfying condition \rf{S-EXT}. Then each set $G_i$ and $H_i$, $1\le i\le k-2$, is a Sobolev extension domain. Furthermore,
\bel{E-GH}
e(\LM(G_i))\le C(m,p)\,\CE~~~\text{and}~~~
e(\LM(H_i))\le C(m,p)\,\CE,~~~1\le i\le k-2.
\ee
Here $\CE$ is the parameter from condition \rf{S-EXT}.
\end{lemma}
\par {\it Proof.} Let us show that for every $\alpha\in(0,1)$ the sets $G_i$ and $H_i$ are $\alpha$-subhyperbolic domains. See Definition \reff{ASHD}. More specifically, we shall prove that for every $a,b\in G_i$ there exists a path $\gamma\subset G_i$ joining $a$ to $b$ such that
\bel{LS-G}
\len_{\alpha,G_i}(\gamma)\le C(\alpha)\,\eta_W\,\|a-b\|^\alpha.
\ee
Here $\eta_W$ is the constant from Corollary \reff{ARC-WP}. We also show that the set $H_i$ has the same property.
\par Note that, given $a,b\in G_i$, by representation \rf{ST-RG}, it suffices to consider the following cases.\medskip
\par {\it The first case: $a,b\in A_i^{\cl}\cap G_i$ or $a,b\in A_{i+1}^{\cl}\cap G_i$.}
\par In this case, given $a,b\in A_i^{\cl}\cap G_i$, by part (a) and part (b) of Lemma \reff{EXT-Q}, there exists a path $\gamma$ which joins $a$ to $b$ in $G_i$ such that
$$
\len_{\alpha,A_i}(\gamma)\le
\tfrac{12}{\alpha} \|a-b\|^\alpha\,.
$$
Since $A_i\subset G_i$, we have
$\len_{\alpha,G_i}(\gamma)\le\len_{\alpha,A_1}(\gamma)$ proving \rf{LS-G} with $C=12/\alpha$.
\par In the same way we treat the case where $a,b\in A_{i+1}^{\cl}\cap G_i$.
\medskip
\par {\it The second case: $\#(Q_i^{\cl}\cap Q_{i+1}^{\cl})=1$, $a\in \hS_i$ and $b\in Q_{i+2}\cup T_{i+1}$. See \rf{T-I}.}
\par Let $\{c\}=Q_i^{\cl}\cap Q_{i+1}^{\cl}$. By Lemma  \reff{Q-DA}, there exists a path $\gamma_{ac}$ connecting $a$ to $c$ in $\hS_i$ such that $\len_{\alpha,\hS_i}(\gamma_{ac})\le \tfrac{3}{\alpha} \|a-c\|^\alpha$. Since $\hS_i\subset G_i$, we have
$\len_{\alpha,G_i}(\gamma_{ac})\le \tfrac{3}{\alpha} \|a-c\|^\alpha$.
\par Note that $c\in A_{i+1}^{\cl}\cap G_i$. As we have proved in the preceding case, there exists a path $\gamma_{cb}$ joining $c$ to $b$ in in $G_i$ such that
$\len_{\alpha,G_i}(\gamma_{cb})\le \tfrac{12}{\alpha} \|b-c\|^\alpha$.
\par Let $\gamma=\gamma_{ac}\cup \gamma_{cb}$. Then
$$
\len_{\alpha,G_i}(\gamma)=\len_{\alpha,G_i}(\gamma_{ac})+
\len_{\alpha,G_i}(\gamma_{cb})\le \tfrac{12}{\alpha}
(\|a-c\|^\alpha+\|c-b\|^\alpha).
$$
By Lemma \reff{H-IJ}, $(2\hS_i)\cap(S_{i+2}^{\cl}\cup\hS_{i+1}^{\cl})=\emp$. Since $Q_{i+2}\subset S_{i+2}$ and $Q_{i+2}\cup T_{i+1}\subset Q_{i+2}^{\cl}\cup \hS_{i+1}^{\cl}$, see \rf{T-I}, we conclude that $b\notin 2\hS_i$.
\par Since $a,c\in \hS_i$, we obtain that $\|a-c\|\le \|a-b\|$. Hence $\|b-c\|\le 2\|a-b\|$ proving that $\len_{\alpha,G_i}(\gamma)\le \tfrac{9}{\alpha} \|a-b\|^\alpha$. Thus in the case under consideration \rf{LS-G} holds.
\par In the same way we treat the case where
$\#(Q_{i+1}^{\cl}\cap Q_{i+2}^{\cl})=1$, $a\in Q_{i}\cup T_{i}$ and $b\in \hS_{i+1}$.
\par It remains to consider\smallskip
\par {\it The third case:
$a\in Q_{i}$, $b\in Q_{i+2}$.}
\par Since the point $a\in Q_i\subset S_i$ and $b\in Q_{i+2}\subset S_{i+2}$, by Lemma \reff{NP-C1}, there exists a point $z\in Q_{i+1}\cup\hS_i\cup\hS_{i+1}$ such that $\|z-a\|\le 2\eta_W\|a-b\|$. Since $a,z\in A_i\cup\hS_{i+1}$ and $z,b\in A_{i+1}\cup\hS_i$, from the results proven in the previous cases it follows the existence of paths $\gamma_1\subset G_i$ and $\gamma_2\subset G_i$ connecting $a$ to $z$ and $z$ to $b$ respectively such that
$$
\len_{\alpha,G_i}(\gamma_1)\le C(\alpha)\|a-z\|^\alpha
~~~\text{and}~~~
\len_{\alpha,G_i}(\gamma_2)\le C(\alpha)\|z-b\|^\alpha\,.
$$
\par Let $\gamma:=\gamma_1\cup\gamma_2$. Then
$$
\len_{\alpha,G_i}(\gamma)=\len_{\alpha,G_i}(\gamma_1)+
\len_{\alpha,G_i}(\gamma_2)\le C(\alpha)
(\|a-z\|^\alpha+\|z-b\|^\alpha).
$$
Since $\|z-a\|\le 2\eta_W\|a-b\|$, we obtain
\be
\len_{\alpha,G_i}(\gamma)&\le& C(\alpha)
(\|a-z\|^\alpha+(\|b-a\|^\alpha+\|a-z\|^\alpha))\nn\\
&\le&  C(\alpha)(1+4\eta_W)\|a-b\|^\alpha\le 5\,C(\alpha)\eta_W\|a-b\|^\alpha\nn
\ee
proving \rf{LS-G} for {\it all} $a,b\in G_i$.\medskip
\par It remains to apply Theorem \reff{SOB-EXT-RN} to $G_i$ and the first inequality in \rf{E-GH} follows.
\par In the same fashion we prove the Sobolev extension property for each $H_i$, $1\le i\le k-2$. We only notice that the main point in this proof is an analog of the third case whose proof is based on Lemma \reff{NP-C1}. But this lemma holds for every $a\in S_i$ and $b\in S_{i+2}$ as well proving the existence of the required point $z\in Q_{i+2}$ in this case.
\par The proof of the lemma is complete.\bx
\medskip
\par The next theorem presents the main result of this section.
\begin{theorem} \lbl{NP-EXT} Let $p>2$, $m\in\N$, and let $\Omega$ be a simply connected bounded domain in $\RT$. Suppose that $\Omega$ is a Sobolev $\LM$-extension domain satisfying the hypothesis of Theorem \reff{MAIN-NEC}. Let $\x,\y\in\Omega$ and let $\Nc=\NPT^{(\x,\y)}$ be a ``Narrow Path'' joining $\x$ to $\y$ in $\Omega$.
\par Then every function $f\in\LM(\Nc)$ extends to a function $F\in\LMPO$ such that
\bel{RF-4}
\|F\|_{\LMPO}\le C(m,p)\,\CE^2\|f\|_{\LM(\Nc)}
\ee
\end{theorem}
\par {\it Proof.}  We prove the theorem in two steps.
\par {\it The first step.} At this step we extend $f$ from ``The Narrow Path'' $\Nc$ to a wider domain $\NW\subset \Wc$. Let
$$
I_{odd}:=\{i: 1\le i\le k-2,~i~~\text{is an odd number}\}\,.
$$
For every $i\in I_{odd}$ we put 
\bel{DGI-1}
\tG_i:=
\left(\,Q_i^{\cl}\,\bigcup \, S_{i+1}^{\cl}\,\bigcup \,Q_{i+2}^{\cl}\bigcup \,\hS_i\,\bigcup \,\hS_{i+1}\right)^{\circ}\,.
\ee
\par Let
\bel{NW-GW}
\NW:=\Nc\bigcup \,
\left\{\bigcup_{i\in I_{odd}}\,\tG_i\right\}.
\ee
\par Comparing this definition with representation \rf{R-NP} we conclude that
\bel{NW-G1}
\NW=\bigcup_{i\in I_{odd}}\,\tG_i ~~~\text{whenever}~~k ~~\text{is odd,}
\ee
and
$$
\NW=\left\{\bigcup_{i\in I_{odd}}\,\tG_i\right\}\bigcup Y_{k-1}\bigcup S_k  ~~~\text{if}~~k ~~\text{is even}.
$$
\par Since $Q_i\subset S_i$, by \rf{WP-DEF}, $\tG_i\subset\Wc$. Hence
\bel{IM-NW-1}
\Nc\subset\NW\subset \Wc.
\ee
\par By Proposition \reff{NP-REP}, ``The Narrow Path'' $\Nc$ is a connected set. The reader can easily see that each set $\tG_i$ is a connected set as well. Clearly,
$\tG_i\cap\Nc\ne\emp$ (because this intersection contains $Q_i$) so that $\NW$ is a {\it connected set}. Since
$\NW$ is open, this set is a {\it domain} in $\RT$.
\par Let $V:=\Nc$, $U:=\NW$, and $\Gc:=\{\tG_i:i\in I_{odd}\}$. Prove that $U,V$ and $\Gc$ satisfy conditions of Proposition \reff{EXT-DM}.
\par First we notice that covering multiplicity of the family $\Gc$ is bounded by $3$. This directly follows from \rf{SIJ-L}, \rf{SH-JN}, and the fact that the squares $\{\hS_i\}$ are pairwise disjoint. See Lemma \reff{H-IJ}.\smallskip
\par Let us show that the members of the family $\{\tG_i\sm \Nc:i\in I_{odd}\}$ are pairwise disjoint.  Let $i,j\in I_{odd}$, $i\ne j$. Hence $|i-j|>1$. By \rf{DGI-1}, \rf{GI-FB} and \rf{GI-NP},
$\tG_i\sm\Nc\subset S_{i+1}^{\cl}\,\ism\,\Omega$.
But, by part (ii) of Lemma \reff{AD-1-2}, the sets $S_{i+1}^{\cl}\cap\Omega$ and $S_{j+1}^{\cl}\cap\Omega$ are disjoint so that the sets $\tG_i\sm\Nc$ and $\tG_j\sm\Nc$ are disjoint as well.\smallskip
\par Prove that
\bel{IM-3}
\tG_i\,\ism\,\Nc=G_i,~~~i\in I_{odd}.
\ee
Clearly, $G_i\subset\tG_i\cap\Nc$, cf. \rf{GI-FB} and \rf{DGI-1}. Note that if $G_i\cap G_j=\emp$ then, by \rf{AI-RG} and \rf{ST-RG}, $|i-j|>2$.
We also notice that, by \rf{DGI-1} and \rf{GI-FB}, $\tG_i\sm G_i\subset S_{i+1}^{\cl}$. On the other hand, by \rf{GI-FB}, for every $j, 1\le j\le k-2$, we have
\bel{IN-T2}
G_j\subset
(S_j^{\cl}\,\usm \, S_{j+1}^{\cl}\,\usm \,S_{j+2}^{\cl}\,\usm \,\hS_j\,\usm \,\hS_{j+1})
\,\ism \,\Omega.
\ee
\par Since $|i-j|>2$, we have $|(i+1)-j|>1$ so that, by part (ii) of Lemma \reff{AD-1-2},
\bel{IN-T1}
S_{i+1}^{\cl}\,\ism\,S_{n}^{\cl}\,\ism\,\Omega=\emp
~~~\text{for every}~~~n=j,j+1,j+2.
\ee
Also, since $j,j+1\ne i+1$, by \rf{SH-JN},
$$
\hS_j\,\ism\,S_{i+1}^{\cl}=
\hS_{j+1}\,\ism\,S_{i+1}^{\cl}=\emp\,.
$$
Combining this with \rf{IN-T2} and  \rf{IN-T1} we conclude that
$$
S_{i+1}^{\cl}\,\ism\, G_j=\emp~~~\text{provided}~~~G_i\,\ism\, G_j=\emp\,.
$$
\par Since $\tG_i\subset S_{i+1}^{\cl}\cup G_i$, see \rf{DGI-1} and \rf{GI-FB}, we obtain that
$$
\tG_i\,\ism\, G_j=\emp
~~~\text{whenever}~~~G_i\,\ism\, G_j=\emp\,.
$$
This property and representation \rf{GI-NP} imply that the set $\Nc\sm G_i$ and the set $\tG_i\cap \Nc$ are disjoint.
Combining this property with the inclusion $G_i\subset \tG_i\cap \Nc$ we obtain the required equality \rf{IM-3}.\smallskip
\par Finally, we notice that, by Lemma \reff{M-PE}, each set $G_i$ is a Sobolev $\LM$-extension domain satisfying inequality \rf{E-GH}.\medskip
\par Now applying Proposition \reff{EXT-DM} to the sets $V$, $U$, and the family $\Gc$ defined above we conclude that the function $f\in\LM(\Nc)$ can be extended to a function $\tF\in\LM(\NW)$ such that
\bel{FE-1}
\|\tF\|_{\LM(\NW)}\le C(m,p)\, \CE\,\|f\|_{\LM(\Nc)}\,.
\ee
\smallskip
\par {\it The second step.} At this step we extend the function $\tF\in\LM(\NW)$ to a function $\hF\in\LM(\Wc)$ with the norm
\bel{FE-2}
\|\hF\|_{\LM(\Wc)}\le C(m,p)\,\CE\, \|\tF\|_{\LM(\NW)}\,.
\ee
\par We construct the extension $\hF$ following the approach suggested at the fist step. Let
$$
I_{even}:=\{i: 1\le i\le k-2,~i~~\text{is an even number}\}.
$$
For every $i\in I_{even}$ we put 
\bel{HW-DF}
\tH_i:=
\left(\,S_i^{\cl}\,\usm \, \, S_{i+1}^{\cl}\,\usm \,S_{i+2}^{\cl}\,\usm \,\hS_i\,\usm \,\hS_{i+1}\right)^{\circ}\,.
\ee
Let
$$
\Hc:=\bigcup_{i\in I_{even}} \tH_i\,.
$$
\par Let $V:=\NW$, $U:=\Wc$, and $\Gc:=\{\tH_i:i\in I_{even}\}$. Prove that these objects satisfy conditions of Proposition \reff{EXT-DM}.
\par First let us prove that \rf{U-V} holds, i.e.,
$$
\Wc=\NW\,\,\usm\,\Hc\,.
$$
This equality is based on the following representation of $\tH_i$:
$$
\tH_i=
\left(\,S_i^{\cl}\,\usm \, S_{i+1}^{\cl}\,\usm \,
\hS_i\right)^{\circ}\,\usm \,\left(\,S_{i+1}^{\cl}\,\usm \, S_{i+2}^{\cl}\,\usm \, \hS_{i+1}\right)^{\circ}\,.
$$
\par This and representation \rf{R-WP} imply the inclusion $\tH_i\subset\Wc$. Since $\NW\subset\Wc$, see \rf{IM-NW-1}, we obtain that $\Wc\supset\NW\bigcup\Hc$.
\par Prove that
\bel{W-H}
\Wc\subset\NW\,\,\usm\,\,\Hc\,.
\ee
By \rf{RW-E}, for every $z\in \Wc$ there exists $i\in\{1,...,k\}$ such that $z\in S_i\cup T_i$. See \rf{T-I}. If $i=1$, then, by  part (1) of Proposition \reff{3-SQ}, $S_1=Q_1$ so that
$$
S_1\,\usm\, Q_1\subset
\left(\,Q_1^{\cl}\,\usm\,S_{2}^{\cl}\,\usm\,
\hS_1\right)^{\circ}\subset \tG_1.
$$
See \rf{DGI-1}. Combining this inclusion with \rf{NW-GW}, we obtain that $z\in \NW$.
\par Let $i=k$. Then $T_k=\emp$, and, by part (1) of Proposition \reff{3-SQ}, $S_k=Q_k$. Hence
$$
z\in S_k\,\usm\, T_k=S_k=Q_k\subset\Nc\subset\NW.
$$
\par Let $k$ be an {\it odd} number, and let $i=k-1$. Then
$$
S_{k-1}\,\usm\, T_{k-1}\subset
\left(\,S_{k-1}^{\cl}\,\usm\,S_{k}^{\cl}\,\usm\,
\hS_{k-1}\right)^{\circ}=
\left(\,S_{k-1}^{\cl}\,\usm\,Q_{k}^{\cl}\,\usm\,
\hS_{k-1}\right)^{\circ}
$$
so that
$$
S_{k-1}\,\usm\, T_{k-1}\subset \tG_{k-1}\subset\NW\,.
$$
See \rf{DGI-1}. Hence $z\in\NW$.
\par Let $1<i<k-1$ or $i=k-1$ and $k$ is {\it even}. If $i$ is even, then $i\le k-2$ so that $i\in I_{even}$. Furthermore, 
$$
S_{i}\,\usm\, T_{i}\subset
\left(\,S_{i}^{\cl}\,\usm\,S_{i+1}^{\cl}\,\usm\,
\hS_{i}\right)^{\circ}\subset \tH_i\,.
$$
If $i$ is odd, then $i-1\in I_{even}$ and
$$
S_{i}\,\usm\, T_{i}\subset
\left(\,S_{i-1}^{\cl}\,\usm\,S_{i}^{\cl}\,\usm\,
\hS_{i-1}\right)^{\circ}\subset \tH_{i-1}\,.
$$
Thus in each case $z\in\Hc$ proving \rf{W-H}.
\par Note that covering multiplicity of the family $\Gc=\{\tH_i:i\in I_{even}\}$ is bounded by $3$. As in the first case, this directly follows from \rf{SIJ-L}, \rf{SH-JN}, and the fact that the squares $\{\hS_i\}$ are pairwise disjoint. See Lemma \reff{H-IJ}.\smallskip
\par Prove that the members of the family $\{\tH_i\sm \NW:i\in I_{even}\}$ are pairwise disjoint.
Let $i,j\in I_{even}$, $i\ne j$. Hence $|i-j|>1$. By \rf{DGI-1}, \rf{NW-GW} and \rf{HW-DF},
$\tH_i\sm\NW\subset S_{i+1}^{\cl}\,\ism\,\Omega$.
By part (ii) of Lemma \reff{AD-1-2}, the sets $S_{i+1}^{\cl}\cap\Omega$ and $S_{j+1}^{\cl}\cap\Omega$ are disjoint so that the sets $\tH_i\sm\NW$ and $\tH_j\sm\NW$ are disjoint as well.\smallskip
\par Prove that
\bel{IN-A}
\tH_i\,\ism\,\NW= H_i,~~~i\in I_{even}.
\ee
See \rf{HI-DF} and \rf{HI-R2}. Clearly, $H_i\subset\tH_i$, cf. \rf{HI-R2} and \rf{HW-DF}. On the other hand, for each $i\in I_{even}$, by \rf{HI-DF} and \rf{DGI-1}, $H_i=B_i\bigcup C_i\subset\tG_{i-1}\bigcup\tG_{i+1}$. Since $i-1$ and $i+1$ are odd numbers, by definition \rf{NW-GW}, $H_i\subset\NW$. Hence $H_i\subset\tH_i\bigcap\NW$.
\par Prove that $\tH_i\bigcap\NW\subset H_i$. Note that if $\tH_i\bigcap \tG_j=\emp$, then, by \rf{DGI-1}, either $j<i-2$ or $i+4<j$. These properties and part (ii) of Lemma \reff{AD-1-2} imply the following:
\bel{L-HG}
\tH_i\,\ism\,\tG_j=\emp~~~\text{provided}~~~
H_i\,\ism\,\tG_j=\emp\,.
\ee
\par This and representation \rf{NW-G1} show that
\bel{DS-N}
\text{the set}~~~\NW\sm H_i~~~\text{and the set}
~~~\tH_i\,\ism\,\NW ~~~\text{are disjoint}
\ee
whenever $k$ is an odd number. If $k$ is even, then $\NW$ is represented by equality \rf{NW-G1}. In this case  $Y_{k-1}\bigcup S_k\subset S_k^{\cl}\bigcap\Omega$ so that, by \rf{HW-DF}, part (ii) of Lemma \reff{AD-1-2} and \rf{SH-JN}, the following is true:
$$
\text{if}~~~\tH_i\,\ism\,(Y_{k-1}\,\usm\, S_k)\ne\emp~~~\text{then}~~~i=k-2\,.
$$
\par Clearly, $H_{k-2}\supset Y_{k-1}\bigcup S_k$. This inclusion, \rf{L-HG} and representation \rf{NW-G1} show that \rf{DS-N} holds for odd number $k$ as well.
\par Now combining \rf{DS-N} with the inclusion $H_i\subset \tH_i\bigcap\NW$ we obtain \rf{IN-A}.
\par Finally, we notice that, by Lemma \reff{M-PE}, each set $H_i$ is a Sobolev $\LM$-extension domain satisfying inequality \rf{E-GH}.
\par These properties of the sets
$\{\tH_i:i\in I_{even}\}$ enable us to apply Proposition \reff{EXT-DM} to the sets $V,U$ and the family $\Gc$ defined at this step. By this proposition, the function $\tF\in\LM(\NW)$ can be extended to a function $\hF\in\LM(\Wc)$ satisfying inequality \rf{FE-2}.
\par Finally we apply Theorem \reff{WP-OM-EXT} to the function $\hF$. By this theorem the function $\hF$ can be extended to a function $F\in\LMPO$ satisfying the following inequality:
$$
\|F\|_{\LMPO}\le C(m,p)\|\hF\|_{\LM(\Wc)}
$$
\par Combining this inequality with inequalities \rf{FE-1} and \rf{FE-2} we obtain the required inequality\rf{RF-4}.
\par The proof of Theorem \reff{NP-EXT} is complete.\bx
\SECT{6. The ``\FG growing'' function}{6}
\setcounter{equation}{0}
\addtocontents{toc}{6. The ``\FG growing'' function.\hfill \thepage\par\VST}
\indent\par Let $\Omega$ be a simply connected bounded domain satisfying the assumption \rf{S-EXT}. In this section, given $\x,\y\in\Omega$ we construct the ``\FG growing'' function $$F_{m}=F_m(z:\x,\y)\in\LM(\Omega)$$ satisfying conditions \rf{T-A}, \rf{T-0} and \rf{T-1}.
For some technical reason it will be more convenient for us to work with a function $H_m=H_m(z:\x,\y)$ which we introduce below than with the function $F_m$. The function $H_m$ is defined by
$$
H_m(z:\x,\y):=\left(\smed_{|\beta|=m-1}
|D^\beta F_m(\y)|\right)^{\frac{1}{p-1}}
\,\cdot\,F_m(z:\x,\y)\,.
$$
Clearly,
\bel{F-HM}
F_m(z:\x,\y):=\left(\smed_{|\beta|=m-1}
|D^\beta H_m(\y)|\right)^{-\frac{1}{p}}
\,\cdot\,H_m(z:\x,\y)\,.
\ee
\par We put this expression in \rf{T-A}, \rf{T-0} and \rf{T-1}, and obtain the following conditions for the function $H_m$:
\bel{T-AH}
D^\beta H_{m}(\x)=0~~~\text{for every multiindex}~~~\beta
~~\text{with}~~|\beta|=m-1,
\ee
\bel{T-0H}
\|H_{m}\|_{\LM(\Omega)}^p\le C_1(m,p,\CE)\,\smed_{|\beta|=m-1}|D^\beta H_{m}(\y)|
\ee
and
\bel{T-1H}
d_{\alpha,\Omega}(\x,\y) \le\, C_2(m,p,\CE)\smed_{|\beta|=m-1}
|D^\beta H_{m}(\y)|.
\ee
Recall that $\alpha=\frac{p-2}{p-1}$ and $\CE$ is the constant from the hypothesis of Theorem \reff{MAIN-NEC}. \smallskip
\par We construct $H_m$ following the approach suggested in Section 1. Thus first we construct a function  $h_{m}\in \LM(\Nc)$ such that
\bel{T-AFS}
D^\beta h_{m}(\x)=0,~~~\text{for every multiindex}~~~\beta
~~\text{with}~~|\beta|=m-1,
\ee
\bel{T-0FS}
\|h_{m}\|_{\LM(\Nc)}^p\le C(m,p)\,\smed_{|\beta|=m-1}|D^\beta h_{m}(\y)|
\ee
and
\bel{T-1FS}
d_{\alpha,\Omega}(\x,\y) \le\, C(m,p)\smed_{|\beta|=m-1}
|D^\beta h_{m}(\y)|.
\ee
Recall that $ \Nc:=\NPT^{(\x,\y)}$ is ``The Narrow Path'' joining $\x$ to $\y$ in $\Omega$. See \rf{NP}.\medskip
\par Then using the Sobolev extension properties  of ``The Wide Path'' and ``The Narrow Path'' proven in Theorems \reff{WP-EXT} and \reff{NP-EXT} respectively, we extend $h_{m}$ to a function $H_{m}\in\LMPO$ such that
$$
\|H_{m}\|_{\LMPO}\le C(m,p,\CE)\,\|h_{m}\|_{\LM(\Nc)}.
$$
\par The function $H_{m}$ satisfies \rf{T-AH}, \rf{T-0H} and \rf{T-1H} so that the function $F_m=F_m(\cdot:\x,\y)$ defined by \rf{F-HM} satisfies \rf{T-A}, \rf{T-0} and \rf{T-1} proving that $F_m$ is the required ``\FG growing ''function.
\medskip
\par Thus the objective of this section is to determine a function $h_{m}\in\LM(\Nc)$ satisfying conditions \rf{T-AFS}, \rf{T-0FS} and \rf{T-1FS}.
\par We define the function $h_m$ with the help of a certain weight function $w:\Nc\to[0,\infty)$.
\begin{definition} \lbl{W-NC} {\em For every $i\in\{1,...,k\}$ and every $z\in Q_i$ we put
\bel{W-ONQ}
w(z):=(\diam Q_i)^{\frac{1}{1-p}}\,.
\ee
\par In turn, we put
\bel{W-ZR}
w(z):=0~~~\text{for every}~~~z\in\Nc\sm \usm\limits_{i=1}^{k}\,Q_i\,.
\ee
}
\end{definition}
\par Thus, in view of representation \rf{RN-E}, $w(z)=0$ provided
$$
z\in \usm\,\{(s_i,t_i): ~\#(Q_i^{\cl}\cap Q_{i+1}^{\cl})>1\}
$$
or
$$
z\in \usm\,\{\hS_i\sm(Q_i\ism Q_{i+1}): ~\#(Q_i^{\cl}\cap Q_{i+1}^{\cl})=1\}.\medskip
$$
Recall that $[s_i,t_i]=Q_i^{\cl}\cap Q_{i+1}^{\cl}$, see \rf{ST-DF}.
\begin{definition} \lbl{L-PH} {\em Given $u,v\in\Nc$ we let $\Lc(u,v)$ denote the family of all paths joining $u$ to $v$ in $\Nc$ {\it with edges parallel to the coordinate axes}. For each path $\gamma\in\Lc(u,v)$ we put
$$
\len_{w,j}(\gamma):=\intl_{\gamma}\,w(z)\,|dz_j|,~~~~j=1,2.
$$
We refer to $\len_{w,j}(\gamma)$ as a $w$-length of $\gamma$ in the direction of the $z_j$-axis.
}
\end{definition}
\par Clearly, for every $\gamma\in\Lc(u,v)$
\bel{LN-V}
\len_{w,j}(\gamma)=\intl_{\gamma^{(j)}}\,w(z)\,ds
\ee
where $\gamma^{(j)}$ is the union of all edges of the path $\gamma$ parallel to the $z_j$-axis.
\medskip
\par Definition \reff{L-PH} motivates us to introduce two important pseudometrics on $\RT$.
\begin{definition} \lbl{RH-12} {\em Let $j\in\{1,2\}$. We introduce a pseudometric $\rj:\Nc\times\Nc\to[0,\infty)$ generated by the $w$-length in the direction of the $z_j$-axis as follows:
$$
\rj(u,v):=\inf\, \len_{w,j}(\gamma),~~~~u,v\in\Nc,
$$
where the infimum is taken over all paths $\gamma\in\Lc(u,v)$.}
\end{definition}
\begin{remark} \lbl{C-LP} {\em Note that $\rj$ a symmetric non-negative function on $\Nc\times\Nc$ satisfying the triangle inequality. But, of course, $\rj(u,v)$ may take the value $0$ for distinct points $u,v\in\Nc$. Thus for each $j=1,2$ the function $\rj$ is a {\it pseudometric} on ``The Narrow Path'' $\Nc=\NPT^{(\x,\y)}$.
\par In particular, for every line segment $[a,b]\subset \Nc$ such that $[a,b]\parallel Oz_2$ and every point $s\in\Nc$ we have $\rp(s,a)=\rp(s,b)$. Correspondingly, $\rv(s,a)=\rv(s,b)$ provided $[a,b]$ is an arbitrary line segment in $\Nc$ parallel to $Oz_1$.\rbx}
\end{remark}
\par Let
\bel{PHI-12}
\pj(z):=\rj(z,\x),~~~~z\in\Nc,~j=1,2.
\ee
\begin{lemma} \lbl{R-W} For each $j\in\{1,2\}$ the function $\pj$ is a locally Lipschitz function on $\Nc$ which belongs to $L^1_p(\Nc)$ and satisfies the following inequality:
$$
\|\pj\|^p_{L^1_p(\Nc)}\,\le\, \smed_{i=1}^k\,(\diam Q_i)^\alpha\,.
$$
\end{lemma}
\par {\it Proof.} Let $j=1$ (the same proof holds for $j=2$). Since $\rp$ satisfies the triangle inequality, for every $u,v\in\Nc$ we have
\bel{F-R}
|\pp(u)-\pp(v)|=|\rp(u,\x)-\rp(v,\x)|\le\rp(u,v)
=\inf_{\gamma\in\Lc(u,v)}\len_{w,1}(u,v)\,.
\ee
\par Let
\bel{W-MX}
w_{max}:=\max_{z\in\Nc} w(z)=\max_{1\le i\le k}(\diam Q_i)^{\frac{1}{1-p}}\,.
\ee
Then, by \rf{F-R} and Definition \reff{L-PH},
$$
|\pp(u)-\pp(v)|\le\, w_{max}\,d_{1,\Nc}(u,v).
$$
Recall that $d_{1,\Nc}$ denotes the geodesic metric on $\Nc$, see \rf{DEF-D}.
\par Applying this inequality to an arbitrary square $K\subset\Nc$ and $u,v\in K$ we obtain the following inequality
\bel{D-UV-N}
|\pp(u)-\pp(v)|\le\, w_{max}\|u-v\|.
\ee
\par Thus $\pp\in \Lip_{loc}(\Nc)$ so that every point $z\in\Nc$ has an open neighborhood where the first order distributional partial derivatives of $\pp$ exist. Hence, by Proposition \reff{LOC-W}, $\pp$ has the first order distributional partial derivatives on all of the set $\Nc$.
\par Let us estimate the norm $\|\pp\|_{L^1_p(\Nc)}$.
As we have noted in Remark \reff{C-LP}, the function $\pp(z)=\rp(z,\x)$, $z\in\Nc$, is constant along straight lines parallel to the axis $Oz_2$. Hence,
\bel{D-Z2}
\frac{\partial\pp}{\partial z_2}(z)\equiv 0~~~\text{on}~~~\Nc.
\ee
\par We also notice that, by \rf{W-ZR},
$$
\|\pp\|_{L^1_p(\Nc)}=\|\pp\|_{L^1_p(U)}
$$
where
$$
U:=\,\usm\limits_{i=1}^k\,Q_i\,.
$$
\par By \rf{F-R} and \rf{LN-V}, for every $u,v\in Q_i$ the following inequality
$$
|\pp(u)-\pp(v)|\le\,M_i\,\|u-v\|
$$
holds. Here $M_i:=\max\{w(x): x\in Q_i\}=(\diam Q_i)^{\frac{1}{1-p}}$, see \rf{W-ONQ}. Hence
$$
\left|\frac{\partial\pp}{\partial z_1}(z)\right|\le (\diam Q_i)^{\frac{1}{1-p}}~~~\text{a.e. on}~~Q_i\,.
$$
By this inequality and \rf{D-Z2},
\be
\|\pp\|^p_{L^1_p(\Nc)}&=&
\intl_{\Nc}\,\left|\frac{\partial\pp}{\partial z_1}(z)\right|^p\,dz
\le\,\smed_{i=1}^k\,\intl_{Q_i}\,
\left|\frac{\partial\pp}{\partial z_1}(z)\right|^p\,dz\nn\\
&\le&
\,\smed_{i=1}^k\,(\diam Q_i)^{\frac{1}{1-p}}\,|Q_i|
=
\,\smed_{i=1}^k\,(\diam Q_i)^{\frac{p-2}{p-1}}\nn
\ee
proving the lemma.\bx
\begin{lemma} \lbl{EST-B} The following inequality
$$
\smed_{n=1}^k\,(\diam Q_n)^\alpha\le\, 8\,\{\pp(\y)+\pv(\y)\}
$$
holds.
\end{lemma}
\par {\it Proof.} By \rf{PHI-12} and Definitions \reff{L-PH} and \reff{RH-12}, the statement of the lemma is equivalent to the following fact: {\it Let $\gamma_1,\gamma_2\in \Lc(\x,\y)$, i.e., $\gamma_1,\gamma_2$ are paths with edges parallel to the coordinate axes each connecting $\x$ to $\y$ in $\Nc$. Then}
$$
\smed_{n=1}^k\,(\diam Q_n)^\alpha \leq 8\left\{\intl_{\gamma_1}w(z)|dz_1|
+\intl_{\gamma_2}w(z)|dz_2|\right\}\,.
$$
\par Let us apply Lemma \reff{C-KZ-2} to the paths $\gamma_j~,~j=1,2$. By this lemma, there exist points $s_n^{(j)},t_n^{(j)}\in \gamma_j~,~1\leq n\leq k,$ such that:\medskip
\par (1). $s_1^{(j)}=\x\,,~t_k^{(j)}=\y$,
$$
s_n^{(j)}\in \gamma\cap Y_{n-1}^{\cl}~~\text{for all}~~2\le n\le k, ~~\text{and}~~t_n^{(j)}\in \gamma\cap Y_{n}^{\cl}~~\text{for all}~~ 1\le n\le k-1\,,~j=1,2\,;
$$
\par (2). Let $\gamma_n^{(j)}$ be a subarc of $\gamma$ with the ends in $s_n^{(j)}$ and $t_n^{(j)}$, $1\le n\le k$, $j=1,2$. Then
$$
\gamma_n^{(j)}\subset Q_n^{\cl}~~\text{for all}~~1\leq n\leq k\,;
$$
\par (3). For each $j=1,2$, the sets $\{\gamma_n^{(j)}\setminus \{s_n^{(j)},t_n^{(j)}\}:1\leq n\leq k\}$ are pairwise disjoint.\bigskip
\par Prove that for every $n$, $1\le n\le k-2$, such that $Q_n^{\cl}\cap S_{n+2}^{\cl}=\emp$ the following inequality
\bel{Q-GZ}
(\diam Q_{n+1})^{\alpha}\leq 4\left\{\, \intl_{\gamma_{n+1}^{(1)}}w(z)|dz_1|+
\intl_{\gamma_{n+1}^{(2)}}w(z)|dz_2|\right\}
\ee
holds. In fact, by Lemma \reff{C-KZ}, in this case
\bel{DQY}
\diam Q_{n+1}\leq 4\dist (Y_n,Y_{n+1}).
\ee
\par Note that, by property (1), for each $j\in\{1,2\}$
\bel{IN-J}
\gamma_n^{(j)}\cap Y_{n-1}^{\cl}\ne \emp~~\text{for all}~~2\le n\le k, ~~\text{and}~~
\gamma_n^{(j)}\cap Y_n^{\cl}\ne\emp~~\text{for all}~~1\le n\le k-1.
\ee
We also notice that, by definition \rf{Y-I}, each set $Y_n$, is either a line segment parallel to one of the coordinate axis, or a square. For such sets the following formula
$$
\dist(Y_n,Y_{n+1})=\max \left\{\dist(\PR_1(Y_n),\PR_1(Y_{n+1})),
\dist(\PR_2(Y_n),\PR_2(Y_{n+1}))\right\}
$$
holds. Here $\PR_j(A)$ denotes the orthogonal projection of a set $A$ on the $z_j$-axis, $j=1,2$\,.
\par By this formula and \rf{DQY}, there exists  $j\in\{1,2\}$ such that
$$
\diam Q_{n+1}\leq 4\dist(\PR_j(Y_n),\PR_j(Y_{n+1})).
$$
\par For simplicity, let us suppose that $j=1$ so that
\bel{D-LJ}
\diam Q_{n+1}\leq 4\dist(\PR_1(Y_n),\PR_1(Y_{n+1}))\,.
\ee
By \rf{IN-J},
$$
\gamma_{n+1}^{(1)}\cap Y_n^{\cl}\ne \emp~~~and~~~
\gamma_{n+1}^{(1)}\cap Y_{n+1}^{\cl}\ne \emp.
$$
Since $\gamma_{n+1}^{(1)}$ is {\it continuous curve}, we have
$$
\dist (\PR_1(Y_n),\PR_1 (Y_{n+1}))\le \length (\PR_1(\gamma_{n+1}^{(1)}))
$$
so that, by \rf{D-LJ}, $\diam Q_{n+1}\leq 4\length (\PR_1(\gamma_{n+1}^{(1)}))$. On the other hand,
$$
\length (\PR_1(\gamma_{n+1}^{(1)}))\leq \intl_{\gamma_{n+1}^{(1)}}|dz_1|
$$
so that
$$
\diam Q_{n+1}\leq 4 \intl_{\gamma_{n+1}^{(1)}}|dz_1|.
$$
\par By property (2) of the present lemma, the path $\gamma_{n+1}^{(1)} \subset Q_{n+1}^{\cl}$, and, by Definition \reff{W-NC},
$w(z)=(\diam Q_{n+1})^{\frac{1}{p-1}}$, $z\in Q_{n+1}$.
Hence,
$$
(\diam Q_{n+1})^{\alpha}=\diam Q_{n+1}\cdot \diam Q_{n+1}^{\frac{1}{1-p}}\leq 4\diam Q_{n+1}^{\frac{1}{1-p}}\intl_{\gamma_{n+1}^{(1)}}|dz_1|=
4\intl_{\gamma_{n+1}^{(1)}}w(z)|dz_1|
$$
proving \rf{Q-GZ}.
\par In the same fashion, using inequalities \rf{DL-VX}, we prove \rf{Q-GZ} for $n=0$ and $n=k-1$.\medskip
\par Now let us consider those numbers $n$ , $1\le n< k-2$, for which $Q_n^{\cl}\cap S_{n+2}^{\cl}\ne \emp$. Then, by \rf{DQ-2} and property (5) of Lemma \reff{3-SQ},
$\diam Q_{n+1}\le \diam Q_{n+2}$
and $Q_{n+1}\cap S_{n+3}=\emp$. As we have proved, in this case
$$
(\diam Q_{n+2})^{\alpha}\leq 4\left\{\, \intl_{\gamma_{n+2}^{(1)}}w(z)|dz_1|+
\intl_{\gamma_{n+2}^{(2)}}w(z)|dz_2|\right\}
$$
so that
\bel{Q-EZ1}
(\diam Q_{n+1})^{\alpha}\leq 4\left\{\, \intl_{\gamma_{n+2}^{(1)}}w(z)|dz_1|+
\intl_{\gamma_{n+2}^{(2)}}w(z)|dz_2|\right\}.
\ee
\par It remains to consider the last case where $n=k-2$ and $Q_{k-2}^{\cl}\cap Q_k^{\cl}\ne\emp$. In this case, by \rf{DQ-2}, $\diam Q_{k-1}\leq \diam Q_k$.
\par As we have noted above, for the case $n=k-1$ inequality \rf{Q-GZ} holds. Hence,
\bel{Q-LZ2}
(\diam Q_{k-1})^{\alpha}\leq \diam Q_k^{\alpha}\leq
4\left\{ \intl_{\gamma_k^{(1)}}w(z)|dz_1|+
\intl_{\gamma_k^{(2)}}w(z)|dz_2|\right\}.
\ee
\par Summarizing inequalities \rf{Q-GZ}, \rf{Q-EZ1} and \rf{Q-LZ2}, we obtain the following:
$$
I=\smed_{n=1}^k\,(\diam Q_n)^\alpha \leq 8\,\smed_{n=1}^k\,\left(\,\,\intl_{\gamma_n^{(1)}}w(z)|dz_1|+
\intl_{\gamma_n^{(2)}}w(z)|dz_2|\right).
$$
But, by property (3) of the present lemma, for each $j=1,2$, the sets $\{\gamma_n^{(j)}\setminus \{s_n^{(j)},t_n^{(j)}\}\}$ are pairwise disjoint. Hence,
$$
I\leq 8\smed_{n=1}^k\,\left(\,\,\intl_{\gamma_n^{(1)}}w(z)|dz_1|+
\intl_{\gamma_n^{(2)}}w(z)|dz_2|\right)\leq 8 \left(\,\,\intl_{\gamma^{(1)}}w(z)|dz_1|+
\intl_{\gamma^{(2)}}w(z)|dz_2|\right).
$$
\par The proof of the lemma is complete.\bx
\begin{lemma} \lbl{D-AXY} The following inequality
$$
d_{\alpha,\Omega}(\x,\y)\le ({12}/{\alpha})\,
\smed_{n=1}^k\,(\diam Q_n)^\alpha
$$
holds.
\end{lemma}
\par{\it Proof.} Let $c_n$ be the center of the square $Q_n$, $n=1,...,k$, and let
$$
G_n=Q_n\cup Q_{n+1}\cup Y_n.
$$
We know that $G_n$ is an open subset of $\Nc$. See \rf{R-NP}, \rf{QST-I} and \rf{Q-BG}.
\par By part (i) of Lemma \reff{EXT-Q}, there exists a path $\gamma_n$, $n=1,...,k-1$, connecting $c_n$ to $c_{n+1}$ in $G_n$ such that
$$
\len_{\alpha,G_n}(\gamma_n)\leq \tfrac{6}{\alpha}\,
\|c_n-c_{n+1}\|^\alpha\,.
$$
See \rf{SH-LN}. Since $Q_n$ and $Q_{n+1}$ are touching squares,
$$
\|c_n-c_{n+1}\|=\tfrac{1}{2}(\diam Q_n +\diam Q_{n+1})
$$
In addition, since $G_n\subset \Omega$, we have $\len_{\alpha,\Omega}\leq \len_{\alpha,G_n}$ so that 
$$
\len_{\alpha,\Omega}(\gamma_n)\leq \len_{\alpha,G_n}(\gamma_n)\leq \tfrac{6}{\alpha\,2^{\alpha}}(\diam Q_n +\diam Q_{n+1})^{\alpha}\leq
\tfrac{6}{\alpha}\left\{(\diam Q_n)^{\alpha}+(\diam Q_{n+1})^{\alpha}\right\}.
$$
\par In turn, by Lemma \reff{2-SQ}, there exists a path $\gamma_k$ joining $c_k$ to $\y$ in $Q_k$ such that
$$
\len_{\alpha,Q_k}(\gamma_k)\leq \ \tfrac{6}{\alpha}\|c_k-\y\|^{\alpha}\,.
$$
Since $Q_k\subset \Omega$ and $\y\in Q_k^{\cl}$, we obtain
$$
\len_{\alpha,\Omega}(\gamma_k)\leq \, \tfrac{6}{\alpha 2^{\alpha}}(\diam Q_k)^{\alpha}\leq \tfrac{6}{\alpha}
(\diam Q_k)^{\alpha}\,.
$$
Let
$$
\gamma := \usm\limits_{n=1}^k \gamma_n\, .
$$
Then
\be
\len_{\alpha,\Omega}(\gamma)
=\smed_{n=1}^k \len_{\alpha,\Omega}(\gamma_n)
&\leq& ({6}/{\alpha})\left\{(\diam Q_k)^{\alpha}+\smed_{n=1}^{k-1}
((\diam Q_n)^{\alpha}+(\diam Q_{n+1})^{\alpha})\right\}\nn\\
&\leq& ({12}/{\alpha})\,\smed_{n=1}^k\,(\diam Q_n)^{\alpha}\,.\nn
\ee
\par But, by \rf{DEF-D}, $d_{\alpha,\Omega}(\x,\y)\leq \len_{\alpha,\Omega}(\gamma)$, and the proof of the lemma is complete.\bx\bigskip
\par We are in a position to define the ``\FG growing'' function on ``The Narrow Path'' $\Nc$.
\begin{definition} \lbl{FGR-F} {\em Let $m\ge 1$, $p>2$, and let $\Omega\subset\RT$ be a simply connected bounded domain. Given $\x,\y\in\Omega$ we put
\bel{F-M1}
h_{1}(z):=\pp(z)+\pv(z),~~~z\in\Nc,
\ee
and
\bel{FGF-D}
h_{m}(z):=\intl_{\gamma}\pp(u)(z_1-u_1)^{m-2}\,du_1+
\pv(u)(z_2-u_2)^{m-2}\,du_2\,,~~~~z=(z_1,z_2)\in\Nc,
\ee
whenever $m>1$. Here $\gamma\in\Lc(\x,z)$ is an arbitrary path joining $\x$ to $\y$ in $\Nc$ with edges parallel to the coordinate axes.}
\end{definition}
\par Recall that the functions $\pj$, $j=1,2$, are defined by \rf{PHI-12}.
\begin{remark} {\em As is customary,
$$
\intl_{\gamma}P_{1,z}(u)\,du_1+P_{2,z}(u)\,du_2
$$
where
\bel{PQ-DF}
P_{j,z}(u):=\pj(u)(z_j-u_j)^{m-2},~~~~j=1,2,
\ee
denotes the standard line integral of the vector field $\vec{F}:=(P_{1,z},P_{2,z})$ along the path $\gamma$.\rbx}
\end{remark}
\begin{lemma} \lbl{FGF-1} (i). The function $h_{m}$, $m>1$, is well defined, i.e., its definition does not depend on the choice of the path $\gamma\in\Lc(\x,z)$ in formula \rf{FGF-D};\medskip
\par (ii). Let $n\in\{1,...,m-2\}$ and let $j\in\{1,2\}$. Then for every path $\gamma\in\Lc(\x,z)$ and every $z=(z_1,z_2)\in\Nc$ the following equality
\bel{DR-J}
\frac{\partial^n h_{m}}{\partial z_j^n}(z) =\frac{(m-2)!}{(m-2-n)!}
\intl_{\gamma}\pj(u)(z_j-u_j)^{m-2-n}\,du_j
\ee
holds. Furthermore,
\bel{DRM-O}
\frac{\partial^{m-1}h_{m}}{\partial z_j^{m-1}}(z) =(m-2)!\,\pj(z),~~~~z\in\Nc,
\ee
and for every $\beta_1,\beta_2>0$, $\beta_1+\beta_2\le m-1$ 
\bel{DRZ-Z}
\frac{\partial^{\beta_1+\beta_2}h_{m}}
{\partial z_1^{\beta_1}\partial z_2^{\beta_2}}\equiv 0 ~~~\text{on}~~~\Nc.
\ee
\end{lemma}
\par {\it Proof.} (i) Let us consider the components $P_1:=P_{1,z}$ and $P_2:=P_{2,z}$ of the vector field $\vec{F}:=(P_{1,z},P_{2,z})$ defined by \rf{PQ-DF}. By this definition and Remark \reff{C-LP}, the function $P_1$ is constant on each interval in $\Nc$ parallel to the $z_2$-axis. In turn, the function $P_2$ is constant on each interval in $\Nc$ parallel to the $z_1$-axis. Hence,
$$
\frac{\partial P_1}{\partial u_2} \equiv 0~~~~\text{and}~~~~
\frac{\partial P_2}{\partial u_1} \equiv 0~~~~\text{on}~~\Nc
$$
proving that
$$
\frac{\partial P_1}{\partial u_2}=
\frac{\partial P_2}{\partial u_1}~~~\text{on}~~\Nc\,.
$$
\par By Proposition \reff{NP-REP} and Lemma \reff{NP-FC}, ``The Narrow Path'' $\Nc$ is a {\it simply connected} plane domain with a piecewise smooth boundary. Therefore, by Green's Theorem, the value of the function $h_m$ in formula \rf{FGF-D} does not depend on the choice of the path $\gamma$ in this formula.
\par In the same fashion we prove that the integral in the right hand side of formula \rf{DR-J} does not depend on the choice of the path $\gamma\in\Lc(\x,z)$.\medskip
\par Prove (ii). We begin with the formulae \rf{DR-J} and \rf{DRM-O}. Let us prove these formulae for $j=1$ (in the same way we prove \rf{DR-J} and \rf{DRM-O} for $j=2$).
\par Let
$$
\thw_{m}(z):=\intl_{\gamma}\pp(u)(z_1-u_1)^{m-2}\,du_1,
~~~~z=(z_1,z_2)\in\Nc.
$$
Prove that for every $n\in\{0,...,m-2\}$, every path $\gamma\in\Lc(\x,z)$ and every $z=(z_1,z_2)\in\Nc$ the following equality
\bel{TH-J1}
\frac{\partial^n \thw_{m}}{\partial z_1^n}(z) =\frac{(m-2)!}{(m-2-n)!}
\intl_{\gamma}\pp(u)(z_j-u_j)^{m-2-n}\,du_1
\ee
holds. Clearly, this equality implies \rf{DR-J}, see \rf{FGF-D}.
\par We prove \rf{TH-J1} by induction on $n$. For $n=0$ nothing to prove. Suppose that \rf{TH-J1} holds for given $n$, $0\le n<m-2$, and prove this statement for $n+1$.
\par Let $z_0=(z_1^{(0)},z_2^{(0)})\in\Nc$ and let $h_t=(t,0)$, $t\in\R$. Let $\gamma\in\Lc(\x,z_0)$ and let
$$
\gamma_t:=\gamma\cup[z_0,z_0+h_t].
$$
Then for $t$ small enough we have:
$$
\frac{\partial^{n+1} \thw_m}{\partial z_1^{n+1}}(z)=
A_{n,m}\,\lim_{t\to 0}\frac{1}{t}\left
\{\intl_{\gamma_t}\pp(u)(z_1+t-u_1)^{m-2-n}du_1-
\intl_{\gamma}\pp(u)(z_1-u_1)^{m-2-n}du_1\right\}
$$
where $A_{n,m}:=(m-2)!/(m-2-n)!$\,. Hence,
\bel{DN-1}
\frac{\partial^{n+1} \thw_m}{\partial z_1^{n+1}}(z)=
A_n\,\lim_{t\to 0}\,(I_1(t)+I_2(t))
\ee
where
$$
I_1(t):=\intl_{\gamma}\pp(u)
\frac{(z_1+t-u_1)^{m-2-n}-(z_1-u_1)^{m-2-n}}{t}\,du_1
$$
and
$$
I_2(t):=\frac{1}{t}\intl_{[z,z+h_t]}
\pp(u)(z_1+t-u_1)^{m-2-n}\,du_1\,.
$$
\par Since the function $\pp(z)=\rp(z,\x)$ is continuous and $n<m-2$, the standard limit theorem for the Riemann integral lead us to the following formula:
$$
\frac{\partial^{n+1} \thw_m}{\partial z_1^{n+1}}(z) =\frac{(m-2)!}{(m-3-n)!}
\intl_{\gamma}\pp(u)(z_1-u_1)^{m-3-n}\,du_1\,.
$$
This proves \rf{TH-J1} for $n+1$.\medskip
\par In particular, for $n=m-2$, we have
$$
\frac{\partial^{m-2} \thw_m}{\partial z_1^{m-2}}(z)=
(m-2)!\,\intl_{\gamma}\pp(u)\,du_1
$$
where $\gamma\in \Lc(\x,z)$ is an arbitrary path. Applying formula \rf{DN-1} to this case with $n=m-2$ we obtain:
$$
\frac{\partial^{m-1} \thw_m}{\partial z_1^{m-1}}(z)
=(m-2)!\,
\lim_{t\to 0}\,\,\frac{1}{t}
\intl_{[z,z+h_t]}\pp(u)\,du_1\,.
$$
Since $\pp$ is a continuous function, we have
$$
\frac{\partial^{m-1}\thw_{m}}{\partial z_1^{m-1}}(z) =(m-2)!\,\pp(z),~~~~z\in\Nc.
$$
Clearly, this equality implies \rf{DRM-O} for $j=1$.\medskip
\par The remaining identity \rf{DRZ-Z} directly follows from the fact that, by formula \rf{DR-J}, for every $n\in\{1,...,m-2\}$ the partial derivative $\frac{\partial^{n} h_m}{\partial z_1^{n}}$
is constant on each interval in $\Nc$ parallel to the $z_2$-axis, and $\frac{\partial^{n} h_m}{\partial z_2^{n}}$ is constant on each interval in $\Nc$ parallel to the $z_1$-axis.
\par The proof of the lemma is complete.\bx
\medskip
\par The results obtained in this section lead us to the following
\begin{proposition} \lbl{P-HSA} The function $h_m=h_m(z:\x,\y)$, $z\in \Nc$, defined by formulae \rf{F-M1} and \rf{FGF-D} belongs to $L_p^m(\Nc)$ and satisfies conditions \rf{T-AFS}, \rf{T-0FS} and \rf{T-1FS}.
\end{proposition}
\par {\it Proof.} Clearly, \rf{T-AFS} follows from \rf{PHI-12}, \rf{DRM-O} and \rf{DRZ-Z}. Prove \rf{T-0FS}. By formulae \rf{DR-J}, \rf{DRM-O} and \rf{DRZ-Z}, $h_m\in C^{m-1}(\Nc)$. Furthermore, by Lemma \reff{R-W}, the functions $\pj$, $j=1,2$, are locally Lipschitz on $\Nc$, so that, by \rf{DRM-O} and \rf{DRZ-Z}, the function $h_m$ belongs to the space $C_{loc}^{m-1,1}(\Nc)$ of functions whose classical partial derivatives of order $m-1$ are locally Lipschitz functions on $\Nc$. It is well known that this space coincides with the space $L_{\infty,loc}^m(\Nc)$ so that the function $h_m$ has (locally) the distributional partial derivatives of all orders up to $m$. Hence, by Proposition \reff{LOC-W}, $h_m$ {\it possesses such derivatives on all of the set $\Nc$}.
\par Furthermore, by \rf{DRM-O},
$$
\frac{\partial^m h_m}{\partial z_j^{m}}(z)=(m-2)!\,\,\frac{\partial \pj}{\partial z_j}(z),~~~~z\in\Nc,~j=1,2\,.
$$
Combining this equality with \rf{DRZ-Z}, we obtain:
$$
\|h_m\|_{L_p^m(\Nc)}=(m-2)!\,(\|\pp\|_{L_p^1(\Nc)}+
\|\pv\|_{L_p^1(\Nc)})\,.
$$
Hence, by Lemma \reff{R-W},
$$
\|h_m\|_{L_p^m(\Nc)}^p\leq (2(m-2)!)^p\,\,\smed_{i=1}^k\,(\diam Q_i)^{\alpha}
$$
so that, by Lemma \reff{EST-B},
$$
\|h_m\|_{L_p^m(\Nc)}^p\leq 8(2(m-2)!)^p\,\,\{\pp(\y)+\pv(\y)\}\,.
$$
In turn, by Lemma \reff{FGF-1},
\bel{F-Y2}
\smed_{|\beta |=m-1}\,|(D^{\beta}h_m)(\y)|=(m-2)!\,(\pp(\y)+\pv(\y))
\ee
so that
$$
\|h_m\|_{L_p^m(\Nc)}^p\leq\,\, \frac{8(2(m-2)!)^p}{(m-2)!}\smed_{|\beta |=m-1}\,|(D^{\beta}h_m)(\y)|
$$
proving \rf{T-0FS}.
\par The remaining inequality \rf{T-1FS} directly follows from Lemma \reff{D-AXY}, Lemma \reff{EST-B} and \rf{F-Y2}. The proposition is completely proved. \bx\medskip
\par Let us construct the function $H_m(z)=H_m(z:\x,\y)$ mentioned at the beginning of the section.
\begin{proposition} \lbl{P-HS} There exists a function $H_m=H_m(z:\x,\y)$, $z\in \Omega$, satisfying conditions \rf{T-AH}, \rf{T-0H} and \rf{T-1H} with constants $C_1= C(m,p)\,\CE^{2p}$ and $C_2=C(m,p)$.
\end{proposition}
\par {\it Proof.} By Theorem \reff{NP-EXT}, the function $h_m(z)=h_m(z:\x,\y)$, $z\in\Nc$, extends to a function
$H_m=H_m(z:\x,\y)$, $z\in \Omega$, such that $H_m\in\LMPO$ and
$$
\|H_m\|_{\LMPO}\le C(m,p)\,\CE^2\,\|h_m\|_{\LM(\Nc)}\,.
$$
This inequality and \rf{T-0FS} imply the following:
$$
\|H_m\|_{\LMPO}^p\le C(m,p)\,\CE^{2p}\, \smed_{|\beta |=m-1}\,|D^{\beta}h_m(\y)|\,.
$$
\par Since $H_m|_{\Nc}=h_m$, we obtain 
\bel{DH-M}
\smed_{|\beta |=m-1}\,|D^{\beta}h_m(\y)|=\smed_{|\beta |=m-1}\,|D^{\beta}H_m(\y)|
\ee
proving \rf{T-0H} with $C_1= C(m,p)\,\CE^{2p}$.
\par Since $h_m$ and $H_m$ coincide on $\Nc$, \rf{T-AFS} implies \rf{T-AH} as well. Finally, \rf{DH-M} and \rf{T-1FS} imply inequality \rf{T-1H} with a constant $C_2=C(m,p)$ proving the proposition.\bx
\medskip
\par {\it Proofs of Theorem \reff{Q-V} and Theorem \reff{MAIN-NEC}.} As we have mentioned in Section 1, the first inequality in \rf{C-MC} follows from Theorem \reff{SOB-EXT-RN}.
\par Let us prove the second inequality in \rf{C-MC} which is equivalent to the statement of Theorem \reff{MAIN-NEC}. We use the approach suggested in Section 1 (after formulation of Theorem \reff{MAIN-NEC}).
\par Let $\Omega$ be a domain satisfying the hypothesis of Theorem \reff{MAIN-NEC}. Since $H_m\in\LMPO$, this function extends to a function $H\in\LMT$ such that
$$
\|H\|_{\LMT}\le \CE\, \|H_m\|_{\LMPO}\,.
$$
Hence, by \rf{T-0H},
\bel{H-W}
\|H\|_{\LMT}^p\le C(m,p)\,\CE^p\cdot\,\CE^{2p}\,\cdot\, D
\ee
where
$$
D:=\smed_{|\beta|=m-1}|D^\beta H_{m}(\y)|\,.
$$
\par On the other hand, since $H|_{\Omega}=H_m$, by \rf{T-AH},
$$
D^p=
\left(\smed_{|\beta|=m-1}|D^\beta H_{m}(\y)-
D^\beta H_{m}(\x)|\right)^p
=
\left(\smed_{|\beta|=m-1}|D^\beta H(\y)-D^\beta H(\x)|\right)^p
$$
so that, by the Sobolev-Poincar\'{e} inequality, see \rf{SP-IN}, and by \rf{H-W},
$$
D^p\le C(m,p)\,\|H\|_{\LMT}^p\,\|\x-\y\|^{p-2}
\le C(m,p)\,\CE^{3p}\, D\,\|\x-\y\|^{p-2}\,.
$$
Hence,
$$
D^{p-1}\le C(m,p)\,\CE^{3p}\,\|\x-\y\|^{p-2}\,.
$$
\par Finally, by inequality \rf{T-1FS},
$$
d_{\alpha,\Omega}(\x,\y) \le\, C(m,p)\, D\le C(m,p)\,\CE^{\frac{3p}{p-1}} \, \|\x-\y\|^{\alpha}\,.
$$
\par The proofs of Theorems \reff{Q-V} and \reff{MAIN-NEC} are complete.\bx
\SECT{7. Further results and comments}{7}
\addtocontents{toc}{7. Further results and comments.\hfill \thepage\par\VST}
\indent
\par {\bf 7.1. Equivalent definitions of Sobolev extension domains.}
\addtocontents{toc}{~~~~7.1. Equivalent definitions of Sobolev extension domains. \hfill \thepage\par}
Let $p\in[1,\infty]$, $m\in \N$, and let $\Omega$ be a domain in $\RN$. We recall that $\Omega$ is said to be a Sobolev $\W$-extension domain if there exists a constant $\theta\ge 1$ such that every function $f\in\WMPO$ extends to a function $F\in\WMP$ such that $\|F\|_{\WMP}\le\, \theta\,\|f\|_{\WMPO}$.
\par Note that this definition is equivalent to the isomorphism of the Banach spaces $\WMP$ and $\WMP|_{\Omega}$, i.e., to the equality
\bel{IS-SM}
\WMP|_{\Omega}=\WMPO.
\ee
Here $\WMP|_{\Omega}$ denotes the trace space of all restrictions of $\WMP$-functions to $\Omega$:
$$
\WMP|_{\Omega}:=\{f:\Omega\to\R:~~ \text{there exists}~~ F\in\WMP ~~ \text{such that}~~F|_\Omega=f\}.
$$
$\WMP|_{\Omega}$ is equipped with the standard quotient space norm
$$
\|f\|_{\WMP|_{\Omega}}:=\inf\{\|F\|_{\WMP}: F\in\WMP,~F|_\Omega=f\}.
$$
\par It is well known that the above definition can be slightly weakened. Namely, we may assume that the trace space $\WMP|_{\Omega}$ and the Sobolev space $\WMPO$ coincide {\it as sets}. In other words, we assume that the restriction operator $R_{\Omega}:\WMP\to\WMPO$ is {\it surjective}, i.e., that every function $f\in\WMPO$ extends to a Sobolev function $F\in\WMP$. Then $\Omega$ is a Sobolev $\W$-extension domain.
\par In fact, let
$$
\ker(R_\Omega):=\{F\in\WMP: F|_\Omega=0\}
$$
be the kernel of the restriction operator, and let $T:\WMP/\ker(R_\Omega)\to \WMPO$ be the projection operator which every equivalence class $[F]\in \WMP/\ker(R_\Omega)$, $F\in\WMP$, assigns the function $f=F|_\Omega$. Clearly,
$T$ is a well defined bounded linear injection (whose  operator norm is bounded by $1$). Since each $f\in\WMP$ extends to a function from $\WMP$, the operator $T$ is a bijection so that, by Banach Inverse Mapping Theorem, it has bounded inverse $T^{-1}:\WMPO\to\WMP/\ker(R_\Omega)$. Hence we conclude that isomorphism \rf{IS-SM} holds proving that $\Omega$ is a Sobolev $\W$-extension domain (with $\theta=\|T^{-1}\|$). \smallskip
\par For various equivalent definitions of Sobolev extension domains we refer the reader to \cite{HKT}. Here we notice the following important statement proven in this paper:
\medskip
\par {\it Let $\Omega\subset\RN$ be an arbitrary domain, $1 <p<\infty$ and $m$ a positive integer. Then $\Omega$ is a Sobolev $\W$-extension domain if and only if there exists a bounded linear extension operator $\Ec:\WMPO\to\WMP$.}
\medskip
\par More specifically, in \cite{HKT} it is proven
that {\it if $\Omega$ is a Sobolev $\W$-extension domain then $\Omega$ is a regular set}, i.e., that there exists a constant $\eta\ge 1$ such that for every $x\in \Omega$
and $0<r\le 1$ the following inequality
$$
|B(x,r)|\le\, \eta\,|B(x,r)\cap \Omega|
$$
holds. Here $B(x,r)$ is the Euclidean ball centered at $x$ with radius $r$.
\par Rychkov \cite{Ry} proved the existence of a bounded linear extension operator $$\Ec:\WMP|_E\to\WMP$$ provided $E$ is an arbitrary {\it regular} subset of $\RN$. See also Shvartsman \cite{S4} where a description of the trace space  $\WMP|_E$ in terms of sharp maximal functions is given. For further results related to characterizations of Sobolev spaces on subsets of $\RN$ we refer the reader to \cite{S7,S-L2}.
\par Finally, we notice that the existence of a bounded linear extension operator from the trace space $\WMP|_E$ into $\WMP$ whenever $E\subset\RN$ is an {\it arbitrary} closed set and $p>n$ has been proven in papers \cite{S7} ($m=1$, $p\in(n,\infty)$), \cite{Isr}, \cite{S-L2} ($n=2$, $m=2$, $p>2$), and \cite{FIL} (arbitrary $m,n$, $p>n$).
\bigskip
\par {\bf 7.2. Self-improvement properties of Sobolev extension domains.}
\addtocontents{toc}{~~~~7.2. Self-improvement properties of Sobolev extension domains. \hfill \thepage\par}
We turn to the proof of the ``open ended property'' of planar Sobolev extension domains, see Theorem \reff{OP-END}. We will also present some other results related to self-improvement properties of Sobolev extension domains and subhyperbolic domains. Proofs of these properties rely on the construction of the ``\FG growing'' function suggested in Section 6, and the following improvement of Theorem \reff{SOB-EXT-RN}.
\begin{theorem}\lbl{SIM-A}(Shvartsman \cite{S10}) Let $m\in\N$, $n<p<\infty$, $\alpha=\tfrac{p-n}{p-1}$, and let $\Omega$ be an $\alpha$\,-\,subhyperbolic domain in $\RN$.
There exists $\tp\in (n,p)$ depending only on $n$, $p$, $m$ and $\Omega$, such that the following is true: every function $f\in L^m_{p,loc}(\Omega)$ extends to a function $F\in\LMPW$ such that
\bel{NA-FAL}
\|F\|_{\LMPW}\le C\,\|f\|_{\LMPWO}
\ee
where $C$ is a positive constant depending only on $n$, $p$, $m$ and $\Omega$.
\end{theorem}
\smallskip
\par This result enables us to prove the following stronger version of Theorem \reff{MAIN-NEC}.
\begin{theorem} \lbl{IM-8} Let $2<p<\infty$, $m\in\N$, and let $\alpha=\frac{p-2}{p-1}$. Let $\Omega\subset\RT$ be a bounded simply connected domain. Suppose that $\Omega$ is a Sobolev $\LM$-extension domain.
\par Then $\Omega$ is an $\tal$-subhyperbolic domain where $\tal\in(0,\alpha)$ is a constant depending only on $p$, $m$ and $\Omega$.
\end{theorem}
\par {\it Proof.} Since $\Omega$ is a Sobolev $\LM$-extension domain, by Theorem \reff{Q-V}, $\Omega$ is an $\alpha$-subhyperbolic domain. Furthermore, $s_\alpha(\Omega)\le C(p,m,\Omega)$.
\par Hence, by Theorem \reff{SIM-A}, there exists a constant $\tp=\tp(p,m,\Omega)$ such that $\tp\in (2,p)$
and every function $f\in \LMPO$ extends to a function $F\in\LMPW$ satisfying inequality \rf{NA-FAL}.
\par Let $\tal=(\tp-2)/(\tp-1)$. Then $0<\tal<\alpha$. Prove that $\Omega$ is an $\tal$-subhyperbolic domain, i.e., that for every $\x,\y\in\Omega$ the following inequality
\bel{D-AW}
d_{\tal,\Omega}(\x,\y)\le \,C\,\|\x-\y\|^{\tal}
\ee
holds. Here $C=C(p,m,\Omega)$.
\par We prove this property by constructing corresponding ``rapidly growing'' function for the exponent $\tp$. In other words, we repeat the definitions related to the
``rapidly growing'' function for $\x,\y$ replacing in these definitions the exponent $p$ with $\tp$. See Section 6.
\par In particular, we modify Definition \reff{W-NC} by letting
$$
\tw(z):=(\diam Q_i)^{\frac{1}{1-\tp}},~~~z\in Q_i,
$$
and $\tw \equiv 0$ on $\Nc\setminus \cup\{Q_i:i=1,...,k\}$.
\par Then we define a pseudometric $\rho_{\tw,j}$, $j=1,2$, by replacing in $\rf{LN-V}$ and Definition \reff{RH-12} the weight $w$ with the new weight $\tw$.\smallskip
\par At the next step we define functions $\tfi_j$, $j=1,2$, by letting
$$
\tfi_j(z):=\rho_{\tw,j}(z,\x),~~~~z\in\Nc.
$$
Cf., \rf{PHI-12}. Note that repeating the proof of  Lemma \reff{R-W} we obtain an analogue of \rf{D-UV-N} for $\tfi_j$, i.e., an inequality
\bel{UV-W}
|\pp(u)-\pp(v)|\le\, \tw_{max}\|u-v\|.
\ee
Here $u,v$ are two arbitrary points of a square $K\subset\Omega$ and
$$
\tw_{max}:=\max_{z\in\Nc} \tw(z)=\max_{1\le i\le k}(\diam Q_i)^{\frac{1}{1-\tp}}\,.
$$
Cf., \rf{W-MX}.
\par Now we are able to define a function $\tlh_m$ on $\Nc$ by replacing in Definition \reff{FGR-F} the function $\varphi_j$ with $\tfi_j$, $j=1,2.$ Then the following analogues of \rf{DRM-O} and \rf{DRZ-Z} hold:
$$
\frac{\partial^{m-1}\tlh_{m}}{\partial z_j^{m-1}}(z) =(m-2)!\,\tfi_j(z),~~~~z\in\Nc,
$$
and 
$$
\frac{\partial^{\beta_1+\beta_2}\tlh_{m}}
{\partial z_1^{\beta_1}\partial z_2^{\beta_2}}\equiv 0 ~~~\text{on}~~~\Nc,~\beta_1,\beta_2>0, \beta_1+\beta_2\le m-1.
$$
\par Combining these properties of $\tlh_m$ with inequality \rf{UV-W}, we conclude that for every multiindex $\beta, |\beta|=m-1$, and every square $K\subset\Omega$ the partial derivative $D^\beta\tlh_m$ is a Lipschitz function on $K$ with the norm $\|D^\beta\tlh_m\|_{\Lip(K)}\le \tw_{max}$. In other words, on each square $K\subset\Omega$ the function $\tlh_m\in C^{m-1,1}(K)$ so that $\tlh_m\in L^m_\infty(K)$
and $\|\tlh_m\|_{L^m_\infty(K)}\le \tw_{max}$. By this inequality and Proposition \reff{LOC-W},
$\tlh_m\in L^m_\infty(\Omega)$ and $\|\tlh_m\|_{L^m_\infty(\Omega)}\le \tw_{max}$.
\par At the next step we construct an analogue of the function $H_m$, a function $\tH_m$, for which analogues of \rf{T-AH}, \rf{T-0H} and \rf{T-1H} hold. Thus
\bel{TW-AH}
D^\beta \tH_{m}(\x)=0~~~\text{for every multiindex}~~~\beta
~~\text{with}~~|\beta|=m-1,
\ee
\bel{TW-0H}
\|\tH_{m}\|_{L^m_{\tp}(\Omega)}^{\tp}\le \tC_1\,\smed_{|\beta|=m-1}|D^\beta \tH_{m}(\y)|
\ee
and
\bel{TW-1H}
d_{\tal,\Omega}(\x,\y) \le\, \tC_2\smed_{|\beta|=m-1}
|D^\beta \tH_{m}(\y)|.
\ee
Here $\tC_j$, $j=1,2$, are positive constants depending only on $p,m$ and $e(\LMPO)$.
\par Furthermore, we claim that 
\bel{HW-LMI}
\tH_{m}\in L^m_\infty(\Omega)\,.
\ee
\par Prove this property of $\tH_m$. Recall that $\tH_m$ is an extension of $h_m$ from $\Nc$ to $\Omega$.
The corresponding extension operator $T_\Omega:\LM(\Nc)\to\LMPO$, $p>2$, is constructed in the proof of Theorem \reff{NP-EXT}. This construction relies on the extensions of $\LM$-functions from certain family $\Dc$ of domains $D\subset\Omega$ of very special geometrical structure. We mean the domains $G_i$ and $H_i$, see \rf{GI-FB} and \rf{HI-R2}, and the domains $\{G\}$ defined in Lemma \reff{EXT-Q}. See also Lemma \reff{M-PE} and Proposition \reff{LM-P}.
\par Note that each special domain $D\in\Dc$ is a union of at most three squares. Furthermore, we have proved that every special domain is an $\alpha$-subhyperbolic set for all $\alpha\in(0,1]$.
\par Recall that an extension operator $\Ec_D:\LM(D)\to\LM(\RT)$ where $D$ is a subhyperbolic domain has been  constructed in \cite{S10}. This operator is a Whitney-type extension operator, and its definition does not depend on $p$. Various approximation properties of $\Ec_D$ have been studied in \cite{S10}. In particular, one of them, Theorem 3.1 proven in \cite{S10}, and standard estimates for the Whitney extension operators, see, e.g.,
Stein \cite{St}, Ch. 6, provide  the required property
$T_\Omega:L^m_\infty(\Nc)\to L^m_\infty(\Omega)$ proving
\rf{HW-LMI}.\medskip
\par We finish the proof of Theorem  \reff{IM-8} following the scheme of the proof of Theorem \reff{MAIN-NEC}. More specifically, we use an analogue of formula \rf{F-HM} and define a function $\tF_m=\tF_m(z;\x,\y)$, $z\in\Omega$, by 
$$
\tF_m(z:\x,\y):=\left(\smed_{|\beta|=m-1}
|D^\beta \tH_m(\y)|\right)^{-\frac{1}{\tp}}
\,\cdot\,\tH_m(z:\x,\y)\,.
$$
Then properties \rf{TW-AH}, \rf{TW-0H} and \rf{TW-1H} imply corresponding properties of $\tF_m$ which are analogues of \rf{T-A}, \rf{T-0} and \rf{T-1}. Thus
$D^\beta \tF_{m}(\x)=0$ for all $\beta$, $|\beta|=m-1$, the norm $\|\tF_{m}\|_{\LMPO}\le C_1$, and
$$
d_{\tal,\Omega}(\x,\y)^{1-\frac{1}{\tp}} \le\, C_2\,\smed_{|\beta|=m-1}
|D^\beta \tF_{m}(\y)|.
$$
\par Furthermore, by \rf{HW-LMI},
$\tF_{m}\in L^m_\infty(\Omega)$. Since $\Omega$ is bounded, the function $\tF_{m}\in\LMPO$. This enables us to apply Theorem \reff{SIM-A} to $\tF_m$. By this theorem, $\tF_m$ extends to a function $\FF\in L^m_{\tp}(\RT)$ such that
\bel{FF-F}
\|\FF\|_{L^m_{\tp}(\RT)}\le C\|\tF_m\|_{L^m_{\tp}(\Omega)}\le C\,C_1.
\ee
\par Inequality \rf{FF-F} is an analogue of inequality \rf{TF-4}. We follow the scheme suggested after this inequality, and replace in estimates \rf{SP-IN}, \rf{SM-ED} and \rf{DA-OM} the function $\TF$ with $\FF$, $p$ with $\tp$, and $F_m$ with $\tF_m$. As a result, we obtain an analogue of inequality \rf{DA-OM} which  states that
$$
d_{\tal,\Omega}(\x,\y)^{1-\frac{1}{\tp}} \le\, C\,\|\x-\y\|^{1-\frac{2}{\tp}}
$$
proving  \rf{D-AW} and the theorem.\bx
\bigskip
\par {\it Proof of Theorem \reff{OP-END}.} By Theorem \reff{IM-8}, the domain $\Omega$ is an $\tal$-subhyperbolic domain for some $\tal\in(0,\alpha)$ depending only on $p$, $m$ and $\Omega$. Let $\tp:=\tfrac{2-\tal}{1-\tal}$,  so that $\tal=\tfrac{\tp-2}{\tp-1}$. Since $0<\tal<\alpha$, we have $2<\tp<p$.
\par Let $\tp<q<\infty$ and let $\alpha^*=\tfrac{q-2}{q-1}$. Clearly, $0<\tal<\alpha^*<1$.
Since $\Omega$ is an $\tal$-subhyperbolic domain, by a result proven in \cite{BKos} (see there Proposition 2.4)
{\it $\Omega$ is an $\alpha^*$-subhyperbolic domain}. Furthermore, $s_{\alpha^*}(\Omega)\le C(\tal,\alpha^*, s_\alpha(\Omega))$. See \rf{SH-MS}.
\par Finally, using Theorem \reff{SOB-EXT-RN} we conclude that $\Omega$ is an $L^k_q$-extension domain proving the theorem.\bx
\par Now we are able to prove the following property of subhyperbolic domains.
\begin{theorem} Let $\alpha\in(0,1)$  and let $\Omega\subset\RT$ be a bounded simply connected $\alpha$-subhyper\-bolic domain. Then $\Omega$ is a  $\beta$-subhyperbolic domain where $\beta\in(0,\alpha)$ is a constant depending only on $\alpha$ and $\Omega$.
\end{theorem}
\par {\it Proof.} Let $p:=\tfrac{2-\alpha}{1-\alpha}$ so that $\alpha=\tfrac{p-1}{p-2}$. Then, by Theorem \reff{SOB-EXT-RN}, $\Omega$ is an $L^1_p$-extension domain so that, by Theorem \reff{IM-8}, $\Omega$ is an $\tal$-subhyperbolic domain where $\tal\in(0,\alpha)$ is a constant depending only on $p$ and $\Omega$.\bx\medskip
\par For a discussion related to this self-improvement property of subhyperbolic domains we refere the reader to \cite{S10}, p. 2209--2210.

\end{document}